\RequirePackage{tikz-cd}
\tikzset{Rightarrow/.style={double equal sign distance,>={Implies},->},
triple/.style={-,preaction={draw,Rightarrow}},
quadruple/.style={preaction={draw,Rightarrow,shorten >=0pt},shorten >=1pt,-,double,double
distance=0.2pt}}

\documentclass[a4paper, 12pt]{article}

\usepackage{amsmath} 
\usepackage{theorem}
\usepackage{amssymb}
\usepackage{amscd}
\usepackage{diagramb}
\usepackage{enumerate}
\usepackage{graphicx}  
\usepackage{cancel}
\usepackage{upgreek}
\usepackage{pxfonts}  
\usepackage{adjustbox}

\DeclareMathAlphabet{\mathpzc}{OT1}{pzc}{m}{it}

\usepackage{xspace,colortbl}
\usepackage{color}

\definecolor{rojo}{rgb}{1,0,0}

\theoremstyle{plain}
\theoremheaderfont{\normalfont\bfseries}

\newtheorem{thm}{Theorem}[section]

\newtheorem{lma}[thm]{Lemma}
\newtheorem{cor}[thm]{Corollary}
\newtheorem{defn}[thm]{Definition}

\theorembodyfont{\rmfamily}

\newtheorem{rem}[thm]{Remark}
\newtheorem{prop}[thm]{Proposition}

\newcommand{\x}{\mbox{-}}

\newcommand{\Ob}{{\mathcal Ob}}
\newcommand{\Cat}{\operatorname {Cat}}

\newcommand{\Mnd}{{\rm Mnd}}
\newcommand{\Mat}{{\rm Mat}}
\newcommand{\Nat}{{\rm Nat}}
\newcommand{\Set}{{\rm Set}}

\newcommand{\Bicat}{{\rm Bicat}}
\newcommand{\I}{{\mathcal I}}
\newcommand{\Span}{\operatorname {Span}}

\newcommand{\Dd}{{\mathbb D}}
\newcommand{\Ee}{\mathbb{E}}
\newcommand{\Ii}{\mathbb{I}}
\newcommand{\Terminal}{\mathbf{1}}

\newcommand{\comp}{\circ}
\newcommand{\iso}{\cong}
\newcommand{\ot}{\otimes}

\newcommand{\C}{{\mathcal C}}
\newcommand{\Tau}{{\mathcal T}}
\newcommand{\M}{{\mathcal M}}
\newcommand{\D}{{\mathcal D}}

\newcommand{\HH}{{\mathcal H}}

\newcommand{\A}{{\mathcal A}}
\newcommand{\B}{{\mathcal B}}
\newcommand{\E}{{\mathcal E}}
\newcommand{\J}{{\mathcal J}}

\def\S{{\mathcal S}}

\newcommand{\crta}{\overline}

\newcommand{\Fi}{\varphi}
\newcommand{\Id}{\operatorname {Id}}

\newcommand{\gama}{\gamma}
\newcommand{\Epsilon}{\varepsilon}

\def\K{{\mathcal K}}  
\def\V{{\mathcal V}}  

\def\Dd{{\mathbb D}}

\def\u#1{\underline{#1}}

\def\Vv{{\mathbb V}}

\newcommand{\Int}{\operatorname{Int}}
\newcommand{\En}{\operatorname{En}}

\newcommand{\cref}[1]{C.~\ref{c:#1}}

\newcommand{\lelabel}[1]{\label{le:#1}}
\newcommand{\leref}[1]{Lemma~\ref{le:#1}}
\newcommand{\eqlabel}[1]{\label{eq:#1}}
\newcommand{\equref}[1]{(\ref{eq:#1})}

\newcommand{\delabel}[1]{\label{de:#1}}
\newcommand{\deref}[1]{Definition~\ref{de:#1}}
\newcommand{\prlabel}[1]{\label{pr:#1}}
\newcommand{\prref}[1]{Proposition~\ref{pr:#1}}
\newcommand{\colabel}[1]{\label{co:#1}}
\newcommand{\coref}[1]{Corollary~\ref{co:#1}}

\newcommand{\rmlabel}[1]{\label{rm:#1}}
\newcommand{\rmref}[1]{Remark~\ref{rm:#1}}
\newcommand{\selabel}[1]{\label{se:#1}}
\newcommand{\seref}[1]{Section~\ref{se:#1}}
\newcommand{\sslabel}[1]{\label{ss:#1}}
\newcommand{\ssref}[1]{Subsection~\ref{ss:#1}}

\newcommand*{\threefrac}[3]{%
  \begin{array}{@{\,}c@{\,}}%
    #1\\
    \hline
    #2\\
    \hline
    #3%
  \end{array}%
}

\usepackage{booktabs}

\begin{document}

\title{Internalization and enrichment \\ via spans and matrices in a tricategory}

\author{Bojana Femi\'c \vspace{4pt} \\
{\small Mathematical Institute of  \vspace{-2pt}}\\
{\small Serbian Academy of Sciences and Arts } \vspace{-2pt}\\
{\small  Kneza Mihaila 36,} \vspace{-2pt}\\
{\small  11 000 Belgrade, Serbia} \vspace{-2pt}\\
{\small  femicenelsur@gmail.com} \vspace{0,4cm}
\and
Enrico Ghiorzi \vspace{4pt} \\
{\small Humanoid Sensing and Perception  \vspace{-2pt}}\\
{\small Istituto Italiano di Tecnologia } \vspace{-2pt}\\
{\small  Via Morego 30} \vspace{-2pt}\\
{\small  16 163 Genoa, Italy} \vspace{-2pt}\\
{\small  enrico.ghiorzi@iit.it}}

\date{}
\maketitle


\begin{abstract}
We introduce categories $\M$ and $\S$ internal in the tricategory $\Bicat_3$ of bicategories, pseudofunctors, pseudonatural transformations and modifications, for matrices and spans in a 1-strict tricategory $V$. Their horizontal tricategories are the tricategories of matrices and spans in $V$. Both the internal and the enriched constructions are tricategorifications of the corresponding constructions in 1-categories. Following \cite{FGK} 
we introduce monads and their vertical morphisms in categories internal in tricategories. We prove an equivalent condition for when the internal 
categories for matrices $\M$ and spans $\S$ in a 1-strict tricategory $V$ are equivalent, and deduce that in that case their corresponding categories of 
(strict) monads and vertical monad morphisms are equivalent, too. We prove that the latter categories are isomorphic to those of categories enriched and discretely internal in $V$, respectively. As a byproduct of our tricategorical constructions we recover some results  
from \cite{Fem}. 
Truncating to 1-categories we recover results from \cite{CFP} and \cite{Ehr} on the equivalence of enriched and discretely internal 1-categories. 
\end{abstract}



\section{Introduction}

When defining a (small) category, we can equivalently require a set of objects and a set of arrows, or instead require a set of objects and, for every pair of objects, a set of arrows between them.
The first notion is that of small category which generalizes to the notion of internal category, and indeed a small category is an internal category in \(\Set\),
the category of sets and functions.
The second notion resembles that of locally small category (although in this case even the collection of objects is small) which generalizes to the notion of enriched category, and indeed a locally small category is a category enriched in \(\Set\) (with the cartesian monoidal structure).
Clearly, though, the two above definitions of a category are equivalent,
so it is to be expected that the notions of internal and enriched categories must be strictly related as well.

As recalled in \seref{Power},
for a Cartesian closed category \(\V\) with finite limits and small coproducts,
we can construct the bicategory \(\V\x Mat\) of matrices on \(\V\), and the bicategory \( Span_d(\V)\) of discrete spans in \(\V\),
i.e., of spans over discrete objects of \(\V\), where an object is discrete if it is a set-indexed coproduct of the terminal object.
Then there is an adjunction between \(Span_d(\V)\) and \(\V\x Mat\),
which becomes a biequivalence when \(\V\) is extensive.
Observe that internal categories are monads in \(Span(\V)\),
and internal categories with a discrete object of objects are monads in \(Span_d(\V)\).
Moreover, enriched categories are monads in \(\V\x Mat\).
This fact is hinted at in \cite{CFP}, where it is observed that the monad morphisms between monads on the bicategories \(Span_d(\V)\) and \(\V\x Mat\) are not, respectively, functors of internal categories in \(\V\) or functors of enriched \(\V\)-categories.
The authors observe that for that to be the case, one would have to use 2-categorical structures. Though  
they do not pursue this direction, in favor of a 1-categorical approach.

In \cite{FGK} it was observed that many distinct mathematical structures can be considered as monads in appropriate 2-categories, whereas though 
their morphisms are not monad morphisms. This was the motivation for the authors to switch to the {\em double category} instead of a 2-category, 
and then to consider the {\em double category of monads in that double category}, rather than the well-known 2-category of monads in a 2-category. 
Besides, it is clear that oftentimes in bicategories the existence of certain additional structures is assumed but it is only in the corresponding 
pseudo-double category that this additional data is indeed contained (think of the bicategory of algebras and their bimodules). As argued in \cite{Shul, Dou}, these additional data should not be neglected, and 
it is often more convenient to consider the setting of internal categories.

In view of  these ideas our goal in the present article is twofold. 
In the first place, we carry out the construction hinted at in \cite{CFP} to make a 2-categorical proof of the characterization of 
equivalences of bicategories of matrices and spans in a 1-category $\V$, on one hand, and of the categories enriched and discretely internal in $\V$, 
on the other. We do this in \prref{double equiv} and \coref{deduce int-en}, 
using constructions of \cite{FGK}. 
Namely, in \cite[Example 2.1]{FGK} a pseudo-double category $\Span(\V)$ of spans in $\V$ was introduced, whose horizontal bicategory is precisely the bicategory 
$Span(\V)$. Then a double category $\Mnd(\Dd)$ of monads in a double category $\Dd$ is 
introduced in \cite[Definition 2.4]{FGK}, so that when $\Dd=\Span(\V)$, the vertical 1-cells in $\Mnd(\Span(\V))$ are morphisms of internal categories in $\V$. 
Inspired by this we define the pseudo-double category $\Span_d(\V)$ by modifying accordingly $\Span(\V)$, and introduce a pseudo-double category 
$\V\x\Mat$  that allows us to extend the biequivalence of bicategories from \cite{CFP} to an equivalence of pseudo-double categories. Then consequently 
the equivalence of categories discretely internal and enriched in $\V$ is recovered, under the corresponding conditions. 

Our second challenge was to generalize the latter construction to tricategories, whose 1-cells obey strict associativity and unitality laws. We call such tricategories {\em 1-strict}. For a 1-strict tricategory $V$ we first set up the notions of tricategroical limits, via the weighted 3-limits, \ssref{weight}. 
For $V$ having 3-pullbacks, 3-products and 3-coproducts we then construct tricategorical analogue of the pseudo-double categories of matrices and spans 
$\V\x\Mat$ and $\Span_d(\V)$ from above. In the terminology of \cite{Shul} they are $(1\times 2)$-categories, respectively $\M$ and $\S$. This means 
that they are internal (bi)categories in a 1-strict tricategory, 
in our case the tricategory $\Bicat_3$ of bicategories, pseudofunctors, pseudonatural transformations and modifications. 
(Internal categories in 1-strict tricategories we introduced in \cite{Fem}.) The (vertical) bicategory of objects in both cases is $\Cat_2$, the 2-category of categories, while the bicategories of morphisms are suitable ones made for matrices and spans in $V$. We also construct 
a lax and a colax internal functor in $\Bicat_3$ between these two $(1\times 2)$-categories $\M$ and $\S$. 

Following the idea of \cite{FGK}, we introduce monads in a tricategory, and then we define a monad in a $(1\times 2)$-category $\Vv$ as a monad in the horizontal tricategory $\HH(\Vv)$ of $\Vv$. We deduce, analogously to the classical 1-categorical case, that (strict) monads in the 
$(1\times 2)$-categories of matrices $\M$ and spans $\S$ are categories enriched and discretely internal in $V$, in the sense of \cite[Definitions 8.1 and 6.2]{Fem}, respectively. 
We then introduce {\em vertical morphisms} of monads in $\Vv$ and the corresponding category of monads $\Mnd(\Vv)$. We prove that if $(1\times 2)$-categories 
$\Vv_1$ and $\Vv_2$ are equivalent, then their respective categories of strict monads $\Mnd(\Vv_1)$ and $\Mnd(\Vv_2)$ are equivalent, too (\prref{lift to monads}).  
We prove equivalence conditions for the $(1\times 2)$-categories $\M$ and $\S$ in the spirit of \cite{CFP} 
(\coref{mainThm}), and deduce that under those conditions their categories of monads $\Mnd(\M)$ and $\Mnd(\S)$ are equivalent, too 
(\coref{equiv-monads}). Consequently, we obtain the equivalence of categories discretely internal and enriched 
in $V$, under those conditions. On the other hand, we also prove a sufficient condition to have a functor from the category of enriched 
to that of internal categories in $V$ (\prref{I preserves monads}). This recovers \cite[Proposition 8.4]{Fem} which here we obtain as a consequence. 
Truncating to 1-categories our results recover those from \cite[Section 4]{CFP} and \cite[Appendix]{Ehr}. 

We point out that although in our \coref{mainThm} an equivalent condition for the $(1\times 2)$-categories $\M$ and $\S$ to 
be equivalent is stated through a triequivalence trifunctor \(\coprod : V^{\u\D} \to V / (\u\D \bullet 1)\), by the construction 
of $\coprod$ and its (co)domain tricategories it is essentially a 2-dimensional functor between sub-bicategories of the bicategories of morphisms $D_1$ and $C_1$ constituting internal categories $\M$ and $\S$. Thus we indeed relate a tricategory triequivalence with a bicategory biequivalence. 
In this sense it is a proper generalization to a higher dimension of one of the two main results in \cite{CFP} (recalled in our \prref{bieq bicats}), 
by which the bicategories of matrices \(\V\x Mat\) and discrete spans \( Span_d(\V)\) are biequivalent if and only if a 1-dimensional 
functor $\V^I \stackrel{\amalg}{\longrightarrow} \V/(I\bullet 1)$ is an equivalence. 
Mind that although we use the same notation, our pseudofunctors $Int$ and $En$ act between bicategories $D_1$ and $C_1$ (different than 
the bicategories \(\V\x Mat\) and \( Span_d(\V)\)!) as parts of internal functors in $\Bicat_3$.

\medskip

The paper is structured as follows. In the Second Section we set up conventions for notations in tricategories. In the Third one we introduce weighted 3-limits, 
and specialize to 3-pullbacks and 3-(co)products. In Section 4 we present the characterization from \cite{CFP} of the equivalence of 
bicategories of spans and matrices in 1-categories in terms of extensivity, 
we extend it to their respective double categories, and considering the double categories of monads in the latter double categories, 
we deduce that if the former double categories are equivalent, then so are the latter. In Sections 5 and 6 we construct the $(1\times 2)$-categories $\M$ and $\S$ of matrices and spans in $V$, respectively, and in Section 7 we define the (co)lax functors between them, 
defining pseudofunctors between their respective bicategories of morphisms $D_1$ and $C_1$. In the last Section we introduce monads in tricategories and the category of monads and vertical monad morphisms in $(1\times 2)$-categories. We also give an equivalent condition for 
$\M$ and $\S$ to be equivalent, consequently their categories of monads $\Mnd(\M)$ and $\Mnd(\S)$, and categories enriched and discretely internal in $V$ are equivalent under those conditions.

\section{Preliminaries: notational conventions and computing with 2-cells} \selabel{notations}

We assume that the reader is familiar with the definition of a tricategory. For the reference we recommend \cite{GPS} and \cite{Gu}. 
By a {\em trifunctor} we mean a trihomomorphism from \cite[Section 3]{GPS}. 
We introduce notational conventions.


Throughout $V$ will be a {\em 1-strict} tricategory, meaning that its 1-cells obey strict associativity and unitality laws, 
and also {\em 2v-strict}, meaning that the associativity and unitality laws for the vertical composition of 2-cells will hold strictly. 
Composition of 1-cells and consequently horizontal composition of 2- and 3-cells we denote by $\ot$, where $y\ot x$ means that first $x$ is applied, then $y$. 
For the horizontal composition we will also use more intuitive notation: $[x\vert y]=y\ot x$. 
Vertical composition of 2- and 3-cells we denote by $\frac{\alpha}{\beta}$, read from top to bottom. Transversal composition of 3-cells we denote by $\cdot$ and read it from right to left.


\bigskip

We are going to use diagrammatic and formulaic notation. 
When we write 2-cells 
in the form of square diagrams, we will usually read them as in the first diagram below, but sometimes also as in the right one: 
$$
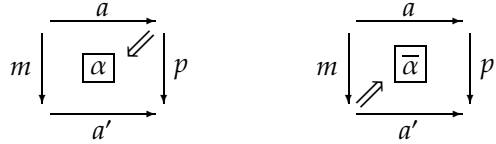
\scalebox{0.86}{\bfig
\putmorphism(0,200)(1,0)[``a]{450}1a
\putmorphism(0,220)(0,-1)[\phantom{Y_2}``m]{400}1l
\putmorphism(450,220)(0,-1)[\phantom{Y_2}``p]{400}1r
\put(150,0){\fbox{$\alpha$}}
\put(310,90){$\Swarrow$}
\putmorphism(0,-150)(1,0)[``a']{450}1b
\efig}
\qquad\qquad
 \scalebox{0.86}{\bfig
\putmorphism(0,200)(1,0)[``a]{450}1a
\putmorphism(0,220)(0,-1)[\phantom{Y_2}``m]{400}1l
\putmorphism(450,220)(0,-1)[\phantom{Y_2}``p]{400}1r
\put(170,0){\fbox{$\crta\alpha$}}
\put(20,-100){$\Nearrow$}
\putmorphism(0,-150)(1,0)[``a']{450}1b
\efig}
$$
This means that $\alpha:pa\Rightarrow a'm$ and $\crta\alpha:a'm\Rightarrow pa$. Usually we will suppress the double arrow labels, and we will write explicitely 
which order of mapping we mean. 

Such written 2-cells can be concatenated both horizontally and vertically. Concatenating them horizontally with 
the two directions of mapping corresponds to the formulaic notations as follows:
\begin{equation} \eqlabel{concatenate-h}
\bfig
\putmorphism(0,200)(1,0)[``a]{450}1a
\putmorphism(460,200)(1,0)[``b]{450}1a
\putmorphism(0,220)(0,-1)[\phantom{Y_2}``m]{400}1l
\putmorphism(450,220)(0,-1)[\phantom{Y_2}``p]{400}1r
\putmorphism(890,220)(0,-1)[\phantom{Y_2}``q]{400}1r
\put(150,0){\fbox{$\alpha$}}
\put(600,0){\fbox{$\beta$}}
\put(310,90){$\Swarrow$}
\put(760,90){$\Swarrow$}
\putmorphism(0,-150)(1,0)[``a']{450}1b
\putmorphism(460,-150)(1,0)[``b']{450}1b
\efig
=\frac{\beta\ot\Id_a}{\Id_{b'}\ot\alpha}
\qquad\qquad
\bfig
\putmorphism(0,200)(1,0)[``a]{450}1a
\putmorphism(460,200)(1,0)[``b]{450}1a
\putmorphism(0,220)(0,-1)[\phantom{Y_2}``m]{400}1l
\putmorphism(450,220)(0,-1)[\phantom{Y_2}``p]{400}1r
\putmorphism(890,220)(0,-1)[\phantom{Y_2}``q]{400}1r
\put(170,0){\fbox{$\crta\alpha$}}
\put(600,0){\fbox{$\crta\beta$}}
\put(20,-100){$\Nearrow$}
\put(470,-100){$\Nearrow$}
\putmorphism(0,-150)(1,0)[``a']{450}1b
\putmorphism(460,-150)(1,0)[``b']{450}1b
\efig
=\frac{\Id_{b'}\ot\crta\alpha}{\crta\beta\ot\Id_a}
\end{equation}
and for vertical concatenation we have:
\begin{equation} \eqlabel{concatenate-v}
\bfig
\putmorphism(0,350)(1,0)[``a]{450}1a
\putmorphism(0,370)(0,-1)[\phantom{Y_2}``m]{400}1l
\putmorphism(450,370)(0,-1)[\phantom{Y_2}``p]{400}1r
\put(150,150){\fbox{$\alpha$}}
\put(150,-210){\fbox{$\gamma$}}
\put(310,240){$\Swarrow$}
\put(310,-140){$\Swarrow$}
\putmorphism(0,10)(0,-1)[\phantom{Y_2}``n]{400}1l
\putmorphism(450,10)(0,-1)[\phantom{Y_2}``r]{400}1r
\putmorphism(0,0)(1,0)[``a']{450}1b
\putmorphism(0,-360)(1,0)[``a'']{450}1b
\efig
=\frac{\Id_r\ot\alpha}{\gamma\ot\Id_m}
\qquad\qquad
\bfig
\putmorphism(0,350)(1,0)[``a]{450}1a
\putmorphism(0,370)(0,-1)[\phantom{Y_2}``m]{400}1l
\putmorphism(450,370)(0,-1)[\phantom{Y_2}``p]{400}1r
\put(170,150){\fbox{$\crta\alpha$}}
\put(20,50){$\Nearrow$}
\putmorphism(0,0)(1,0)[``a']{450}1b
\put(150,-250){\fbox{$\gamma$}}
\put(310,-160){$\Nearrow$}
\putmorphism(0,10)(0,-1)[\phantom{Y_2}``n]{400}1l
\putmorphism(450,10)(0,-1)[\phantom{Y_2}``r]{400}1r
\putmorphism(0,-360)(1,0)[``a'']{450}1b
\efig
=\frac{\crta\gamma\ot\Id_m}{\Id_r\ot\crta\alpha}.
\end{equation}
The two ways of considering 2-cells with respect to the direction of mapping, are clearly inverse to each other. Nevertheless, since 
depending on an occasion we will use both of them, we record these properties. The following simple rule can now easily be deduced, it will 
be useful to us in further computations.

\begin{lma} \lelabel{equiv-2-cells}
Let $V$ be a 1-strict tricategory. 
Given equivalence 2-cells $\alpha, \alpha'$ with their quasi-inverses $\alpha^{-1}, (\alpha')^{-1}$ and consider them as in the diagrams:
$\bfig
\putmorphism(60,200)(1,0)[``a]{380}1a
\putmorphism(60,240)(0,-1)[\phantom{Y_2}``m]{360}1l
\putmorphism(450,240)(0,-1)[\phantom{Y_2}``p]{360}1r
\put(160,40){\fbox{$\alpha$}}
\put(310,90){$\Swarrow$}
\putmorphism(60,-80)(1,0)[``a']{380}1b
\efig$ 
and 
$\bfig
\putmorphism(60,200)(1,0)[``a]{380}1a
\putmorphism(60,240)(0,-1)[\phantom{Y_2}``m]{360}1l
\putmorphism(450,240)(0,-1)[\phantom{Y_2}``p]{360}1r
\put(210,40){\fbox{$\alpha^{-1}$}}
\put(90,-30){$\Nearrow$}
\putmorphism(60,-80)(1,0)[``a']{380}1b
\efig$
and similarly for $\alpha':p'a\Rightarrow a'm'$. Moreover, suppose that we are given 2-cells $\lambda:m\Rightarrow m'$ and $\rho: p\Rightarrow p'$. 
Then to give a 3-cell 
$$
\bfig
  \putmorphism(0,200)(1,0)[``\Id]{450}1a
 \putmorphism(340,200)(1,0)[\phantom{F(B)}` `a]{600}1a

\putmorphism(0,220)(0,-1)[\phantom{Y_2}``m']{400}1l
\putmorphism(450,220)(0,-1)[\phantom{Y_2}``m]{400}1r
\putmorphism(900,220)(0,-1)[\phantom{Y_2}``p]{400}1r
\put(130,0){\fbox{$\lambda$}}
\put(600,0){\fbox{$\alpha$}}
\putmorphism(0,-150)(1,0)[``\Id]{450}1b
 \putmorphism(340,-150)(1,0)[\phantom{F(B)}` `a']{580}1b
\put(310,90){$\Swarrow$}
\put(760,90){$\Swarrow$}
\efig
\quad\stackrel{\Sigma}{\Rrightarrow}\quad
\bfig
  \putmorphism(0,200)(1,0)[``a]{450}1a
 \putmorphism(350,200)(1,0)[\phantom{F(B)}` `\Id]{550}1a

\putmorphism(0,220)(0,-1)[\phantom{Y_2}``m']{400}1l
\putmorphism(450,220)(0,-1)[\phantom{Y_2}``p']{400}1r
\putmorphism(900,220)(0,-1)[\phantom{Y_2}``p]{400}1r
\put(120,0){\fbox{$\alpha'$}}
\put(600,20){\fbox{$\rho$}}
\put(310,90){$\Swarrow$}
\put(760,90){$\Swarrow$}
\putmorphism(0,-150)(1,0)[``a']{450}1b
 \putmorphism(350,-150)(1,0)[\phantom{F(B)}` `\Id]{580}1b
\efig
$$
is equivalent to giving a 3-cell 
$$
\bfig
  \putmorphism(0,200)(1,0)[``\Id]{450}1a
 \putmorphism(340,200)(1,0)[\phantom{F(B)}` `a]{600}1a

\putmorphism(0,220)(0,-1)[\phantom{Y_2}``m]{400}1l
\putmorphism(450,220)(0,-1)[\phantom{Y_2}``m']{400}1l
\putmorphism(900,220)(0,-1)[\phantom{Y_2}``p']{400}1r
\put(150,0){\fbox{$\lambda$}}
\put(600,0){\fbox{$\alpha'^{-1}$}}

\putmorphism(0,-150)(1,0)[``\Id]{450}1b
 \putmorphism(340,-150)(1,0)[\phantom{F(B)}` `a']{580}1b
\put(20,-100){$\Nearrow$}
\put(470,-100){$\Nearrow$}
\efig
\quad\stackrel{\tilde\Sigma}{\Rrightarrow}\quad
\bfig
  \putmorphism(0,200)(1,0)[``a]{450}1a
 \putmorphism(340,200)(1,0)[\phantom{F(B)}` `\Id]{550}1a

\putmorphism(0,220)(0,-1)[\phantom{Y_2}``m]{400}1l
\putmorphism(430,220)(0,-1)[\phantom{Y_2}``p]{400}1r
\putmorphism(860,220)(0,-1)[\phantom{Y_2}``p']{400}1r
\put(150,0){\fbox{$\alpha^{-1}$}}
\put(630,20){\fbox{$\rho$}}
\put(20,-100){$\Nearrow$}
\put(470,-100){$\Nearrow$}

\putmorphism(0,-150)(1,0)[``a']{450}1b
 \putmorphism(340,-150)(1,0)[\phantom{F(B)}` `\Id]{550}1b
\efig. 
$$
In formulas: to give a 3-cell 
$\Sigma: \displaystyle{\frac{[\Id\vert\alpha]}{[\lambda\vert\Id_{a'}]}} \Rrightarrow \displaystyle{\frac{[\Id_a\vert\rho]}{[\alpha'\vert\Id]}}$ 
is equivalent to giving a 3-cell \\
$\tilde\Sigma: \displaystyle{\frac{[\lambda\vert\Id_{a'}]}{[\Id\vert\alpha'^{-1}]}} \Rrightarrow  \displaystyle{\frac{[\alpha^{-1}\vert\Id]}{[\Id_a\vert\rho]}}$. 
\end{lma}

When the direction of mapping of 2-cells written in (square) diagrams is fixed, to shorten notation we will also denote horizontal concatenation 
as in \equref{concatenate-h} by $(\alpha\vert\beta)$, respectively $(\crta\alpha\vert\crta\beta)$.



\section{3-limits}

In this Section we study limits and colimits in 3-categories. In the first Subsection we develope the notion of weighted 3-limits and apply it to deduce 
tricategorical pullbacks, that we will simply call 3-pullbacks. In later Subsections we will introduce 3-(co)products and deduce their properties that 
will be crucial for operating throughout in our proofs.


\subsection{Weighted 3-limits and 3-pullbacks}  \sslabel{weight}

We will need 3-natural transformations among trifunctors, we define them here. 

\begin{defn}
For trifunctors \(F, G \colon \C \to \D\) of 3-categories,
a 3-natural transformation \(\alpha \colon F \Rightarrow G\) is given by
\begin{itemize}
\item For each \(A\) in \(\C\), a 1-cell \(\alpha_A \colon F_0 A \to G_0 A\);
\item For each \(f \colon A \to B\) in \(\C\) an equivalence 2-cell \(\alpha_f\)
    \[
    \begin{tikzcd}
    F_0 A \ar[d, "{\alpha_A}"'] \ar[r, "{F_1 f}"]
    &F_0 B \ar[d, "{\alpha_B}"] \\
    G_0 A \ar[r, "{G_1 f}"'] \ar[ur, Rightarrow, "{\alpha_f}" description]
    &G_0 B
    \end{tikzcd}
    \]
    such that
    \begin{itemize}
    \item For each \(A\) and \(B\) in \(\C\), the \(\alpha_f\) are the components of a 2-natural transformation (equivalence)
    \[
    \alpha_{A, B} 
    \colon (\alpha_A)^* \circ G_{A, B} \Rightarrow (\alpha_B)_* \circ F_{A, B}
    \colon \C(A, B) \to \D(F_0A, G_0B)
    \]
    \item the assignment \(f \mapsto \alpha_f\) is well-behaved with respect to identity and composition.
    \end{itemize}
\end{itemize}
\end{defn}

The above definition implies that, in particular, for each \(\phi \colon f \Rightarrow g \colon A \to B\) in \(\C\) there is an isomorphism 3-cell \(\alpha_\phi\)
\[
\begin{tikzcd}
\alpha_B \circ F_1(f) \ar[d, Rightarrow, "{\Id_{\alpha_B} \otimes F(\phi)}"'] \ar[r, Rightarrow, "{\alpha_f}"]
&G_1(f) \circ \alpha_A \ar[d, Rightarrow, "{G(\phi) \otimes \Id_{\alpha_A}}"] \\
\alpha_B \circ F_1(g)  \ar[r, Rightarrow, "{\alpha_g}"'] \ar[ur, triple, "{\alpha_\phi}" description]
&G_1(g) \circ \alpha_A \\
\end{tikzcd}
\]
such that, for each \(A\), \(B\) and \(f, g \colon A \to B\) in \(\C\), the \(\alpha_\phi\) are the components of a 1-natural transformation (isomorphism)
\[
\alpha_{f, g}
\colon (\alpha_f)^* \circ (G_{A, B})_{f, g} \Rightarrow (\alpha_g)_* \circ (F_{A, B})_{f, g}
\colon \C(A, B)(f, g) \to \D(F_0A, G_0B)(F_1f, G_1g)
\]
Moreover, the assignment \(\phi \mapsto \alpha_{\phi}\) is well-behaved with respect to identity and composition.

\medskip

Let \(\Bicat_3\) be the tricategory of (small) bicategories, pseudofunctors, pseudonatural transformations, 
and modifications between pseudonatural transformations. Observe that it is both 1- and 2v-strict. 

\medskip

Let \(\D = 1 \rightarrow 0 \leftarrow 2\) be the cospan graph, seen as a tricategory, and \(J \colon \D \xrightarrow{!} \Terminal \xrightarrow{\Terminal} \Bicat_3\) the constantly-valued trifunctor on the terminal bicategory.
Let \(\K\) be a tricategory and \(F \colon \D \to \K\) a diagram on \(\K\) with image the cospan \(A \xrightarrow{f} C \xleftarrow{g} B\).
The bicategory \([\D, \K](\Delta_X, F)\) is given by:
\begin{description}
\item[0-cells] 3-natural transformations \(\Delta_X \Rightarrow F\). By definition, that amounts to 1-cells \(p_1 \colon \Delta_X(1) \to F_0(1)\) and \(p_2 \colon \Delta_X(2) \to F_0(2)\), i.e., \(p_1 \colon X \to A\) and \(p_2 \colon X \to B\), such that the square
\[
\begin{tikzcd}
X \ar[d, "{p_1}"'] \ar[r, "{p_2}"]
&B \ar[d, "{g}"] \\
A \ar[r, "{f}"]
&C
\end{tikzcd}
\]
commutes up to an equivalence 2-cell \(\omega \colon f p_1 \Rightarrow g p_2\).
(The third component, \(p_0 \colon \Delta_X(0) \to F_0(0)\), is determined up to equivalence by either \(p_1\) or \(p_2\) via composition with \(f\) or \(g\), respectively.)
To sum up, a 0-cell of \([\D, \K](\Delta_X, F)\) is given by a triple
\[
	(
		p_1 \colon X \to A,
		p_2 \colon X \to B,
		\omega \colon f p_1 \xrightarrow{\equiv} g p_2
	)
\]
\item[1-cells] given 0-cells \((p_1, p_2, \omega)\) and \((q_1, q_2, \sigma)\), a 1-cell \((p_1, p_2, \omega) \to (q_1, q_2, \sigma)\) is a modification of 3-natural transformations, i.e., 2-cells \(\alpha_1 \colon p_1 \Rightarrow q_1\) and \(\alpha_2 \colon p_2 \Rightarrow q_2\) of \(\K\) such that there is an invertible 3-cell \(\Gamma\)
\[
	\begin{tikzcd}
		f p_1 \ar[r, Rightarrow, "{\omega}"] \ar[d, Rightarrow, "{\Id_f \otimes \alpha_1}"']
		& g p_2  \ar[d, Rightarrow, "{\Id_g \otimes \alpha_2}"] \\
		f q_1  \ar[r, Rightarrow, "{\sigma}"] \ar[ur, triple, "{\Gamma}"]
		& g q_2
	\end{tikzcd}
\]
\item[2-cells] given 0-cells \((p_1, p_2, \omega)\) and \((q_1, q_2, \sigma)\), and 1-cells \((\alpha_1, \alpha_2, \Gamma) \colon (p_1, p_2, \omega) \to (q_1, q_2, \sigma)\) and \((\beta_1, \beta_2, \Gamma') \colon (p_1, p_2, \omega) \to (q_1, q_2, \sigma)\),
a 2-cell \((\alpha_1, \alpha_2, \Gamma) \to (\beta_1, \beta_2, \Gamma')\) is a perturbation, i.e., 3-cells \(\Theta_1 \colon \alpha_1 \Rrightarrow \beta_1\) and \(\Theta_2 \colon \alpha_2 \Rrightarrow \beta_2\) in \(\K\) such that
\[
\begin{tikzcd}
	\frac{\Id_f \otimes \alpha_1}{\sigma} \ar[r, triple, "{\Gamma}"] \ar[d, triple, "{\frac{\Id_{\Id_f} \otimes \Theta_1}{\Id_{\sigma}}}"']
	&\frac{\omega}{\Id_g \otimes \alpha_2}
		\ar[d, triple, "{\frac{\Id_{\omega}}{\Id_{\Id_g} \otimes \Theta_2}}"] \\
	\frac{\Id_f \otimes \beta_1}{\sigma} \ar[r, triple, "{\Gamma'}"]
	&\frac{\omega}{\Id_g \otimes \beta_2}
\end{tikzcd}
\]
commutes.
\end{description}

Two 0-cells \((p_1, p_2, \omega)\) and \((q_1, q_2, \theta)\) are equivalent if there exist 1-cells
\begin{align}
	(\alpha_1, \alpha_2, \Gamma) &\colon (p_1, p_2, \omega) \to (q_1, q_2, \theta) \\
	(\beta_1, \beta_2, \Omega) &\colon (q_1, q_2, \theta) \to (p_1, p_2, \omega)
\end{align}
and invertible 2-cells
\begin{align}
	(\Theta_1, \Theta_2) &\colon \left(\beta_1\alpha_1, \beta_2\alpha_2, \frac{(\Id_{g} \otimes \beta_2) \otimes \Gamma}{\Omega \otimes (\Id_f \otimes \alpha_1)}\right) \iso \Id_{(p_1, p_2, \omega)} \\
	(\Omega_1, \Omega_2) &\colon \left(\alpha_1\beta_1, \alpha_2\beta_2, \frac{(\Id_{g} \otimes \alpha_2) \otimes \Omega}{\Gamma \otimes (\Id_f \otimes \beta_1)}\right) \iso \Id_{(q_1, q_2, \sigma)}
\end{align}
meaning that \(p_1 \equiv q_1\) in \(\K(X, A)\) and \(p_2 \equiv q_2\) in \(\K(X, B)\) and the invertible 3-cells \(\Gamma\) and \(\Omega\) exist.

Every 1-cell \(h \colon Y \to X\) induces a pseudofunctor
\[
	[\D, \K](\Delta_{h}, F) \colon [\D, \K](\Delta_{X}, F) \to [\D, \K](\Delta_{Y}, F)
\]
given by precomposition. 
Likewise, every 2-cell \(\phi \colon h \Rightarrow k \colon Y \to X\) induces a pseudonatural transformation
\[
	[\D, \K](\Delta_{\phi}, F) \colon [\D, \K](\Delta_{h}, F) \Rightarrow [\D, \K](\Delta_{k}, F)
\]
defined by
\begin{align}
[\D, \K](\Delta_{\phi}, F)_{(p_1, p_2, \omega)}
	&= (\Id_{p_1} \otimes \phi, \Id_{p_2} \otimes \phi, \upepsilon_{\omega, \phi}) \\
[\D, \K](\Delta_{\phi}, F)_{(\alpha_1, \alpha_2, \Gamma)}
	&= (\upepsilon_{\alpha_1, \phi}, \upepsilon_{\alpha_2, \phi})
\end{align}
for every 0-cell \((p_1, p_2, \omega)\),
and 1-cell \((\alpha_1, \alpha_2, \Gamma) \colon (p_1, p_2, \omega) \to (q_1, q_2, \theta)\),
and where \(\upepsilon\) is the interchange 3-isomorphism; in particular,
\(\upepsilon_{\omega, \phi} \colon \frac{\Id_{f p_1} \otimes \phi}{\omega\otimes\Id_{k}} \Rrightarrow \frac{\omega\otimes\Id_{h}}{\Id_{g p_2}\otimes \phi}\)
and \(\upepsilon_{\alpha_i, \phi} \colon \frac{\alpha_i \otimes \Id_{h}}{\Id_{q_i} \otimes \phi} \Rrightarrow \frac{\Id_{p_i} \otimes \phi}{\alpha_i \otimes \Id_{k}}\).

Finally, every 3-cell \(\Gamma \colon \phi \Rrightarrow \psi \colon h \Rightarrow k \colon Y \to X\) induces a modification
\[
	[\D, \K](\Delta_{\Gamma}, F) \colon [\D, \K](\Delta_{\phi}, F) \Rightarrow [\D, \K](\Delta_{\psi}, F)
\]
defined by
\[
	[\D, \K](\Delta_{\Gamma}, F)(p_1, p_2, \omega)
	= (\Id_{\Id_{p_1}} \otimes \Gamma, \Id_{\Id_{p_2}} \otimes \Gamma).
\]

\medskip

Recall that a biequivalence \(F \colon \B_1 \to \B_2\) of bicategories can be characterized as essentially surjective, fully faithful pseudofunctor (see {\em e.g.} \cite[Definition 2.4.9]{GR}), meaning:
\begin{description}
\item[Essentially surjective] surjective on equivalence classes of objects.
\item[Fully faithful] for each pair of objects \(A\) and \(B\) in \(\B_1\), the component functor
\[
F_{A, B} \colon \B_1(A, B) \to \B_2\big(F_0(A), F_0(B)\big)
\]
is an equivalence of hom-categories, so an essentially surjective, fully faithful functor itself.
\end{description}

\begin{defn}
Let \(\K\) and \(\D\) be tricategories, and \(J \colon \D \to \Bicat_3\) and \(F \colon \D \to \K\) be trifunctors.
A \(J\)-weighted \(3\)-limit of \(F\) is an object \(L\) of \(\K\) equipped with a 3-natural biequivalence
\[
	\epsilon_X \colon \K(X, L) \equiv [\D, \Bicat_3](J, \K(X, F{-}))
\]
where \([\D, \Bicat_3](J, \K(X, F{-}))\) is the bicategory of 3-natural transformations between \(J\) and \(\K(X, F{-})\), modifications between them and their perturbations.
\end{defn}

\begin{prop}
Set \(\D = (1 \rightarrow 0 \leftarrow 2)\) and let \(A \xrightarrow{f} C \xleftarrow{g} B\) be the image of \(F\) in \(\K\).
Also let \(J \colon \D \xrightarrow{!} \Terminal \xrightarrow{\Terminal} \Bicat_3\) be the constantly-valued trifunctor on the terminal bicategory.
A 3-natural transformation \(J \Rightarrow \K(X, F{-})\) amounts to a 3-natural transformation \(\Delta_X \Rightarrow F\).
\end{prop}
\begin{proof}
A 3-natural transformation \(\alpha \colon J \Rightarrow \K(X, F{-})\) is given by a 3-natural family of pseudofunctors \(\alpha_D \colon J_0(D) \to \K(X, F_0(D))\), i.e., \(\alpha_D \colon 1 \to \K(X, F_0(D))\), such that the naturality squares
\[
\begin{tikzcd}
1 \ar[d, "{\Id}"'] \ar[r, "{\alpha_D}"]
&\K(X, F_0(D)) \ar[d, "{\K(X, F_1(f))}"] \\
1 \ar[r, "{\alpha_{D'}}"'] \ar[ur, Rightarrow, "{\alpha_f}"]
&\K(X, F_0(D'))
\end{tikzcd}
\]
for \(f \colon D \to D'\) commute up to a 2-natural equivalence \(\alpha_f\).
That amounts to a family of 1-cells \(\alpha_D \colon X \to F_0(D)\) in \(\K\) and 2-natural equivalence \(\alpha_f \colon F_1(f) \alpha_D \equiv \alpha_{D'}\), i.e., a 3-natural transformation \(\Delta_X \Rightarrow F\).
\end{proof}

Since 3-natural transformations \(J \Rightarrow \K(X, F{-})\) are the same as 3-natural transformations \(\Delta_X \Rightarrow F\), then the 3-pullback of \(g\) along \(f\) is an object \(A \times_{C} B\) of \(\K\) equipped with a 3-natural biequivalence of bicategories
\[
	\epsilon_X \colon \K(X, A \times_{C} B) \equiv [\D, \K](\Delta_X, F)
\]
(i.e., a 3-natural transformation \(\epsilon\) whose components \(\epsilon_X\) are biequivalence of bicategories).
This means that, in particular, there is a biequivalence
\[
	\epsilon_{A \times_{C} B} \colon \K(A \times_{C} B, A \times_{C} B) \equiv [\D, \K](\Delta_{A \times_{C} B}, F)
\]
making the identity on \(A \times_{C} B\) correspond to a 3-natural transformation \(\Delta_{A \times_{C} B} \Rightarrow F\), i.e., to 1-cells
\(p_1 \colon A \times_{C} B \to A\) and
\(p_2 \colon A \times_{C} B \to B\) and to an equivalence 2-cell
\(\omega \colon f p_1 \equiv g p_2\).
So, let \(\epsilon_{A \times_{C} B} (\Id_{A \times_{C} B}) = (p_1, p_2, \omega)\).

For every 0-cell \(X\) in \(\K\), we know that \(\epsilon_X\) is a biequivalence of bicategories, i.e., by the characterization of biequivalences, an essentially surjective, fully faithful pseudofunctor.

%

Consider a 0-cell \(X\) in \(\K\) and a 0-cell
\[
(q_1 \colon X \to A, q_2 \colon X \to B, \sigma \colon f q_1 \equiv g q_2)
\]
in \([\D, \K](\Delta_X, F)\).
Since \(\epsilon_X\) is essentially surjective, there exists a 0-cell
\(u \colon X \to A \times_{C} B\) in \(\K(X, A \times_{C} B)\)
such that \(\epsilon_X(u) \equiv (q_1, q_2, \sigma)\).
By 3-naturality of \(\epsilon_X\) in \(X\),
the naturality square
\[
\begin{tikzcd}
\K(A \times_{C} B, A \times_{C} B)
	\ar[r, "{\epsilon_{A \times_C B}}"]
	\ar[d, "{\K(u, A \times_{C} B)}"']
&{[\D, \K](\Delta_{A \times_{C} B}, F)}
	\ar[d, "{[\D, \K](\Delta_{u}, F)}"] \\
\K(X, A \times_{C} B) \ar[r, "{\epsilon_X}"]
	\ar[ur, Rightarrow, "{\epsilon_u}"]
&{[\D, \K](\Delta_X, F)}
\end{tikzcd}
\]
commutes up to a 2-natural equivalence \(\epsilon_u\).
Thus,
\[
\begin{split}
	\epsilon_X(u)
		&= \epsilon_X(\K(u, A \times_{C} B)(\Id_{A \times_{C} B})) \\
		&\equiv ([\D, \K](\Delta_{u}, F))(\epsilon_{A \times_C B}(\Id_{A \times_{C} B})) \\
		&= ([\D, \K](\Delta_{u}, F))(p_1, p_2, \omega) \\
		&= (p_1 u, p_2 u, \omega \otimes \Id_u)
\end{split}
\]
and thus \((q_1, q_2, \sigma) \equiv (p_1 u, p_2 u, \omega \otimes \Id_{u})\),
i.e., there exist equivalence 2-cells \(\zeta_1 \colon p_1 u \Rightarrow q_1\) and \(\zeta_2 \colon p_2 u \Rightarrow q_2\), and an invertible 3-cell
\(\frac{\Id_{f} \otimes \zeta_1}{\sigma} \Rrightarrow \frac{\omega \otimes \Id_{u}}{\Id_{g} \otimes \zeta_2}\).

Moreover, since \(\epsilon_X\) is fully faithful, for each pair of objects \(u\) and \(v\) in \(\K(X, A \times_C B)\), the component functor
\[
(\epsilon_X)_{u, v} \colon \K(X, A \times_C B)(u, v)
	\to [\D, \K](\Delta_{X}, F)(\epsilon_X(u), \epsilon_X(v))
\]
is an equivalence of hom-categories, i.e., it is an essentially surjective, fully faithful functor.

Since \((\epsilon_X)_{u, v}\) is essentially surjective, then for each object
\[
(
	\alpha_1 \colon p_1 u \Rightarrow p_1 v,
	\alpha_2 \colon p_2 u \Rightarrow p_2 v,
	\upkappa \colon \frac{\Id_{f}\otimes\alpha_1}{\omega\otimes\Id_{v}} \overset{\iso}{\Rrightarrow} \frac{\omega\otimes\Id_{u}}{\Id_{g}\otimes\alpha_2}
)
\]
in \([\D, \K](\Delta_{X}, F)(\epsilon_X(u), \epsilon_X(v))\),
i.e., \([\D, \K](\Delta_{X}, F)\big((p_1u, p_2u, \omega \otimes \Id_{u}), (p_1v, p_2v, \omega \otimes \Id_{v})\big)\),
there exists a 0-cell \(\gamma \colon u \Rightarrow v\)
in \(\K(X, A \times_{C} B)(u, v)\) (i.e., a 2-cell in \(\K\))
such that \((\epsilon_X)_{u, v} (\gamma) \iso (\alpha_1, \alpha_2, \upkappa)\).

By 2-naturality of \(\epsilon_u\) in \(u\), there is a naturality square
\[
\begin{tikzcd}
{\epsilon_X \circ \K(u, A \times_{C} B)}
	\ar[r, Rightarrow, "{\epsilon_u}"]
	\ar[d, Rightarrow, "{\Id_{\epsilon_X} \otimes \K(\gamma, A \times_{C} B)}"']
&{[\D, \K](\Delta_{u}, F) \circ \epsilon_{A \times_C B}}
	\ar[d, Rightarrow, "{[\D, \K](\Delta_{\gamma}, F) \otimes \Id_{\epsilon_{A \times_C B}}}"] \\
{\epsilon_X \circ \K(v, A \times_{C} B)}
	\ar[r, Rightarrow, "{\epsilon_v}"]
	\ar[ur, triple, "{\epsilon_\gamma}" description]
&{[\D, \K](\Delta_{v}, F) \circ \epsilon_{A \times_C B}}
\end{tikzcd}
\]
commuting up to an invertible 3-cell \(\epsilon_\gamma\) in \(\Bicat_3\), i.e., a modification.
(Observe that in the vertices we have composition of pseudofunctors, and in the edges we have horizontal composition of 2-natural transformations.)
Then the component of \(\epsilon_\gamma\) at \(\Id_{A \times_C B}\) is an invertible 2-cell in \([\D, \K](\Delta_{X}, F)\),
\[
\begin{tikzcd}[column sep = large, row sep = large]
{\epsilon_X (\K(u, A \times_{C} B)(\Id_{A \times_C B}))}
	\ar[r, "{\epsilon_u(\Id_{A \times_C B})}"]
	\ar[d, "{(\Id_{\epsilon_X} \otimes \K(\gamma, A \times_{C} B))(\Id_{A \times_C B})}"']
&{[\D, \K](\Delta_{u}, F) (\epsilon_{A \times_C B}(\Id_{A \times_C B}))}
	\ar[d, "{([\D, \K](\Delta_{\gamma}, F) \otimes \Id_{\epsilon_{A \times_C B}})(\Id_{A \times_C B})}"] \\
{\epsilon_X (\K(v, A \times_{C} B)(\Id_{A \times_C B}))}
	\ar[r, "{\epsilon_v(\Id_{A \times_C B})}"']
	\ar[ur, Rightarrow, "{\epsilon_\gamma(\Id_{A \times_C B})}" description]
&{[\D, \K](\Delta_{v}, F)(\epsilon_{A \times_C B}(\Id_{A \times_C B}))}
\end{tikzcd}
\]
i.e., considering that
\[
\begin{split}
([\D, \K](\Delta_{\gamma}, F) \otimes \Id_{\epsilon_{A \times_C B}})(\Id_{A \times_C B})
&= [\D, \K](\Delta_{\gamma}, F) (\epsilon_{A \times_C B}(\Id_{A \times_C B})) \\
&= [\D, \K](\Delta_{\gamma}, F) (p_1, p_2, \omega) \\
&= (\Id_{p_1} \otimes \gamma, \Id_{p_2} \otimes \gamma, \upepsilon_{\omega, \gamma})
\end{split}
\]
where \(\upepsilon_{\omega, \gamma} \colon \frac{\Id_{f p_1} \otimes \gamma}{\omega\otimes\Id_{v}} \Rrightarrow \frac{\omega\otimes\Id_{u}}{\Id_{g p_2}\otimes \gamma}\) is the interchange 3-isomorphism,
and that
\[
(\Id_{\epsilon_X} \otimes \K(\gamma, A \times_{C} B))(\Id_{A \times_C B})
= (\epsilon_X)_{u, v}(\K(\gamma, A \times_{C} B)(\Id_{A \times_C B}))
= (\epsilon_X)_{u, v}(\gamma)
\]
we have an invertible 2-cell
\[
\begin{tikzcd}[column sep = large, row sep = large]
{\epsilon_X (u)}
	\ar[r, "{\epsilon_u(\Id_{A \times_C B})}"]
	\ar[d, "{(\epsilon_X)_{u, v}(\gamma)}"']
&{(p_1 u, p_2 u, \omega \otimes \Id_{u})}
	\ar[d, "{(\Id_{p_1} \otimes \gamma, \Id_{p_2} \otimes \gamma, \upepsilon_{\omega, \gamma})}"] \\
{\epsilon_X (v)}
	\ar[r, "{\epsilon_v(\Id_{A \times_C B})}"']
	\ar[ur, Rightarrow, "{\epsilon_\gamma(\Id_{A \times_C B})}" description]
&{(p_1 v, p_2 v, \omega \otimes \Id_{v})}
\end{tikzcd}
\]
where the 1-cells \(\epsilon_u(\Id_{A \times_C B})\) and \(\epsilon_v(\Id_{A \times_C B})\) are equivalences.
Thus,
\[
(\alpha_1, \alpha_2, \upkappa)
\iso (\epsilon_X)_{u, v} (\gamma) 
\iso (\Id_{p_1} \otimes \gamma, \Id_{p_2} \otimes \gamma,
\upepsilon_{\omega, \gamma} \colon \frac{\Id_{f p_1} \otimes \gamma}{\omega\otimes\Id_{v}} \Rrightarrow \frac{\omega\otimes\Id_{u}}{\Id_{g p_2}\otimes \gamma}
)
\]
up to equivalence 1-cells \(\epsilon_u(\Id_{A \times_C B})\) and \(\epsilon_v(\Id_{A \times_C B})\).
That is, there exist invertible 3-cells
\(\Theta_1 \colon \Id_{p_1} \otimes \gamma \Rrightarrow \alpha_1\) and
\(\Theta_2 \colon \Id_{p_2} \otimes \gamma \Rrightarrow \alpha_2\)
such that
\[
\begin{tikzcd}
	\frac{\Id_{f p_1} \otimes \gamma}{\omega\otimes\Id_{v}}
		\ar[r, triple, "{\upepsilon_{\omega, \gamma}}"]
		\ar[d, triple, "{\frac{\Id_{f} \otimes \Theta_1}{\omega\otimes\Id_{v}}}"']
	&\frac{\omega\otimes\Id_{u}}{\Id_{g p_2}\otimes \gamma}
		\ar[d, triple, "{\frac{\omega\otimes\Id_{u}}{\Id_{g} \otimes \Theta_2}}"] \\ 
	\frac{\Id_{f}\otimes\alpha_1}{\omega\otimes\Id_{v}}
		\ar[r, triple, "{\upkappa}"]
	&\frac{\omega\otimes\Id_{u}}{\Id_{g}\otimes\alpha_2}
\end{tikzcd}
\]
commutes.

Since \((\epsilon_X)_{u, v}\) is fully faithful,
then for each pair of 0-cells \(\gamma, \gamma' \colon u \Rightarrow v\)
in \(\K(X, A \times_{C} B)(u, v)\) (i.e., 2-cells in \(\K\))
and 1-cell
\[
	\Big( (\epsilon_X)_{u, v}(\gamma) \xrightarrow{(\upchi_1, \upchi_2)} (\epsilon_X)_{u, v}(\gamma') \Big) =
	\Big( (\Id_{p_1} \otimes \gamma, \Id_{p_2} \otimes \gamma, \upepsilon_{\omega, \gamma})
	\to (\Id_{p_1} \otimes \gamma', \Id_{p_2} \otimes \gamma', \upepsilon_{\omega, \gamma'}) \Big)
\]
in \([\D, \K](\Delta_{X}, F) \big(\epsilon_X(u), \epsilon_X(v)\big)\),
where
\(\upchi_1 \colon \Id_{p_1} \otimes \gamma \Rrightarrow \Id_{p_1} \otimes \gamma'\) and
\(\upchi_2 \colon \Id_{p_2} \otimes \gamma \Rrightarrow \Id_{p_2} \otimes \gamma'\) such that
\[
\begin{tikzcd}
	\frac{\Id_{f} \otimes (\Id_{p_1} \otimes \gamma)}{\omega\otimes\Id_{v}}
		\ar[r, triple, "{\upepsilon_{\omega, \gamma}}"]
		\ar[d, triple, "{\frac{\Id_{f} \otimes \upchi_1}{\omega\otimes\Id_{v}}}"']
	&\frac{\omega\otimes\Id_{u}}{\Id_{g} \otimes (\Id_{p_2} \otimes \gamma)}
		\ar[d, triple, "{\frac{\omega\otimes\Id_{u}}{\Id_{g} \otimes \upchi_2}}"] \\ 
	\frac{\Id_{f} \otimes (\Id_{p_1} \otimes \gamma')}{\omega\otimes\Id_{v}}
		\ar[d, triple, "{\frac{\upalpha_{\Id_{f}, \Id_{p_1}, \gamma'}}{\omega\otimes\Id_{v}}}"']
	&\frac{\omega\otimes\Id_{u}}{\Id_{g} \otimes (\Id_{p_2} \otimes \gamma')}
		\ar[d, triple, "{\frac{\omega\otimes\Id_{u}}{\upalpha_{\Id_{g}, \Id_{p_2}, \gamma'}}}"] \\
	\frac{(\Id_{f} \otimes \Id_{p_1}) \otimes \gamma'}{\omega\otimes\Id_{v}}
		\ar[r, triple, "{\upepsilon_{\omega, \gamma'}}"]
	&\frac{\omega\otimes\Id_{u}}{(\Id_{g} \otimes \Id_{p_2}) \otimes \gamma'}
\end{tikzcd}
\]
commutes,
there exists a \emph{unique} 1-cell \(\upchi \colon \gamma \Rrightarrow \gamma'\)
in \(\K(X, A \times_C B)(u, v)\) (i.e., a 3-cell in \(\K\))
such that
\(((\epsilon_X)_{u, v})_{\gamma, \gamma'}(\upchi) = (\upchi_1, \upchi_2)\).
We need to compute \(((\epsilon_X)_{u, v})_{\gamma, \gamma'}(\upchi)\).
By naturality of \(\epsilon_{\gamma}\) in \(\gamma\),
for any 1-cell \(\upchi \colon \gamma \Rrightarrow \gamma'\) we have that the naturality square
\[
\begin{tikzcd}
\frac{\Id_{\epsilon_X} \otimes \K(\gamma, A \times_{C} B))}{\epsilon_v}
	\ar[r, triple, "{\epsilon_{\gamma}}"]
	\ar[d, triple, "{\frac{\Id_{\Id_{\epsilon_X}} \otimes \K(\upchi, A \times_{C} B)}{\Id_{\epsilon_v}}}"']
&\frac{\epsilon_u}{[\D, \K](\Delta_{\gamma}, F) \otimes \Id_{\epsilon_{A \times_C B}}}
	\ar[d, triple, "{\frac{\Id_{\epsilon_u}}{[\D, \K](\Delta_{\upchi}, F) \otimes \Id_{\Id_{\epsilon_{A \times_C B}}}}}"] \\
\frac{\Id_{\epsilon_X} \otimes \K(\gamma', A \times_{C} B)}{\epsilon_v}
	\ar[r, triple, "{\epsilon_{\gamma'}}"']
&\frac{\epsilon_u}{[\D, \K](\Delta_{\gamma'}, F) \otimes \Id_{\epsilon_{A \times_C B}}}
\end{tikzcd}
\]
commutes.
Computing the component of the square at \(\Id_{A \times_C B}\), and considering that \((\epsilon_{\gamma})_{\Id_{A \times_C B}}\) and \((\epsilon_{\gamma'})_{\Id_{A \times_C B}}\) are invertible 2-cells and \((\epsilon_{u})_{\Id_{A \times_C B}}\) and \((\epsilon_{v})_{\Id_{A \times_C B}}\) are equivalence 1-cells, and that
\[
(\Id_{\Id_{\epsilon_X}} \otimes \K(\upchi, A \times_{C} B))_{\Id_{A \times_C B}}
= \epsilon_X (\K(\upchi, A \times_{C} B)_{\Id_{A \times_C B}})
= \epsilon_X (\upchi)
\]
and
\[
\begin{split}
([\D, \K](\Delta_{\upchi}, F) \otimes \Id_{\Id_{\epsilon_{A \times_C B}}})_{\Id_{A \times_C B}}
&= [\D, \K](\Delta_{\upchi}, F)({\epsilon_{A \times_C B}}(\Id_{A \times_C B})) \\
&= [\D, \K](\Delta_{\upchi}, F)(p_1, p_2, \omega) \\
&= (\Id_{\Id_{p_1}} \otimes \upchi, \Id_{\Id_{p_2}} \otimes \upchi)
\end{split}
\]
we get that \(\epsilon_X (\upchi) = (\Id_{\Id_{p_1}} \otimes \upchi, \Id_{\Id_{p_2}} \otimes \upchi)\),
up to invertible 2-cells \((\epsilon_{\gamma})_{\Id_{A \times_C B}}\) and \((\epsilon_{\gamma'})_{\Id_{A \times_C B}}\)
and equivalence 1-cells \((\epsilon_{u})_{\Id_{A \times_C B}}\) and \((\epsilon_{v})_{\Id_{A \times_C B}}\).
Thus, \(\upchi\) is such that
\(\upchi_1 = \Id_{\Id_{p_1}} \otimes \upchi\) and
\(\upchi_2 = \Id_{\Id_{p_2}} \otimes \upchi\).

Based on the above discussion, we characterize 3-pullbacks in \deref{3-pb}.


\subsection{3-(co)products} \sslabel{terminal}

In the characterization of 3-limits through biequivalence of bicategories we will also use an equivalent reformulation of the condition 
of essential fullness, in terms of a pair of 
pseudofunctors $F:\B_1\to\B_2:G$ and pseudonatural transformations $\Psi:GF\Rightarrow\Id$ and $\Phi: \Id\Rightarrow FG$ which are themselves 
equivalences. Namely, essential surjectiveness corresponds to an equivalence 2-cell $\Phi_p:p\Rightarrow FG(p)$ in $V$ for every 0-cell $p\in\B_2$, 
which is a 1-cell in $V$, whereas essential fullness corresponds to the fact 
that for every 1-cell in $\B_2$, that is a 2-cell $\theta:p\Rightarrow p'$ in $V$, there is a family of invertible 2-cells in $\B_2$, 
that is invertible 3-cells 
\begin{equation} \eqlabel{Phi cond}
\Phi_\theta: \frac{\theta}{\Phi_{p'}}\Rrightarrow\frac{\Phi_p}{FG(\theta)}
\end{equation} 
in $V$, which moreover satisfy a naturality condition.

\bigskip

A terminal object in a tricategory $V$ is an object 1 such that for any object $X$ in $V$ the bicategory $V(X,1)$ is biequivalent 
to the terminal bicategory. This means that there is a particular 1-cell $!:X\to 1$ such that for every 1-cell $f:X\to 1$ there is an equivalence 2-cell $\psi: !\Rightarrow f$ unique up to a unique isomorphism, every 2-endocell on ! is isomorphic to identity, and 
the only 3-endocell on the identity 2-cell on ! is the identity one. 

\medskip

In particular, for a terminal object 1 and any 1-cell $f:A\to B$ there is an equivalence 2-cell:
\begin{equation} \eqlabel{terminal}
\bfig
 \putmorphism(0,200)(1,0)[A`B ` f]{550}1a
 \putmorphism(-60,210)(1,-1)[\phantom{A}` `!]{350}1l
 \putmorphism(-70,230)(1,-1)[\phantom{A}`1 `]{350}0l
\putmorphism(610,210)(-1,-1)[`\phantom{M\times N}`!]{350}1r
\put(200,70){\fbox{$\kappa_f$}}
\efig
\end{equation}
unique up to a unique isomorphism, which we will call {\em terminal 2-cell}. 
Moreover, due to \equref{Phi cond}, given any 2-cell $\alpha:f\Rightarrow g$ there is an isomorphism 3-cell 
$K:\frac{\alpha}{\kappa_g}\Rrightarrow \kappa_f$. We call them {\em terminal 3-cells}. 

\medskip

Analogous dual properties we have for an initial object $0$ of $V$, with {\em initial 2- and 3-cells}.

\medskip

We recall here the definition of a 3-pullback from \cite{Fem}.

\begin{defn} \delabel{3-pb}
A 3-pullback with respect to a cospan \(M \xrightarrow{f} S \xleftarrow{g} N\) of 1-cells 
in a tricategory $V$ is given by: 
a 0-cell $P$, 1-cells $p_1:P\to M, p_2:P\to N$ and an equivalence 2-cell $\omega:  gp_2\Rightarrow fp_1$ so that 
   \begin{enumerate}
   \item for every 0-cell $T$, 1-cells $q_1: T\to M, q_2:T\to N$ and equivalence 2-cell $\sigma: gq_2\Rightarrow fq_1$ 
	there exist a 1-cell $u:T\to P$, equivalence 2-cells $\zeta_1: p_1u\Rightarrow q_1$ and $\zeta_2: q_2\Rightarrow p_2u$ and 
	an invertible 3-cell 
	$$\Sigma: \threefrac{\Id_{g}\ot\zeta_2}{\omega\ot\Id_u}{\Id_f\ot\zeta_1} \Rrightarrow\sigma$$
	\[
\begin{tikzcd}
T \ar[dr, "{u}"] \ar[ddr, bend right, "{q_1}"', ""{name=q1, above}] \ar[drr, bend left, "{q_2}", ""'{name=q2, below}] \\
&P \ar[d, "{p_1}"'] \ar[r, "{p_2}"] \ar[Rightarrow, to=q1, "{\zeta_1}"'] \ar[r, "{p_2}"] \ar[Rightarrow, from=q2, "{\zeta_2}"]
&B \ar[d, "{g}"] \ar[dl, Rightarrow, "{\omega}", "{\equiv}"'] \\
&A \ar[r, "{f}"']
&C
\end{tikzcd}
\]
		\item for all 1-cells $u,v: T\to P$, 2-cells $\alpha:p_1u\Rightarrow p_1v, \beta: p_2u\Rightarrow p_2v$ and an invertible 3-cell 
		$\kappa:
		\frac{\Id_{g}\otimes\beta}{\omega\otimes\Id_{v}} \Rrightarrow \frac{\omega\otimes\Id_{u}}{\Id_{f}\otimes\alpha}$
		there are a 2-cell $\gamma: u\Rightarrow v$ and isomorphism 3-cells $\Gamma_1:\Id_{p_1}\ot\gamma\Rightarrow\alpha, 
		\Gamma_2:\Id_{p_2}\ot\gamma\Rightarrow\beta$ such that
	\[
\begin{tikzcd}
	\frac{(\Id_{g} \otimes \Id_{p_2}) \otimes \gamma}{\omega\otimes\Id_{v}}
		\ar[r, triple, "{\upepsilon_{\omega, \gamma}}"]
		\ar[d, triple, "{\frac{a^{-1}}{\Id}}"']
	&\frac{\omega\otimes\Id_{u}}{(\Id_{f} \otimes \Id_{p_1}) \otimes \gamma}
		\ar[d, triple, "{\frac{\Id}{a^{-1}}}"] \\ 
	\frac{\Id_{g} \otimes (\Id_{p_2} \otimes \gamma)}{\omega\otimes\Id_{v}}
		\ar[d, triple, "{\frac{\Id_{\Id_{g}} \otimes \Gamma_2}{\Id}}"']
	&\frac{\omega\otimes\Id_{u}}{\Id_{f} \otimes (\Id_{p_1} \otimes \gamma)}
		\ar[d, triple, "{\frac{\Id}{\Id_{\Id_{f}} \otimes \Gamma_1}}"] \\ 
	\frac{\Id_{g}\otimes\beta}{\omega\otimes\Id_{v}}
		\ar[r, triple, "{\upkappa}"]
	&\frac{\omega\otimes\Id_{u}}{\Id_{f}\otimes\alpha}
\end{tikzcd}
\]
commutes (where \(\upepsilon_{\omega, \gamma}\) is the interchange 3-cell).
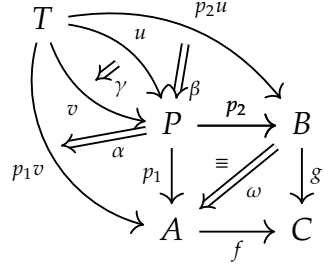
\[
\begin{tikzcd}
T
	\ar[dr, bend right, "{v}"', ""{name="u", above}]
	\ar[dr, bend left, "{u}", ""'{name="v", below}]
	\ar[ddr, bend right=45, "{p_1 v}"', ""{name=p1v, above}]
	\ar[drr, bend left, "{p_2 u}", ""'{name=p2u, below}] \\
&P
	\ar[d, "{p_1}"'] \ar[r, "{p_2}"]
	\ar[Rightarrow, to=p1v, "{\alpha}"]
	\ar[r, "{p_2}"]
	\ar[Rightarrow, from=p2u, "{\beta}"]
&B \ar[d, "{g}"] \ar[dl, Rightarrow, "{\omega}", "{\equiv}"']\\
&A \ar[r, "{f}"']
&C
	\ar[from="v", to="u", Rightarrow, "{\gamma}"]
\end{tikzcd}
\]
		
		\item for all 2-cells $\gamma,\gamma': u\Rightarrow v$ and 3-cells
\(\upchi_1 \colon \Id_{p_1} \otimes \gamma \Rrightarrow \Id_{p_1} \otimes \gamma'\) and
\(\upchi_2 \colon \Id_{p_2} \otimes \gamma \Rrightarrow \Id_{p_2} \otimes \gamma'\) such that
\[
\begin{tikzcd}
	\frac{\Id_{g} \otimes (\Id_{p_2} \otimes \gamma)}{\omega\otimes\Id_{v}}
		\ar[r, triple, "{\upepsilon_{\omega, \gamma}}"]
		\ar[d, triple, "{\frac{\Id_{g} \otimes \upchi_2}{\Id}}"']
	&\frac{\omega\otimes\Id_{u}}{\Id_{f} \otimes (\Id_{p_1} \otimes \gamma)}
		\ar[d, triple, "{\frac{\Id}{\Id_{f} \otimes \upchi_1}}"] \\ 
	\frac{\Id_{g} \otimes (\Id_{p_2} \otimes \gamma')}{\omega\otimes\Id_{v}}
		\ar[d, triple, "{\frac{a}{\Id}}"']
	&\frac{\omega\otimes\Id_{u}}{\Id_{f} \otimes (\Id_{p_1} \otimes \gamma')}
		\ar[d, triple, "{\frac{\Id}{a}}"] \\
	\frac{(\Id_{g} \otimes \Id_{p_2}) \otimes \gamma'}{\omega\otimes\Id_{v}}
		\ar[r, triple, "{\upepsilon_{\omega, \gamma'}}"]
	&\frac{\omega\otimes\Id_{u}}{(\Id_{f} \otimes \Id_{p_1}) \otimes \gamma'}
\end{tikzcd}
\]
commutes, there exists a unique 3-cell \(\upchi \colon \gamma \Rrightarrow \gamma'\) such that
\(\upchi_1 = \Id_{\Id_{p_1}} \otimes \upchi\) and
\(\upchi_2 = \Id_{\Id_{p_2}} \otimes \upchi\).
	\end{enumerate}
\end{defn}


\medskip





For  convenience we write out the definition of a tricategorical product from \cite{Fem} that we call 3-product for short. 
It is the dual of Definition 5.3 from {\em loc.cit.}.

\begin{defn} \delabel{3-pr}
A 3-product of 0-cells $A$ and $B$ in a tricategory $V$ consists of: 
a 0-cell $A\times B$ and 1-cells $p_1:A\times B\to A, p_2:A\times B\to B$, such that
 \begin{enumerate}
  \item for every 0-cell $T$ and 1-cells $f_1: T\to A, f_2: T\to B$ 
	      there are a 1-cell $u: T\to A\times B$ and equivalence 2-cells $\zeta_i: f_i\Rightarrow p_iu, i=1,2$; 
  		 \item 
			for all 1-cells $u,v: T\to A\times B$ and 2-cells $\alpha: p_1u\Rightarrow p_1v$ and $\beta: p_2u\Rightarrow p_2v$, 
	     there are a 2-cell $\gamma:u\Rightarrow v$ and isomorphism 3-cells $\Gamma_1: \alpha\Rrightarrow \Id_{p_1}\ot\gamma$ and 
			$\Gamma_2: \beta\Rrightarrow \Id_{p_2}\ot\gamma$; 
				\item 
				for every two 2-cells $\gamma,\gamma':u\Rightarrow v$ and every two 3-cells 
				$\chi_i: \Id_{p_i}\ot\gamma\Rrightarrow \Id_{p_i}\ot\gamma', i=1,2$ there is a unique 3-cell 
				$\Gamma: \gamma\Rrightarrow\gamma'$ such that $\chi_i=\Id_{\Id_{p_i}}\ot\Gamma, i=1,2$.
       \end{enumerate}
\end{defn}

Corresponding to \equref{Phi cond} and with notations as in items 1) and 2) above, one has that for all 2-cells $(\theta_1,\theta_1): 
(f_1, f_2)\Rightarrow (f_1', f_2')$ there are invertible 3-cells 
\begin{equation} \eqlabel{Omega_i}
\Omega_1:\frac{\zeta_1}{\Id_{p_1}\ot\theta}\Rrightarrow\frac{\theta_1}{\zeta_1'}\quad \text{and}\quad 
\Omega_2:\frac{\zeta_2}{\Id_{p_2}\ot\theta}\Rrightarrow\frac{\theta_2}{\zeta_2'}
\end{equation}
where $\zeta_i':f_i'\Rightarrow p_iv, i=1,2$ and $\theta=G(\theta_1,\theta_2)$ where 
$G: V(T,a)\times V(T,B) \to V(T,A\times B)$ is a biequivalence.

\begin{rem} \rmlabel{uniqueness}
The 3-limits in the above two definitions are unique up to a {\em 1-cell invertible up to a 2-equivalence}. 
We will call such cells {\em biequivalence 1-cells}. 
The latter means a 1-cell $f: L\to L'$ for which there exists a 1-cell $g:L'\to L$ and an equivalence 2-cell $\omega: gf\Rightarrow id_L$. 
Likewise, in the item 1) of the definitions the 1-cells whose existence is claimed are unique up to equivalence 2-cells, and 
the equivalence 2-cells whose existence is claimed are unique up to isomorphism. 
In the item 2) of the definitions the 2-cells $\gamma$ 
whose existence is claimed are unique up to a unique isomorphism. 
We record this in the next Corollary. 
\end{rem}

\begin{cor} \colabel{unique 3-cell Gamma}
Let $A\times B$ be a 3-product in a tricategory $V$. Then in reference to the items in the above Definition it is:
\begin{enumerate}
\item in the item 1) the 1-cell $u$ is unique up to an equivalence 2-cell, and 
the equivalence 2-cells $\zeta_i$ are unique up to isomorphism;
\item if $\alpha, \beta$ as in the item 2) induce 2-cells $\gamma,\gamma':u\Rightarrow v$, then there is a unique isomorphism 3-cell 
$\Gamma: \gamma\Rrightarrow\gamma'$;
\item given 2-cells $\gamma,\gamma':u\Rightarrow v$ and invertible 3-cells 
				$\chi_i: \Id_{p_i}\ot\gamma\Rrightarrow \Id_{p_i}\ot\gamma', i=1,2$, then the unique 3-cell 
				$\Gamma: \gamma\Rrightarrow\gamma'$ from the item 3) is invertible.
      
\end{enumerate}
\end{cor}

An important direct consequence of the definition is: 

\begin{lma} \lelabel{equiv cells}
If the 2-cells $\alpha: p_1u\Rightarrow p_1v$ and $\beta: p_2u\Rightarrow p_2v$ in the part 2) in \deref{3-pr} are equivalence 2-cells, 
then so is $\gamma: u\Rightarrow v$. Namely, quasi-inverses $\alpha^{-1}$ and $\beta^{-1}$ induce a quasi-inverse $\gamma^{-1}$. 
\end{lma}

Analogous results to the above two hold also for 3-pullbacks. 
In view of these results, for any quasi-inverse of an equivalence 2-cell $\zeta$ obtained in the context of these 3-limits 
we will write simply $\zeta^{-1}$, throughout, without any further reference to a choice in the respective isomorphism class.

\subsection{Inducing 3-products on higher cells}

In this Subsection we will induce 1-cells (2-cells) $x\times y$ for given 1-cells (2-cells) $x$ and $y$. Also, given certain 3-cells 
$P_{\alpha,\alpha'}$ and $P_{\beta, \beta'}$, where $P_{\alpha,\alpha'}$ can be thought of as specific (transversal) prism whose bases 
is a 2-cell $\alpha$ (vertically in the back) and the opposite face is a 2-cell $\alpha'$ (vertically in the front) and analogously for 
$P_{\beta, \beta'}$, we will induce a 3-cell $P_{\alpha\times\beta, \alpha'\times\beta'}$. In order not to make this lengthy paper even longer, 
we will skip the proofs in this Subsection and will present only the results that we obtained.

\begin{lma} \lelabel{lema}
Let $V$ be a 1-strict tricategory. Given two 3-products $M\times N$ and $P\times Q$ of 0-cells in $V$, one has: 
\begin{enumerate}[a)]
\item given objects $T,S$ and 1-cells $q_1:T\to M, q_2:T\to N$ and $s_1:S\to P, s_2:S\to Q$, then there are 1-cells $t:T\to M\times N$ and 
$s:S\to P\times Q$ and equivalence 2-cells $\zeta_1:q_1\Rightarrow p_1t, \zeta_2:q_2\Rightarrow p_2t$ and $\theta_1:s_1\Rightarrow p'_1s, 
\theta_2:s_2\Rightarrow p'_2s$ as in the picture below; 
\item given 1-cells $g:M\to P$ and $h:N\to Q$, then there is a 1-cell that we will denote by $g\times h$ acting between 3-products 
$M\times N\to P\times Q$ and there are equivalence 2-cells $\omega_1:gp_1\Rightarrow p'_1(g\times h)$ and $\omega_2:hp_2\Rightarrow p'_2(g\times h)$ 
as in the picture below; 
\item additionally to the data from a) and b), given 1-cell $f:T\to S$ and 2-cells $\alpha:gq_1\Rightarrow s_1f$ and $\beta:hq_2\Rightarrow s_2f$, then 
there is a 2-cell $\gamma:(g\times h)t\Rightarrow sf$ and invertible 3-cells  
$$\Gamma^1_\gamma:\threefrac{\Id_g\ot\zeta_1}{\omega_1\ot\Id_t}{\Id_{p_1'}\ot\gamma}\Rrightarrow\frac{\alpha}{\theta_1\ot\Id_f},\qquad
 \Gamma^2_\gamma:\threefrac{\Id_h\ot\zeta_2}{\omega_2\ot\Id_t}{\Id_{p_2'}\ot\gamma}\Rrightarrow\frac{\beta}{\theta_2\ot\Id_f};$$
\item if $\alpha$ and $\beta$ in c) are equivalence 2-cells, then so is $\gamma$. 
\end{enumerate}
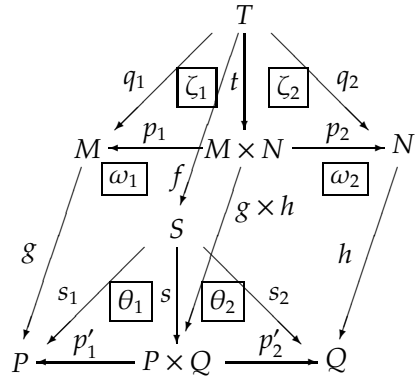
$$
 \scalebox{0.86}{\bfig
\putmorphism(400,300)(-1,-1)[``q_1]{500}1l
 \putmorphism(480,300)(1,-1)[\phantom{\u\C \bullet 1}` `q_2]{500}1r
 \putmorphism(-130,-200)(1,0)[M`\phantom{M\times N} ` p_1]{600}{-1}a
 \putmorphism(490,-200)(1,0)[\phantom{\u\C \bullet 1}`N ` p_2]{550}1a
\putmorphism(450,300)(0,-1)[T`M\times N`]{500}1l
\putmorphism(470,300)(0,-1)[`\phantom{M\times N}`t]{500}0l
\put(200,0){\fbox{$\zeta_1$}}
\put(540,0){\fbox{$\zeta_2$}}
\put(-80,-330){\fbox{$\omega_1$}}
\put(740,-330){\fbox{$\omega_2$}}
\putmorphism(450,300)(-1,-3)[``]{260}1l
\putmorphism(420,80)(-1,-3)[``f]{260}0l
\putmorphism(-130,-200)(-1,-3)[``g]{260}1l
\putmorphism(1040,-200)(-1,-3)[``h]{260}1l
\putmorphism(460,-200)(-1,-3)[``]{250}1r
\putmorphism(480,-50)(-1,-3)[``g\times h]{250}0r
\putmorphism(150,-500)(-1,-1)[``s_1]{500}1l
 \putmorphism(230,-500)(1,-1)[\phantom{\u\C \bullet 1}` `s_2]{500}1r
 \putmorphism(-380,-1000)(1,0)[P`\phantom{M\times N} ` p'_1]{600}{-1}a
 \putmorphism(240,-1000)(1,0)[\phantom{\u\C \bullet 1}`Q ` p'_2]{550}1a
\putmorphism(200,-500)(0,-1)[S`P\times Q`]{500}1l
\putmorphism(220,-500)(0,-1)[`\phantom{P\times Q}`s]{500}0l
\put(-50,-800){\fbox{$\theta_1$}}
\put(290,-800){\fbox{$\theta_2$}}
\efig}
$$
\end{lma}

{\bf Convention.} The 3-cells $\Gamma_1$ and $\Gamma_2$ we will call informally and for short {\em prisms $P_i$ from $\zeta_i$ to $\theta_i$}, 
$i=1,2$. 

\medskip

\begin{rem} \rmlabel{opposite 2-cells}
Since in the item c) the 1-cells $f,g,h$ in our drawing go transversally from back to front, we considered 2-cells $\alpha$ and $\beta$ in the mapping order we did (from $q_i$ to $\theta_i, i=1,2$, so to say). Though, this way we get: 
$$
\bfig
\putmorphism(0,130)(1,0)[``f]{420}1a
\putmorphism(0,160)(0,-1)[\phantom{Y_2}``q_1]{350}1l
\putmorphism(400,160)(0,-1)[\phantom{Y_2}``s_1]{350}1r
\put(170,-20){\fbox{$\alpha$}}
\put(20,-110){$\Nearrow$}
\putmorphism(0,-150)(1,0)[``g]{420}1b
\efig
\qquad\text{and}\qquad
\bfig
\putmorphism(0,130)(1,0)[``f]{420}1a
\putmorphism(0,160)(0,-1)[\phantom{Y_2}``q_2]{350}1l
\putmorphism(400,160)(0,-1)[\phantom{Y_2}``s_2]{350}1r
\put(170,-20){\fbox{$\beta$}}
\put(20,-110){$\Nearrow$}
\putmorphism(0,-150)(1,0)[`.`h]{420}1b
\efig
$$
However, since we usually consider 2-cells written in squares in the direction $\Swarrow$, we make the following remark. 
If one considers the 2-cells $\alpha, \beta$ in c) in the reversed order, namely: $\alpha:s_1f\Rightarrow gq_1$ and $\beta:s_2f\Rightarrow hq_2$, 
whereas $\zeta_i, \theta_i$ and $\omega_i, i=1,2$ maintain the same order, then 
one gets a 2-cell $\gamma$ accordingly, namely $\gamma:sf\Rightarrow (g\times h)t$. 
\end{rem}

\begin{rem} \rmlabel{converse for 3-pr}
In all the three parts of \deref{3-pr} the existence of the announced cells is subject to the existence of previously mentioned cells that determine them. 
Thus given a 2-cell $\gamma$ as in \leref{lema} c) it is understood that $\gamma$ comes together with some 2-cells $\alpha$ and $\beta$ that determine it. 
\end{rem}

We will next construct a prism with basis a 2-cell $\gamma$ as in the part c) of \leref{lema} whose opposite face is an analogous 2-cell 
$\gamma'$, out of prisms with bases $\alpha$ and $\beta$.

\begin{prop} \prlabel{cutting 3-cells}
Let $V$ be a 1-strict tricategory.  
Let $\gamma: (g\times h)t\Rightarrow sf$ be as in \leref{lema} c) with its assigned 2-cells $\zeta_i, \theta_i, \omega_i, i=1,2, 
\alpha$ and $\beta$, 
and consider another $\gamma': (g'\times h')t\Rightarrow sf'$ induced by $\alpha': g'q_1\Rightarrow s_1f'$ and $\beta':h'q_2\Rightarrow s_2f'$ with 
$\omega_1': g'p_1\Rightarrow p'_1(g'\times h')$ and $\omega_2':h'p_2\Rightarrow p'_2(g'\times h')$ and the same 2-cells $\zeta_i, \theta_i, i=1,2$. 
Suppose there are 2-cells $\xi:f\Rightarrow f', \lambda: g\Rightarrow g', \rho: h\Rightarrow h'$ and 3-cells 
$$
\bfig
  \putmorphism(0,200)(1,0)[``\Id]{450}1a
 \putmorphism(340,200)(1,0)[\phantom{F(B)}` `q_1]{600}1a

\putmorphism(0,220)(0,-1)[\phantom{Y_2}``f']{400}1l
\putmorphism(450,220)(0,-1)[\phantom{Y_2}``f]{400}1l
\putmorphism(900,220)(0,-1)[\phantom{Y_2}``g]{400}1r
\put(130,0){\fbox{$\xi$}}
\put(600,0){\fbox{$\alpha$}}
\putmorphism(0,-150)(1,0)[``\Id]{450}1b
 \putmorphism(340,-150)(1,0)[\phantom{F(B)}` `s_1]{580}1b
\put(310,90){$\Swarrow$}
\put(760,90){$\Swarrow$}
\efig
\quad\stackrel{P_\alpha}{\Rrightarrow}\quad
\bfig
  \putmorphism(0,200)(1,0)[``q_1]{450}1a
 \putmorphism(340,200)(1,0)[\phantom{F(B)}` `\Id]{550}1a

\putmorphism(0,220)(0,-1)[\phantom{Y_2}``f']{400}1l
\putmorphism(450,220)(0,-1)[\phantom{Y_2}``g']{400}1l
\putmorphism(890,220)(0,-1)[\phantom{Y_2}``g]{400}1r
\put(130,0){\fbox{$\alpha'$}}
\put(610,0){\fbox{$\lambda$}}
\put(310,100){$\Swarrow$}
\put(760,100){$\Swarrow$}
\putmorphism(0,-150)(1,0)[``s_1]{450}1b
 \putmorphism(340,-150)(1,0)[\phantom{F(B)}` `\Id]{550}1b
\efig
$$ 
and
$$
\bfig
  \putmorphism(0,200)(1,0)[``\Id]{450}1a
 \putmorphism(340,200)(1,0)[\phantom{F(B)}` `q_2]{600}1a

\putmorphism(0,220)(0,-1)[\phantom{Y_2}``f']{400}1l
\putmorphism(450,220)(0,-1)[\phantom{Y_2}``f]{400}1l
\putmorphism(900,220)(0,-1)[\phantom{Y_2}``h]{400}1r
\put(130,0){\fbox{$\xi$}}
\put(600,0){\fbox{$\beta$}}
\putmorphism(0,-150)(1,0)[``\Id]{450}1b
 \putmorphism(340,-150)(1,0)[\phantom{F(B)}` `s_2]{580}1b
\put(310,90){$\Swarrow$}
\put(760,90){$\Swarrow$}
\efig
\quad\stackrel{P_\beta}{\Rrightarrow}\quad
\bfig
  \putmorphism(0,200)(1,0)[``q_2]{450}1a
 \putmorphism(340,200)(1,0)[\phantom{F(B)}` `\Id]{550}1a

\putmorphism(0,220)(0,-1)[\phantom{Y_2}``f']{400}1l
\putmorphism(450,220)(0,-1)[\phantom{Y_2}``h']{400}1l
\putmorphism(890,220)(0,-1)[\phantom{Y_2}``h]{400}1r
\put(130,0){\fbox{$\beta'$}}
\put(610,0){\fbox{$\rho$}}
\put(310,100){$\Swarrow$}
\put(760,100){$\Swarrow$}
\putmorphism(0,-150)(1,0)[``s_2]{450}1b
 \putmorphism(340,-150)(1,0)[\phantom{F(B)}` .`\Id]{550}1b
\efig
$$ 
Then there is a unique 3-cell 
$$
\bfig
  \putmorphism(0,200)(1,0)[``\Id]{450}1a
 \putmorphism(340,200)(1,0)[\phantom{F(B)}` `t]{600}1a

\putmorphism(0,220)(0,-1)[\phantom{Y_2}``f']{400}1l
\putmorphism(450,220)(0,-1)[\phantom{Y_2}``f]{400}1l
\putmorphism(900,220)(0,-1)[\phantom{Y_2}``g\times h]{400}1r
\put(130,0){\fbox{$\xi$}}
\put(600,0){\fbox{$\gamma$}}
\putmorphism(0,-150)(1,0)[``\Id]{450}1b
 \putmorphism(340,-150)(1,0)[\phantom{F(B)}` `s]{580}1b
\put(310,90){$\Swarrow$}
\put(760,90){$\Swarrow$}
\efig
\quad\stackrel{\Gamma}{\Rrightarrow}\quad
\bfig
  \putmorphism(0,200)(1,0)[``t]{570}1a
 \putmorphism(450,200)(1,0)[\phantom{F(B)}` `\Id]{620}1a

\putmorphism(0,220)(0,-1)[\phantom{Y_2}``f']{400}1l
\putmorphism(550,220)(0,-1)[\phantom{Y_2}``]{400}1l
\putmorphism(570,200)(0,-1)[\phantom{Y_2}``g'\times h']{400}0l
\putmorphism(1050,220)(0,-1)[\phantom{Y_2}``g\times h]{400}1r
\put(80,0){\fbox{$\gamma'$}}
\put(610,0){\fbox{$\lambda\times\rho$}}
\put(410,100){$\Swarrow$}
\put(930,100){$\Swarrow$}
\putmorphism(0,-150)(1,0)[``s]{570}1b
 \putmorphism(450,-150)(1,0)[\phantom{F(B)}` `\Id]{620}1b
\efig
$$
{\em i.e.} $\Gamma: \displaystyle{\frac{[\Id\vert\gamma]}{[\xi\vert\Id_s]}} \Rrightarrow \displaystyle{\frac{[\Id_t\vert\lambda\times\rho]}{[\gamma'\vert\Id]}}$. Moreover, if $P_\alpha, P_\beta$ are invertible, then so is $\Gamma$. 
\end{prop}

As a consequence we may formulate:

\begin{cor} \colabel{cutting 3-cell-times}
Given 1-cells $u,v: T\to A\times B$ and 2-cells $\alpha: p_1u\Rightarrow p_1v$ and $\beta: p_2u\Rightarrow p_2v$, which induce a 2-cell 
$\gamma: u\Rightarrow v$, and similarly assume that $\alpha': p_1u'\Rightarrow p_1v'$ and $\beta: p_2u'\Rightarrow p_2v'$ induce 
$\gamma':u\Rightarrow v'$, for $u',v': T\to A\times B$. Suppose there are 2-cells $\xi:u\Rightarrow u'$ and $\zeta:v\Rightarrow v'$ 
and 3-cells $P_{\tilde\alpha}: \frac{\tilde\alpha}{[\zeta\vert\Id_{p_1}]}\Rrightarrow\frac{[\xi\vert\Id_{p_1}]}{\tilde\alpha'}$ and 
$P_{\tilde\beta}: \frac{\tilde\beta}{[\zeta\vert\Id_{p_2}]}\Rrightarrow\frac{[\xi\vert\Id_{p_2}]}{\tilde\beta'}$. Then there is a unique 
3-cell $\Gamma: \frac{\gamma}{\zeta}\Rrightarrow \frac{\xi}{\gamma'}$.  If $P_{\tilde\alpha}, P_{\tilde\beta}$ are invertible, so is $\Gamma$. 
\end{cor}

Let us formulate the dualization of \prref{cutting 3-cells}, whereas $\alpha$ and $\beta$ are considered in the reversed order, as indicated in 
\rmref{opposite 2-cells}, and in the setting of a 3-coproduct $\amalg_{i\in I}A_i$ for a set $I$ (instead of a 3-coproduct $A\amalg B$). One has:

\begin{prop} \prlabel{dual cutting}
Let $V$ be a 1-strict tricategory. 
Let $\alpha_i: t_{A_i'}g_i\Rightarrow fs_{A_i}$, for $i\in I$, induce $\gama: t(\amalg_{i\in I}g_i)\Rightarrow fs$ in the diagram below:
$$
 \scalebox{0.86}{\bfig
 \putmorphism(410,300)(1,0)[\phantom{\u\C \bullet 1}` \amalg_iA_i`]{580}1a 
\putmorphism(470,300)(1,-1)[A_i``]{550}1l
\putmorphism(490,300)(1,-1)[``s_{A_i}]{550}0l
\putmorphism(450,300)(-1,-3)[`A_i'`]{260}1l
\putmorphism(430,160)(-1,-3)[``g_i]{260}0l
\putmorphism(1000,-200)(-1,-3)[``f]{240}1l
 \putmorphism(150,-500)(1,0)[\phantom{\u\C \bullet 1}`\amalg_iA_i' `]{580}1r 
\putmorphism(180,-450)(1,-1)[\phantom{P\times Q}`T`]{550}1l
\putmorphism(230,-500)(1,-1)[\phantom{P\times Q}``t_{A_i'}]{550}0l

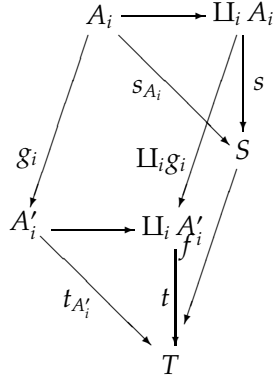
\putmorphism(1000,300)(0,-1)[`S`s]{500}1r
\putmorphism(750,-500)(0,-1)[``]{500}1l
\putmorphism(770,-500)(0,-1)[``t]{500}0l
\putmorphism(1000,300)(-1,-3)[``]{260}1l
\putmorphism(980,150)(-1,-3)[``\amalg_ig_i]{260}0l
\efig}
$$
and similarly let $\alpha_i': t_{A_i'}g_i'\Rightarrow f's_{A_i}, i\in I$, induce $\gama': t(\amalg_{i\in I}g_i')\Rightarrow f's$. 
Suppose moreover that there are 2-cells $\xi:f\Rightarrow f'$ and $\sigma_i:g_i\Rightarrow g_i'$ and 3-cells 
$$
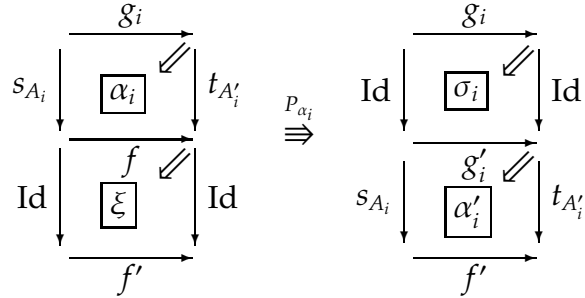
\bfig
\putmorphism(0,350)(1,0)[``g_i]{450}1a
\putmorphism(0,370)(0,-1)[\phantom{Y_2}``s_{A_i}]{400}1l
\putmorphism(430,370)(0,-1)[\phantom{Y_2}``t_{A_i'}]{400}1r
\put(130,140){\fbox{$\alpha_i$}}
\put(130,-230){\fbox{$\xi$}}
\put(310,250){$\Swarrow$}
\put(310,-90){$\Swarrow$}
\putmorphism(0,50)(0,-1)[\phantom{Y_2}``\Id]{450}1l
\putmorphism(430,50)(0,-1)[\phantom{Y_2}``\Id]{450}1r
\putmorphism(0,10)(1,0)[``f]{450}1b
\putmorphism(0,-370)(1,0)[``f']{450}1b
\efig
\quad\stackrel{P_{\alpha_i}}{\Rrightarrow}\quad
\bfig
\putmorphism(0,350)(1,0)[``g_i]{450}1a
\putmorphism(0,370)(0,-1)[\phantom{Y_2}``\Id]{400}1l
\putmorphism(430,370)(0,-1)[\phantom{Y_2}``\Id]{400}1r
\put(130,150){\fbox{$\sigma_i$}}
\put(130,-240){\fbox{$\alpha_i'$}}
\put(310,240){$\Swarrow$}
\put(310,-100){$\Swarrow$}
\putmorphism(0,10)(0,-1)[\phantom{Y_2}``s_{A_i}]{400}1l
\putmorphism(430,10)(0,-1)[\phantom{Y_2}``t_{A_i'}]{400}1r
\putmorphism(0,0)(1,0)[``g_i']{450}1b
\putmorphism(0,-370)(1,0)[``f']{450}1b
\efig
$$
for every $i\in I$. Then there is a unique 3-cell 
$$
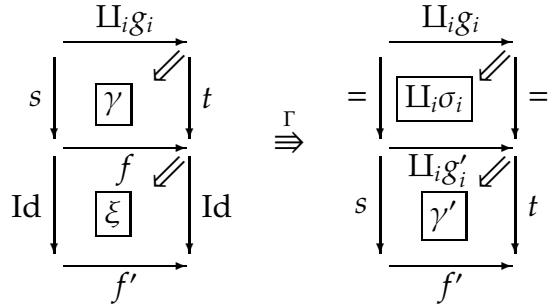
\bfig
\putmorphism(0,350)(1,0)[``\amalg_ig_i]{450}1a
\putmorphism(0,370)(0,-1)[\phantom{Y_2}``s]{400}1l
\putmorphism(430,370)(0,-1)[\phantom{Y_2}``t]{400}1r
\put(130,140){\fbox{$\gamma$}}
\put(130,-230){\fbox{$\xi$}}
\put(310,250){$\Swarrow$}
\put(310,-90){$\Swarrow$}
\putmorphism(0,50)(0,-1)[\phantom{Y_2}``\Id]{450}1l
\putmorphism(430,50)(0,-1)[\phantom{Y_2}``\Id]{450}1r
\putmorphism(0,10)(1,0)[``f]{450}1b
\putmorphism(0,-370)(1,0)[``f']{450}1b
\efig
\quad\stackrel{\Gamma}{\Rrightarrow}\quad
\bfig
\putmorphism(0,350)(1,0)[``\amalg_ig_i]{450}1a
\putmorphism(0,370)(0,-1)[\phantom{Y_2}``=]{400}1l
\putmorphism(430,370)(0,-1)[\phantom{Y_2}``=]{400}1r
\put(50,140){\fbox{$\amalg_i\sigma_i$}}
\put(130,-230){\fbox{$\gamma'$}}
\put(310,250){$\Swarrow$}
\put(310,-90){$\Swarrow$}
\putmorphism(0,50)(0,-1)[\phantom{Y_2}``s]{450}1l
\putmorphism(430,50)(0,-1)[\phantom{Y_2}``t]{450}1r
\putmorphism(0,10)(1,0)[``]{450}1b
\putmorphism(-40,10)(1,0)[``\amalg_ig_i']{450}0b
\putmorphism(0,-370)(1,0)[``f']{450}1b
\efig
$$
If $P_{\alpha_i}, i\in I$ are invertible, then so is $\Gamma$. 
\end{prop}

Observe here that $g,h$ on one hand and $\lambda, \rho$ on the other, from \prref{cutting 3-cells}, pass to $g_i$, respectively $\sigma_i$. 
Consequently, $\lambda\times\rho$ passes to $\amalg_i \sigma_i$. 

The dual of \coref{cutting 3-cell-times} is:

\begin{cor} \colabel{cutting 3-cell-coproduct}
Given 1-cells $u,v: \amalg_i A_i\to T$ and 2-cells $\alpha_i: u\iota_i\Rightarrow v\iota_i$ for every $i\in I$, which induce a 2-cell 
$\gamma: u\Rightarrow v$, and similarly assume that $\alpha_i': u'\iota_i\Rightarrow v'\iota_i$ induce 
$\gamma':u'\Rightarrow v'$, for $u',v': \amalg_i A_i\to T$. Suppose there are 2-cells $\xi:u\Rightarrow u'$ and $\zeta:v\Rightarrow v'$ 
and 3-cells $P_{\alpha_i}: \frac{\alpha_i}{[\Id_{\iota_i}\vert\zeta]}\Rrightarrow\frac{[\Id_{\iota_i}\vert\xi]}{\alpha_i'}, i\in I$, 
then there is a unique 3-cell $\Gamma: \frac{\gamma}{\zeta}\Rrightarrow \frac{\xi}{\gamma'}$. If $P_{\alpha_i}, i\in I$ are invertible, 
so is $\Gamma$. 

\end{cor}

\medskip

We now induce a 3-product 2-cell, {\em i.e.} a 2-cell $\alpha\times\beta$ obtained as a 3-product of two 2-cells $\alpha, \beta$.

\begin{cor} (Corollary of \leref{lema}) \colabel{alfa-x-beta} \\
Let $V$ be a 1-strict tricategory. 
Given four 3-products $M\times N, M'\times N', P\times Q$ and $P'\times Q'$ of 0-cells and 2-cells:
$$
\bfig
\putmorphism(-150,50)(1,0)[M` P `a]{450}1a
\putmorphism(-150,-300)(1,0)[M'` P' `a']{440}1b
\putmorphism(-170,50)(0,-1)[\phantom{Y_2}``m]{350}1l
\putmorphism(280,50)(0,-1)[\phantom{Y_2}``p]{350}1r
\put(0,-140){\fbox{$\alpha$}}
\efig
\qquad
\bfig
\putmorphism(-150,50)(1,0)[N` Q `b]{450}1a
\putmorphism(-150,-300)(1,0)[N'` Q' `b']{440}1b
\putmorphism(-170,50)(0,-1)[\phantom{Y_2}``n]{350}1l
\putmorphism(280,50)(0,-1)[\phantom{Y_2}``q]{350}1r
\put(0,-140){\fbox{$\beta$}}
\efig
$$
there is a 2-cell
$$
\bfig
\putmorphism(-150,50)(1,0)[M\times N` P\times Q `a\times b]{650}1a
\putmorphism(-150,-300)(1,0)[M'\times N'` P'\times Q' ` a'\times b']{640}1b
\putmorphism(-170,50)(0,-1)[\phantom{Y_2}``m\times n]{350}1l
\putmorphism(480,50)(0,-1)[\phantom{Y_2}``p\times q]{350}1r
\put(0,-140){\fbox{$\alpha\times\beta$}}
\efig
$$
and two isomorphism 3-cells, which are the evident transversal prisms from back to front in the diagram:
$$
 \scalebox{0.86}{\bfig
\putmorphism(-100,300)(-1,0)[``q_1]{400}{-1}a
 \putmorphism(510,300)(1,0)[\phantom{\u\C \bullet 1}` N`q_2]{500}1a
 \putmorphism(-160,-200)(1,0)[M'`\phantom{M\times N} ` p_1]{600}{-1}a
 \putmorphism(490,-200)(1,0)[\phantom{\u\C \bullet 1}`N' ` p_2]{500}1a
\putmorphism(470,300)(0,-1)[M\times N`M'\times N'`m\times n]{500}1r
\putmorphism(450,300)(-1,-3)[``]{260}1l
\putmorphism(410,80)(-1,-3)[``a\times b]{260}0l
\putmorphism(-110,-200)(-1,-3)[``]{260}1l
\putmorphism(-180,-340)(-1,-3)[``a']{260}0r
\putmorphism(1000,-200)(-1,-3)[``]{240}1r
\putmorphism(930,-360)(-1,-3)[``b']{260}0r
\putmorphism(460,-200)(-1,-3)[``]{250}1r
\putmorphism(360,-350)(-1,-3)[``a'\times b']{250}0r
\putmorphism(-400,-500)(-1,0)[``s_1]{450}{-1}l
 \putmorphism(250,-500)(1,0)[\phantom{\u\C \bullet 1}`Q `s_2]{500}1r
 \putmorphism(-400,-1000)(1,0)[P'`\phantom{M\times N} ` p'_1]{600}{-1}b
 \putmorphism(240,-1000)(1,0)[\phantom{\u\C \bullet 1}`Q' ` p'_2]{500}1b
\putmorphism(200,-500)(0,-1)[P\times Q`P'\times Q'`]{500}1l
\putmorphism(220,-500)(0,-1)[\phantom{P\times Q}`\phantom{P'\times Q'}`p\times q]{500}0l
\putmorphism(-160,300)(-1,-3)[M`P`]{260}1l
\putmorphism(-210,80)(-1,-3)[``a]{260}0l
\putmorphism(-150,300)(0,-1)[``m]{500}1r
\putmorphism(1000,300)(0,-1)[``n]{500}1r
\putmorphism(-400,-500)(0,-1)[``p]{500}1l
\putmorphism(750,-500)(0,-1)[``]{500}1l
\putmorphism(770,-500)(0,-1)[``q]{500}0l
\putmorphism(1000,300)(-1,-3)[``]{260}1r
\putmorphism(880,50)(-1,-3)[``b]{260}0r
\efig}
$$
If $\alpha$ and $\beta$ are equivalence 2-cells, then so is $\alpha\times\beta$. 
\end{cor}

\begin{rem} \rmlabel{opposit} 
If the 2-cells $\alpha$ and $\beta$ were in the reversed order, that is $\crta\alpha: s_1m\Rightarrow pa, \crta\beta: b'n\Rightarrow qb$,  
most importantly, if one considers the same squares as for $\alpha$ and $\beta$ above but directed in the opposite direction, 
whereas $\zeta_i, \theta_i,\omega_i$ and $\omega_i', i=1,2$ remain the same, 
then one gets a 2-cell $\crta\alpha\times\crta\beta$ and 3-cells 
$$\crta C^1_{\crta\alpha\times\crta\beta}: \threefrac{\Id_{p_1}\ot(\crta\alpha\times\crta\beta)}{\zeta_1^{-1}\ot\Id_{a\times b}}
{\Id_{m}\ot\omega_1^{-1}} \Rrightarrow \threefrac{\omega_1'^{-1}\ot\Id_{p\times q}}{\Id_{a'}\ot\theta_1^{-1}}{\crta\alpha\ot\Id_{s_1}},
\qquad
\crta C^2_{\crta\alpha\times\crta\beta}: \threefrac{\Id_{p_2}\ot(\crta\alpha\times\crta\beta)}{\zeta_2^{-1}\ot\Id_{a\times b}}
{\Id_{n}\ot\omega_2^{-1}} \Rrightarrow \threefrac{\omega_2'^{-1}\ot\Id_{p\times q}}{\Id_{b'}\ot\theta_2^{-1}}{\crta\beta\ot\Id_{s_2}}.
$$
\end{rem}

\begin{rem} \rmlabel{uniqueness2} 
By our remarks on uniqueness, observe that the 1-cell $g\times h$ in \leref{lema} is unique up to an equivalence 2-cell, 
and that the equivalence 2-cells $\omega_1, \omega_2$ from there and $\alpha\times\beta$ from \coref{alfa-x-beta} are unique up to isomorphism. 
\end{rem}

We finally construct a 3-cell with an associated prism with basis $\alpha\times\beta$ out of two given 3-cells with associated prisms 
of basis $\alpha$ and $\beta$.

\begin{prop} \prlabel{middle prism}
Given two 2-cells $\crta\alpha\times\crta\beta$ and $\crta\alpha'\times\crta\beta'$ with notations as in \coref{alfa-x-beta}, but with the 
orientation as in \rmref{opposit}. Suppose that further 2-cells $\lambda_1:m\Rightarrow m', \lambda_2:n\Rightarrow n', \rho_1: p\Rightarrow p', 
\rho_2:q\Rightarrow q'$ are given and two 3-cells 
$$
 \scalebox{0.86}{\bfig
  \putmorphism(0,180)(1,0)[``a]{450}1a
 \putmorphism(330,180)(1,0)[\phantom{F(B)}` `\Id]{600}1a

\putmorphism(0,200)(0,-1)[\phantom{Y_2}``m]{400}1l
\putmorphism(430,200)(0,-1)[\phantom{Y_2}``p]{400}1l
\putmorphism(900,200)(0,-1)[\phantom{Y_2}``p']{400}1r
\put(150,0){\fbox{$\crta\alpha$}}
\put(600,0){\fbox{$\rho_1$}}
\put(20,-100){$\Nearrow$}
\put(470,-100){$\Nearrow$}

\putmorphism(0,-160)(1,0)[``a']{450}1b
 \putmorphism(330,-160)(1,0)[\phantom{F(B)}` `\Id]{580}1b
\efig}
\quad\stackrel{P_1}{\Rrightarrow}\quad
 \scalebox{0.86}{\bfig
  \putmorphism(0,180)(1,0)[``\Id]{450}1a
 \putmorphism(360,180)(1,0)[\phantom{F(B)}` `a]{570}1a

\putmorphism(0,200)(0,-1)[\phantom{Y_2}``m]{400}1l
\putmorphism(460,200)(0,-1)[\phantom{Y_2}``]{400}1l
\putmorphism(490,200)(0,-1)[\phantom{Y_2}``m']{400}0l
\putmorphism(900,200)(0,-1)[\phantom{Y_2}``p']{400}1r
\put(150,0){\fbox{$\lambda_1$}}
\put(620,-20){\fbox{$\crta\alpha'$}}
\putmorphism(0,-160)(1,0)[``\Id]{450}1b
 \putmorphism(370,-160)(1,0)[\phantom{F(B)}` `a']{550}1b
\put(20,-100){$\Nearrow$}
\put(480,-100){$\Nearrow$}
\efig}
$$
and 
$$
 \scalebox{0.86}{\bfig
  \putmorphism(0,180)(1,0)[``b]{450}1a
 \putmorphism(330,180)(1,0)[\phantom{F(B)}` `\Id]{600}1a

\putmorphism(0,200)(0,-1)[\phantom{Y_2}``n]{400}1l
\putmorphism(430,200)(0,-1)[\phantom{Y_2}``q]{400}1l
\putmorphism(900,200)(0,-1)[\phantom{Y_2}``q']{400}1r
\put(150,0){\fbox{$\crta\beta$}}
\put(600,0){\fbox{$\rho_2$}}

\putmorphism(0,-160)(1,0)[``b']{450}1b
 \putmorphism(330,-160)(1,0)[\phantom{F(B)}` `\Id]{580}1b
\put(20,-100){$\Nearrow$}
\put(470,-100){$\Nearrow$}
\efig}
\quad\stackrel{P_2}{\Rrightarrow}\quad
 \scalebox{0.86}{\bfig
  \putmorphism(0,180)(1,0)[``\Id]{450}1a
 \putmorphism(360,180)(1,0)[\phantom{F(B)}` `b]{560}1a

\putmorphism(0,200)(0,-1)[\phantom{Y_2}``m]{400}1l
\putmorphism(450,200)(0,-1)[\phantom{Y_2}``]{400}1l
\putmorphism(480,200)(0,-1)[\phantom{Y_2}``n']{400}0l
\putmorphism(900,200)(0,-1)[\phantom{Y_2}``q']{400}1r
\put(150,0){\fbox{$\lambda_2$}}
\put(620,-20){\fbox{$\crta\beta'$}}

\putmorphism(0,-160)(1,0)[``\Id]{450}1b
 \putmorphism(340,-160)(1,0)[\phantom{F(B)}` `b']{560}1b
\put(20,-100){$\Nearrow$}
\put(480,-100){$\Nearrow$}
\efig}
$$
Then there is a unique 3-cell 
$$
 \scalebox{0.86}{\bfig
  \putmorphism(0,200)(1,0)[``a\times b]{520}1a
 \putmorphism(420,200)(1,0)[\phantom{F(B)}` `\Id]{660}1a

\putmorphism(0,220)(0,-1)[\phantom{Y_2}``m\times n]{420}1l
\putmorphism(500,220)(0,-1)[\phantom{Y_2}``]{420}1l
\putmorphism(1080,220)(0,-1)[\phantom{Y_2}``p'\times q']{420}1r
\put(150,0){\fbox{$\crta\alpha\times\crta\beta$}}
\put(650,0){\fbox{$\rho_1\times\rho_2$}}
\put(20,-100){$\Nearrow$}
\put(520,-100){$\Nearrow$}
\putmorphism(0,-160)(1,0)[``a'\times b']{500}1b
 \putmorphism(420,-160)(1,0)[\phantom{F(B)}` `\Id]{650}1b
\efig}
\quad\stackrel{\Gamma}{\Rrightarrow}\quad
 \scalebox{0.86}{\bfig
  \putmorphism(0,200)(1,0)[``\Id]{580}1a
 \putmorphism(470,200)(1,0)[\phantom{F(B)}` `a\times b]{700}1a

\putmorphism(0,220)(0,-1)[\phantom{Y_2}``m\times n]{420}1l
\putmorphism(580,220)(0,-1)[\phantom{Y_2}``]{420}1l
\putmorphism(1170,220)(0,-1)[\phantom{Y_2}``p'\times q']{420}1r
\put(150,0){\fbox{$\lambda_1\times\lambda_2$}}
\put(760,0){\fbox{$\crta\alpha'\times\crta\beta'$}}
\put(20,-100){$\Nearrow$}
\put(620,-100){$\Nearrow$}
\putmorphism(0,-160)(1,0)[``\Id]{570}1b
 \putmorphism(470,-160)(1,0)[\phantom{F(B)}` .`a'\times b']{700}1b
\efig}
$$
\end{prop}

\begin{cor} \colabel{middle prism}
Given equivalence 2-cells $\alpha, \beta, \alpha', \beta'$ and 3-cells 
$$
 \scalebox{0.86}{\bfig
  \putmorphism(0,200)(1,0)[``\Id]{450}1a
 \putmorphism(340,200)(1,0)[\phantom{F(B)}` `a]{600}1a

\putmorphism(0,220)(0,-1)[\phantom{Y_2}``m']{400}1l
\putmorphism(450,220)(0,-1)[\phantom{Y_2}``m]{400}1r
\putmorphism(900,220)(0,-1)[\phantom{Y_2}``p]{400}1r
\put(130,0){\fbox{$\lambda_1$}}
\put(600,0){\fbox{$\alpha$}}
\putmorphism(0,-150)(1,0)[``\Id]{450}1b
 \putmorphism(340,-150)(1,0)[\phantom{F(B)}` `a']{580}1b
\put(310,90){$\Swarrow$}
\put(760,90){$\Swarrow$}
\efig}
\quad\stackrel{\tilde P_1}{\Rrightarrow}\quad
 \scalebox{0.86}{\bfig
  \putmorphism(0,200)(1,0)[``a]{450}1a
 \putmorphism(350,200)(1,0)[\phantom{F(B)}` `\Id]{550}1a

\putmorphism(0,220)(0,-1)[\phantom{Y_2}``m']{400}1l
\putmorphism(440,220)(0,-1)[\phantom{Y_2}``p']{400}1r
\putmorphism(900,220)(0,-1)[\phantom{Y_2}``p]{400}1r
\put(120,0){\fbox{$\alpha'$}}
\put(590,20){\fbox{$\rho_1$}}
\put(310,90){$\Swarrow$}
\put(780,90){$\Swarrow$}
\putmorphism(0,-150)(1,0)[``a']{450}1b
 \putmorphism(350,-150)(1,0)[\phantom{F(B)}` `\Id]{580}1b
\efig}
$$
and
$$
 \scalebox{0.86}{\bfig
  \putmorphism(0,200)(1,0)[``\Id]{450}1a
 \putmorphism(340,200)(1,0)[\phantom{F(B)}` `b]{600}1a

\putmorphism(0,220)(0,-1)[\phantom{Y_2}``n']{400}1l
\putmorphism(450,220)(0,-1)[\phantom{Y_2}``n]{400}1r
\putmorphism(900,220)(0,-1)[\phantom{Y_2}``q]{400}1r
\put(130,0){\fbox{$\lambda_2$}}
\put(600,0){\fbox{$\beta$}}
\putmorphism(0,-150)(1,0)[``\Id]{450}1b
 \putmorphism(340,-150)(1,0)[\phantom{F(B)}` `b']{580}1b
\put(310,90){$\Swarrow$}
\put(760,90){$\Swarrow$}
\efig}
\quad\stackrel{\tilde P_2}{\Rrightarrow}\quad
 \scalebox{0.86}{\bfig
  \putmorphism(0,200)(1,0)[``b]{450}1a
 \putmorphism(350,200)(1,0)[\phantom{F(B)}` `\Id]{550}1a

\putmorphism(0,220)(0,-1)[\phantom{Y_2}``n']{400}1l
\putmorphism(440,220)(0,-1)[\phantom{Y_2}``q']{400}1r
\putmorphism(900,220)(0,-1)[\phantom{Y_2}``q]{400}1r
\put(120,0){\fbox{$\beta'$}}
\put(590,20){\fbox{$\rho_2$}}
\put(310,90){$\Swarrow$}
\put(780,90){$\Swarrow$}
\putmorphism(0,-150)(1,0)[``b']{450}1b
 \putmorphism(350,-150)(1,0)[\phantom{F(B)}` `\Id]{580}1b
\efig}
$$
with 2-cells $\lambda_i, \rho_i, i=1,2$ as in \prref{middle prism}. 
Then there is a unique 3-cell 
$$
 \scalebox{0.86}{\bfig
  \putmorphism(0,200)(1,0)[``\Id]{580}1a
 \putmorphism(470,200)(1,0)[\phantom{F(B)}` `a\times b]{700}1a

\putmorphism(0,220)(0,-1)[\phantom{Y_2}``m\times n]{420}1l
\putmorphism(580,220)(0,-1)[\phantom{Y_2}``]{420}1l
\putmorphism(1150,220)(0,-1)[\phantom{Y_2}``p'\times q']{420}1r
\put(50,0){\fbox{$\lambda_1\times\lambda_2$}}
\put(650,0){\fbox{$\alpha\times\beta$}}
\put(450,80){$\Swarrow$}
\put(1000,80){$\Swarrow$}
\putmorphism(0,-160)(1,0)[``\Id]{570}1b
 \putmorphism(470,-160)(1,0)[\phantom{F(B)}` `a'\times b']{700}1b
\efig}
\quad\stackrel{\Gamma'}{\Rrightarrow}\quad
 \scalebox{0.86}{\bfig
  \putmorphism(0,200)(1,0)[``a\times b]{570}1a
 \putmorphism(470,200)(1,0)[\phantom{F(B)}` `\Id]{650}1a

\putmorphism(0,220)(0,-1)[\phantom{Y_2}``m\times n]{420}1l
\putmorphism(550,220)(0,-1)[\phantom{Y_2}``]{420}1l
\putmorphism(1120,220)(0,-1)[\phantom{Y_2}``p'\times q']{420}1r
\put(40,0){\fbox{$\alpha'\times\beta'$}}
\put(580,0){\fbox{$\rho_1\times\rho_2$}}
\put(410,80){$\Swarrow$}
\put(980,80){$\Swarrow$}
\putmorphism(0,-160)(1,0)[``a'\times b']{550}1b
 \putmorphism(470,-160)(1,0)[\phantom{F(B)}` `\Id]{650}1b
\efig}
$$
{\em i.e.} $\Gamma': \displaystyle{\frac{[\Id\vert\alpha\times\beta]}{[\lambda_1\times\lambda_2\vert\Id_{a'\times b'}]}} \Rrightarrow 
\displaystyle{\frac{[\Id_{a\times b}\vert\rho_1\times\rho_2]}{[\alpha'\times\beta'\vert\Id]}}$. 
\end{cor}

\subsection{On 3-pullbacks on higher cells}

Analogously to \leref{lema}, for 3-pullbacks we have:

\begin{lma} \lelabel{pullb-lema}
Let $V$ be a 1-strict tricategory. Given two 3-pullbacks $M\times_S N$ and $P\times_S Q$ of cospans $M\stackrel{m}{\to} S\stackrel{n}{\leftarrow}N$ 
and  $P\stackrel{p}{\to} S\stackrel{q}{\leftarrow}Q$ in $V$, respectively, together with equivalence 2-cells
$\zeta_1:p_1t\Rightarrow q_1, \zeta_2:q_2\Rightarrow p_2t$ and $\theta_1:s_1\Rightarrow p'_1s, 
\theta_2:s_2\Rightarrow p'_2s$, as in the diagram below, and corresponding bijective 3-cells 
$$\Sigma_1: \threefrac{\Id_{n}\ot\zeta_2}{\omega\ot\Id_t}{\Id_m\ot\zeta_1} \Rrightarrow\sigma_1 
\quad\text{and}\quad
\Sigma_2: \threefrac{\Id_{q}\ot\theta_2}{\omega'\ot\Id_s}{\Id_p\ot\theta_1} \Rrightarrow\sigma_2.$$
Then one has:
\begin{enumerate}[a)]
\item given 1-cells $g:M\to P$ and $h:N\to Q$ and equivalence 2-cells \(\Fi_1: m \Rightarrow pg\) and \(\Fi_2: qh \Rightarrow n\),
then there exist a 1-cell $g\times_S h: M\times_S N\to P\times_S Q$, equivalence 2-cells $\omega_1:gp_1\Rightarrow p'_1(g\times_S h)$ 
and $\omega_2:hp_2\Rightarrow p'_2(g\times_S h)$, and an invertible 3-cell
$$\Sigma: \threefrac{\Fi_2 \otimes \Id_{p_2}}{\omega}{\Fi_1 \otimes \Id_{p_1}} \Rrightarrow
\threefrac{\Id_q \otimes \omega_2}{\omega' \otimes \Id_{g \times_S h}}{\Id_p \otimes \omega_1^{-1}};$$
\item additionally to the data from a), given 1-cell $f:T\to R$, 2-cells $\alpha:gq_1\Rightarrow s_1f$ and $\beta:hq_2\Rightarrow s_2f$, and a 3-cell
$$\upkappa_0: \threefrac{\sigma_1}{\Fi_1 \otimes \Id_{q_1}}{\Id_p \otimes \alpha} \Rrightarrow
\threefrac{\Fi_2 \otimes \Id_{q_2}}{\Id_q \otimes \beta}{\sigma_2 \otimes \Id_m},$$
then there is a 2-cell $\gamma:(g\times_S h)t\Rightarrow sf$ and invertible 3-cells  
$$\Xi^1_\gamma:\threefrac{\Id_g\ot\zeta_1}{\omega_1\ot\Id_t}{\Id_{p_1'}\ot\gamma}\Rrightarrow\frac{\alpha}{\theta_1\ot\Id_f}\qquad\text{and}\qquad
 \Xi^2_\gamma:\threefrac{\Id_h\ot\zeta_2}{\omega_2\ot\Id_t}{\Id_{p_2'}\ot\gamma}\Rrightarrow\frac{\beta}{\theta_2\ot\Id_f};$$
\item if $\alpha$ and $\beta$ in b) are equivalence 2-cells, then so is $\gamma$. 
\end{enumerate}
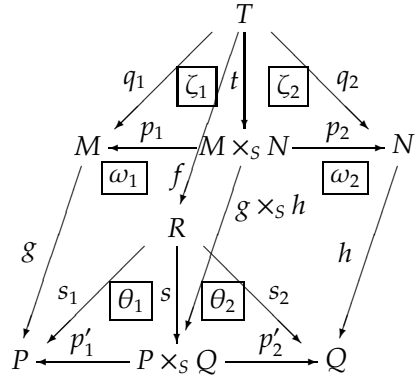
$$
 \scalebox{0.86}{\bfig
\putmorphism(400,300)(-1,-1)[``q_1]{500}1l
 \putmorphism(480,300)(1,-1)[\phantom{\u\C \bullet 1}` `q_2]{500}1r
 \putmorphism(-130,-200)(1,0)[M`\phantom{M\times_S  N} ` p_1]{600}{-1}a
 \putmorphism(490,-200)(1,0)[\phantom{\u\C \bullet 1}`N ` p_2]{550}1a
\putmorphism(450,300)(0,-1)[T`M\times_S  N`]{500}1l
\putmorphism(470,300)(0,-1)[`\phantom{M\times N}`t]{500}0l
\put(200,0){\fbox{$\zeta_1$}}
\put(540,0){\fbox{$\zeta_2$}}
\put(-80,-330){\fbox{$\omega_1$}}
\put(740,-330){\fbox{$\omega_2$}}
\putmorphism(450,300)(-1,-3)[``]{260}1l
\putmorphism(420,80)(-1,-3)[``f]{260}0l
\putmorphism(-130,-200)(-1,-3)[``g]{260}1l
\putmorphism(1040,-200)(-1,-3)[``h]{260}1l
\putmorphism(460,-200)(-1,-3)[``]{250}1r
\putmorphism(480,-50)(-1,-3)[``g\times_S h]{250}0r
\putmorphism(150,-500)(-1,-1)[``s_1]{500}1l
 \putmorphism(230,-500)(1,-1)[\phantom{\u\C \bullet 1}` `s_2]{500}1r
 \putmorphism(-380,-1000)(1,0)[P`\phantom{M\times_S N} ` p'_1]{600}{-1}a
 \putmorphism(240,-1000)(1,0)[\phantom{\u\C \bullet 1}`Q ` p'_2]{550}1a
\putmorphism(200,-500)(0,-1)[R`P\times_S Q`]{500}1l
\putmorphism(220,-500)(0,-1)[`\phantom{P\times_S Q}`s]{500}0l
\put(-50,-800){\fbox{$\theta_1$}}
\put(290,-800){\fbox{$\theta_2$}}
\efig}
$$
\end{lma}

\section{The double categories of matrices and spans in a 1-category} \selabel{Power}

In this Section we review the notion of extensivity used in \cite{CFP} to characterize biequivalence of bicategories of spans and matrices 
in a 1-category $\V$, 
and we show this characterization. In the last two Subsections we then complete the characterization of equivalence of the discretely 
internal and enriched categories using 2-categorical tools, that was announced in \cite{CFP}, but was not carried out this way, as the 
authors chose to do the proof using 1-categories.

\subsection{Review of extensivity} 

Let $\V$ be a category, $I$ a set, and $(X_i)_{i\in I}$ an $I$-indexed family of objects of $\V$. If $\V$ admits all
$I$-indexed coproducts, there is a functor 
\begin{equation} \eqlabel{ext-fun}
\Pi_{i\in I} (\V/X_i) \stackrel{\amalg}{\to} \V/(\amalg_{i\in I} X_i)
\end{equation}
mapping a family $(f_i : A_i \to X_i)_{i\in I}$ to $\amalg_{i\in I} f_i: \amalg_{i\in I} A_i\to \amalg_{i\in I} X_i$. 

Following \cite[Chap 2, 6.3]{BJ} and \cite[4.1]{CV} the definition below is introduced in \cite[Definition 2.1]{CFP}.

\begin{defn}
A category $\V$ is {\em extensive} if $\V$ admits all small coproducts and, for each small
family $(X_i)_{i\in I}$ of objects of $\V$, the functor in \equref{ext-fun} is an equivalence of categories.
\end{defn}

The right adjoint of $\amalg$, if it exists, has $i$-th component $\sigma_i^*: \V/(\amalg_{i\in I} X_i) \to \V/X_i$ mapping
$f: A \to \amalg_{i\in I} X_i$ to $\sigma_i^*(f): \sigma_i^*(A)\to X_i$ defined by the pullback
$$
 \scalebox{0.86}{\bfig
\putmorphism(0,500)(1,0)[\phantom{\sigma_i^*(A)}`\phantom{B}`]{500}1b
\putmorphism(500,500)(0,-1)[A`\amalg_{i\in I} X_i`f]{450}1r
\putmorphism(0,500)(0,-1)[\sigma_i^*(A)`X_i`\sigma_i^*(f)]{450}1l
\putmorphism(0,100)(1,0)[\phantom{X_i}`\phantom{\amalg_{i\in I} X_i}`\sigma_i]{500}1a
\efig}
$$
in $\V$. Then the adjunction has the form
$$
\bfig
\putmorphism(0,540)(1,0)[\phantom{\Pi_{i\in I} (\V/X_i)}`\phantom{\V/(\amalg_{i\in I} X_i)}`\amalg]{1000}1a
\putmorphism(0,500)(1,0)[\Pi_{i\in I} (\V/X_i)`\V/(\amalg_{i\in I} X_i)`]{1000}0a
\putmorphism(0,450)(1,0)[\phantom{\Pi_{i\in I} (\V/X_i)}`\phantom{\V/(\amalg_{i\in I} X_i)}` \perp]{1000}{-1}a
\putmorphism(0,470)(1,0)[\phantom{\Pi_{i\in I} (\V/X_i)}`\phantom{\V/(\amalg_{i\in I} X_i)}` \langle\sigma_i^*\rangle_{i\in I}]{1000}{0}b
\efig \vspace{-1,4cm}
$$

For a set $I$ and an object $X$ of $\V$ let $I\bullet X$ denote the $I$-fold copower of $X$ by $I$. If $\V$ has a
terminal object $1$, setting $X_i = 1$ for all $i\in I$ in \equref{ext-fun}, and observing that $\V/1\iso \V^I$, one obtains the functor
\begin{equation} \eqlabel{ext-I}
\V^I \stackrel{\amalg}{\longrightarrow} \V/(I\bullet 1)
\end{equation}
which is an adjunction 
\begin{equation} \eqlabel{adj-I}
\end{equation}
$$
\bfig
\putmorphism(0,540)(1,0)[\phantom{\V^I}`\phantom{\V/(I\bullet 1)}`\amalg]{1000}1a
\putmorphism(0,500)(1,0)[\V^I`\V/(I\bullet 1)`]{1000}0a
\putmorphism(0,450)(1,0)[\phantom{\V^I}`\phantom{\V/(I\bullet 1)}` \perp]{1000}{-1}a
\putmorphism(0,470)(1,0)[\phantom{\V^I}`\phantom{\V/(I\bullet 1)}` \langle i^*\rangle_{i\in I}]{1000}{0}b
\efig \vspace{-1,4cm}
$$
if for all $i\in I$ the category $\V$ admits all pullbacks along the coprojections $\crta i: 1 \to I \bullet 1$, due to \cite[Proposition 2.2]{CFP}. One sufficient 
condition for $\V$ to be extensive is:

\begin{prop} (\cite[Proposition 2.5]{CFP}, \cite[Proposition 4.1]{CLW}) \\
Let $\V$ be a category with small coproducts and a terminal object. If for every small set $I$, the functor \equref{ext-I} is an equivalence, 
then $\V$ is extensive. 
\end{prop}

\subsection{The bicategories of matrices and spans} \sslabel{bicats}

In \cite{BCSW} the bicategory $\V\x Mat$ of $\V$-matrices was introduced, and in \cite{Be} the bicategory $Span(\V)$ of spans over $\V$. It is immediate 
to see that monads 
in the former bicategory are $\V$-categories (categories enriched over $\V$), while the monads in the latter are categories internal in $\V$. As a matter 
of fact, in \cite[Section 5.4.3]{Be} categories internal in $\V$ are defined this way. For the purposes of examining the relation between enriched and internal 
categories in $\V$, the bicategory of spans over $\V$ is modified in \cite[Section 3.2]{CFP} so that the 0-cells are small sets, rather than objects of $\V$. 
This version of the bicategory is denoted by $Span_d(\V)$, here ``$d$'' stands for {\em discrete}, as monads in $Span_d(\V)$ are those internal categories 
whose object of objects is discrete. To shorten, these internal categories we will call throughout {\em discretely internal categories}, and spans in $Span_d(\V)$ 
we will call {\em discrete spans}. For the biequivalence of the mentioned two bicategories the authors have proved the following.

\begin{prop} \prlabel{bieq bicats} \cite[Proposition 3.2 and Theorem 3.3]{CFP} \\
Let $\V$ be a Cartesian closed category with finite limits and small coproducts. The following are equivalent: 
\begin{enumerate}
\item for every set $I$ the adjunction \equref{adj-I} is an adjoint equivalence;
\item $\V$ is extensive;
\item the oplax functor $Int: \V\x Mat \to Span_d(\V)$ is a biequivalence; 
\item the lax functor $En: Span_d(\V) \to \V\x Mat$ is a biequivalence.
\end{enumerate}
\end{prop}

The functors $Int$ and $En$ are obtained from the adjunction \equref{adj-I} by substituting the set $I$ from the latter by the set $I\times J$. 
Namely, 0-cells for both bicategories are sets $I,J...$ and the hom categories of the bicategories $\V\x Mat$ and $Span_d(\V)$ are given by 
$\V\x Mat(I, J) = \V^{I\times J}$ (1-cells are matrices of dimension $\vert I\vert\times\vert J\vert$ whose entries are objects of $\V$, 
and 2-cells are families of morphisms in $\V$ between the corresponding objects in the matrices) and 
$Span_d(\V)(I, J) = \V/((I \bullet 1) \times (J \bullet 1)) \iso \V/((I \times J)\bullet 1)$ (1-cells are spans of objects in $\V$ of the form: 
$I \bullet 1 \longleftarrow V \longrightarrow J \bullet 1$, 
and 2-cells are morphisms in $\V$ between such objects $V$ making two evident triangle diagrams between two spans commute), respectively. 
The above isomorphism of slice categories is assured by Cartesian closedness of $\V$. 

Concretely, we describe here the actions of $Int$ and $En$ on 1-cells, their actions on 2-cells can then be deduced easily, and they both 
are identities on 0-cells. $Int$ maps a matrix $(M(i;j))_{i\in I,j\in J}$ to the span 
$I \bullet 1 \longleftarrow \amalg_{i\in I,j\in J} M(i,j) \longrightarrow J \bullet 1$.  
Here the arrows to $I \bullet 1$ and $J \bullet 1$ are induced by the following composite with the domain $M(i,j)$ for fixed 
$i\in I, j\in J$: the unique morphism to 1 followed by the coprojections to $I \bullet 1$ and $J \bullet 1$, respectively. $En$ maps a 
span $I \bullet 1 \stackrel{v_I}{\longleftarrow} V \stackrel{v_J}{\longrightarrow} J \bullet 1$ to the matrix $En(V)$ whose $(i,j)$-th 
component is given by the pullback 
\begin{equation} \eqlabel{matrix of A}
\end{equation}
$$ \scalebox{0.86}{\bfig
\putmorphism(0,800)(1,0)[\phantom{En(V)(i,j)}`\phantom{B}`!]{800}1a
\putmorphism(800,800)(0,-1)[1`(I \bullet 1) \times (J \bullet 1)`\langle \crta i,\crta j\rangle]{500}1r
\putmorphism(0,800)(0,-1)[En(V)(i,j)`V`\iota_{i,j}]{500}1l
\putmorphism(0,300)(1,0)[\phantom{V}`\phantom{(I \bullet 1) \times (J \bullet 1)}`\langle v_I,v_J\rangle]{800}1a
\efig }
$$
for each $i\in I, j\in J$, where $\crta i$ and $\crta j$ denote the $i$-th and $j$-th coprojections, respectively.

\bigskip

Although monads in $\V\x Mat$ and $Span_d(\V)$ are $\V$-categories and discretely internal categories in $\V$, respectively, 
the morphisms of monads in these 
bicategories are not morphisms in $\V\x\Cat$, the category of $\V$-categories, and $\Cat_d(\V)$, the category of discretely internal categories in $\V$. Hence, as the authors comment, one can not use the biequivalence of bicategories 
$\V\x Mat$ and $Span_d(\V)$ to conclude the equivalence of categories $\V\x\Cat$ and $\Cat_d(\V)$. 
They prove the latter equivalence in another way avoiding two-dimensional category theory, although they comment that one could proceed by using 
``additional two-dimensional structure, such as that of a pseudo-double category''. For the purpose of our work, we will use the pseudo-double categories 
of $\V$-matrices and ``discrete spans in $\V$''.

\subsection{The double categories of matrices and spans} \sslabel{double}


In \cite[Example 2.1]{FGK} the authors introduced a pseudo-double category $\Span(\V)$ of spans in $\V$ whose horizontal bicategory is precisely the bicategory 
$Span(\V)$. In \cite[Definition 2.4]{FGK} a pseudo-double category $\Mnd(\Dd)$ of monads in a pseudo-double category $\Dd$ is 
introduced, so that when $\Dd=\Span(\V)$, the vertical 1-cells in $\Mnd(\Span(\V))$ are morphisms of internal categories in $\V$ (see \cite[Example 2.6]{FGK}). 
This inspires us to define the pseudo-double category $\Span_d(\V)$ by modifying accordingly $\Span(\V)$, and to introduce a pseudo-double category 
$\V\x\Mat$ so to extend the biequivalence of bicategories from \prref{bieq bicats} to an equivalence of pseudo-double categories. 

We define a pseudo-double category $\Span_d(\V)$ as follows. Its 0-cells are small sets, for sets $I,J$ 1h-cells $I\to J$ are spans 
$I \bullet 1 \stackrel{a_1}{\longleftarrow} A \stackrel{a_2}{\longrightarrow} J \bullet 1$, while 1v-cells $I\to J$ are set maps between $I$ and $J$.  
For spans 
$I \bullet 1 \stackrel{a_1}{\longleftarrow} A \stackrel{a_2}{\longrightarrow} J \bullet 1$ and 
$K \bullet 1 \stackrel{b_1}{\longleftarrow} B \stackrel{b_2}{\longrightarrow} L \bullet 1$, 
and maps $u:I\to K$ and $v: J\to L$, 2-cells are given by the diagrams 
\begin{equation} \eqlabel{span square}
\end{equation}
$$
 \scalebox{0.86}{\bfig
\putmorphism(-50,500)(1,0)[\phantom{\sigma_i^*(A)}`\phantom{B}`a_1]{550}{-1}a
\putmorphism(500,500)(0,-1)[A`B`f]{400}1r
\putmorphism(0,500)(0,-1)[I\bullet 1`K\bullet 1`u\bullet 1]{400}1l
\putmorphism(-30,100)(1,0)[\phantom{\sigma_i^*(A)}`\phantom{B}`b_1]{550}{-1}a
\putmorphism(400,500)(1,0)[\phantom{\sigma_i^*(A)}`\phantom{B}`a_2]{550}1a
\putmorphism(400,100)(1,0)[\phantom{\sigma_i^*(A)}`\phantom{B}`b_2]{550}1a
\putmorphism(1000,500)(0,-1)[J\bullet 1`L\bullet 1`v\bullet 1]{400}1r
\efig}
$$

A pseudo-double category $\V\x\Mat$ we define as follows. Its 0-cells are small sets, for sets $I,J$ 1h-cells $I\to J$ are matrices of dimension 
$\vert I\vert\times\vert J\vert$ whose entries are objects of $\V$, 1v-cells $I\to K$ are maps of sets, and for matrices 
$(M(i,j))_{i\in I, j\in J}$ and $(N(k,l))_{k\in K, l\in L}$ and 1v-cells $u:I\to K$ and $v: J\to L$, 2-cells are given by the families $f$ 
of morphisms in $\V$ determined so that the following diagram commutes: 
\begin{equation} \eqlabel{matrix square}
\end{equation}
$$\scalebox{0.86}{\bfig
\putmorphism(0,800)(1,0)[\phantom{\amalg_{i\in I,j\in J} M(i,j)} `\phantom{I\bullet 1\times J\bullet 1}`]{1100}1a
\putmorphism(0,800)(0,-1)[\amalg_{i\in I,j\in J} M(i,j)`\amalg_{k\in K,l\in L} N(k,l)`f]{400}1l
\putmorphism(1100,800)(0,-1)[I\bullet 1\times J\bullet 1`K\bullet 1\times L\bullet 1.`]{400}0r
\putmorphism(1020,800)(0,-1)[``u\bullet 1\times v\bullet 1]{400}1r
\putmorphism(0,400)(1,0)[\phantom{\amalg_{i\in I,j\in J} M(i,j)}`\phantom{V}`]{850}1a
\efig} \vspace{-1cm}
$$
(the horizontal arrows above are induced by the terminal morphism followed by the map to the product induced by the two corresponding 
coprojections). 
This means that $f$ is given by a family of morphisms 
$$f_{i,j}: M(i,j)\to N(u(i),v(j))$$ 
in $\V$ for every $i\in I,j\in J$.

\begin{rem} \rmlabel{prod-coprod}
By the product property, the 2-cells \equref{span square} can equivalently be described by commutative squares
\equref{one-sided span}. On the other hand, when $\V$ is a Cartesian closed category, the functors $X\times-$ and $-\times X$ 
for $X\in\V$ are left adjoint functors. As such they preserve colimits, implying that there is a natural isomorphism 
$\phi_{I, J} \colon (I \bullet 1) \times (J \bullet 1) \iso (I \times J) \bullet 1$ in \(\V\). 
Then this implies that the squares \equref{one-sided span} can equivalently be described by commutative squares
\equref{coproduct square}.   \vspace{-1cm}
\begin{center} 
\begin{tabular}{p{5.8cm}p{0cm}p{5cm}}
\begin{equation} \eqlabel{one-sided span}
 \scalebox{0.86}{\bfig
\putmorphism(0,400)(1,0)[\phantom{A} `\phantom{(I \bullet 1) \times (J \bullet 1)}`\langle a_1,a_2\rangle]{900}1a
\putmorphism(0,400)(0,-1)[A`B`f]{400}1l
\putmorphism(900,400)(0,-1)[(I \bullet 1) \times (J \bullet 1)`(K \bullet 1) \times (L \bullet 1)`]{400}0r
\putmorphism(720,400)(0,-1)[``u\bullet 1 \times v\bullet 1]{400}1r
\putmorphism(0,0)(1,0)[\phantom{A}`\phantom{V}`\langle b_1,b_2\rangle]{600}1a
\efig }
\end{equation}
&  & 
\begin{equation} \eqlabel{coproduct square}
 \scalebox{0.86}{\bfig
\putmorphism(0,400)(1,0)[\phantom{A} `\phantom{I\bullet 1}`a]{500}1a
\putmorphism(0,400)(0,-1)[A`B`f]{400}1l
\putmorphism(640,400)(0,-1)[(I\times J)\bullet 1`(K\times L)\bullet 1.`]{400}0r
\putmorphism(540,400)(0,-1)[``(u\times v)\bullet 1]{400}1r
\putmorphism(0,0)(1,0)[\phantom{A}`\phantom{V}`b]{420}1a
\efig} 
\end{equation}
\end{tabular}
\end{center} 
\end{rem}

\begin{rem} \rmlabel{other matrices}

Adding natural isomorphism $\phi_{I, J}$ from the above Remark to the square \equref{matrix square} yields that 
the 2-cells of $\V\x\Mat$ can equivalently be described as the squares: 
\begin{equation} 
\end{equation}
$$ \scalebox{0.86}{\bfig
\putmorphism(0,800)(1,0)[\phantom{\amalg_{i\in I,j\in J} M(i,j)} `\phantom{(I\times J)\bullet 1}`]{1100}1a
\putmorphism(0,800)(0,-1)[\amalg_{i\in I,j\in J} M(i,j)`\amalg_{k\in K,l\in L} N(k,l)`f]{400}1l
\putmorphism(1100,800)(0,-1)[(I\times J)\bullet 1`(K\times L)\bullet 1`]{400}0r
\putmorphism(1020,800)(0,-1)[``(u\times v)\bullet 1]{400}1r
\putmorphism(0,400)(1,0)[\phantom{\amalg_{i\in I,j\in J} M(i,j)}`\phantom{V}`]{900}1a
\efig} \vspace{-1cm}
$$
In this case the horizontal arrows are induced by the unique morphism to 1 followed by the corresponding coprojections. 
\end{rem}

\medskip

In the pseudo-double categories $\Span_d(\V)$ and $\V\x\Mat$ 0- and 1h-cells are the same as 0- and 1-cells in the bicategroies 
$Span_d(\V)$ and $\V\x Mat$, respectively, and 1v-cells in both pseudo-double categories are the same. It is immediate to see that 
the lax functor $En$ is compatible with 2-cells (a morphism between pullbacks $En(A)(i,j)$ to $En(B)(k,l)$ exists for $(k,l)=(u(i),v(j))$ precisely because \equref{span square} commutes). 
In the 1v-direction $En$ is a strict functor, thus we get a lax double functor $\En: \Span_d(\V)\to\V\x\Mat$.

Conversely, starting from a 2-cell \equref{matrix square}, 
it clearly induces a 2-cell of the form \equref{span square}. Hence, we have that if the oplax functor $Int: \V\x Mat \to Span_d(\V)$ is a biequivalence, 
then the oplax double functor $\Int: \V\x\Mat \to \Span_d(\V)$ is a double equivalence. The converse is also true (by restriction to identity 1v-cells), 
so we have: $Int$ is a biequivalence if and only if $\Int$ is a double equivalence. 
A similar statement holds for $En$ and $\En$. In view of \prref{bieq bicats}  we have:

\begin{prop} \prlabel{double equiv} 
Let $\V$ be a Cartesian closed category with finite limits and small coproducts. The following are equivalent: 
\begin{enumerate}
\item for every set $I$ the adjunction \equref{adj-I} is an adjoint equivalence;
\item $\V$ is extensive;
\item the oplax double functor $\Int: \V\x\Mat \to \Span_d(\V)$ is a double equivalence; 
\item the lax double functor $\En: \Span_d(\V) \to \V\x \Mat$ is a double equivalence.
\end{enumerate}
\end{prop}




\subsection{Monads in the double categories of matrices and spans} \sslabel{monads in double cats}

Our pseudo-double category $\Span_d(\V)$ is a pseudo-double subcategory of $\Span(\V)$ (our 1h-cells are specific 1h-cells of the latter), and we have that 
the vertical 1-cells in $\Mnd(\Span_d(\V))$ are morphisms of discretely internal categories in $\V$ (see the beginning of \ssref{bicats}). As for our 
pseudo-double category $\V\x\Mat$, we find the following. A monad is given by a matrix $M=(M(i,j))_{i,j\in I}$ (1-endocell over the 0-cell $I$) 
and 2-cells which are given by families of morphisms $\mu^M_{i,k}:\amalg_{j\in I}M(j,k)\times M(i,j)\to M(i,k)$ and 
$\eta^M_{i,k}:\Ii\to M(i,k)$, for every $i,k\in I$, both given by a commutative diagram \equref{matrix square} where the maps $u$ and $v$ are identities 
on $I$, satisfying associativity and unity laws. (Here $\Ii$ is the unit matrix, where $1$ is the terminal object and $0$ the initial object in $\V$.) 
Thus, a monad in $\V\x\Mat$ is precisely a monad in the bicategory $\V\x Mat$. Given another monad $(N=(N(k,l))_{k\in K, l\in L}, \mu^N,\eta^N)$, a 1v-cell 
in $\Mnd(\V\x\Mat)$ 
between $M$ and $N$ is given by a set map $\omega:I\to I$ and a commutative square \equref{matrix square} in which $u$ and $v$ are equal to $\omega$ 
(this square comes down to a morphism $f_{i,j}: M(i,j)\to N(\omega(i),\omega(j))$ in $\V$ for every $i,j\in I$) and which satisfy: 
$$\scalebox{0.86}{\bfig
\putmorphism(0,800)(1,0)[\phantom{\amalg_{j\in I}M(j,k)\times M(i,j)}`\phantom{M(i,k)}`\mu^M_{i,k}]{1400}1a
\putmorphism(0,800)(0,-1)[\amalg_{j\in I}M(j,k)\times M(i,j)`\amalg_{j\in I}N(\omega(j),\omega(k))\times N(\omega(i),\omega(j))`
   ]{400}0l
\putmorphism(110,800)(0,-1)[``  \amalg_{j\in I}f_{j,k}\times f_{i,j}]{400}1l
\putmorphism(1400,800)(0,-1)[M(i,k)`N(\omega(i),\omega(k))`f_{i,k}]{400}1r
\putmorphism(0,400)(1,0)[\phantom{{\amalg_{j\in I}N(\omega(j),\omega(k))\times N(\omega(i),\omega(j))}} ` 
   \phantom{N(\omega(i),\omega(k))}`\mu^N_{i,k}]{1400}1a
\efig} \vspace{-1,4cm}
$$
and $f_{i,j}\comp\eta^M_{i,j}=\eta^N_{i,j}$. Thus vertical 1-cells in $\V\x\Mat$ clearly correspond to enriched functors in $\V$.

\bigskip

It is immediately seen that a biequivalence of two 2-categories induces a biequivalence of the 2-categories of their monads. 
Passing from a strict to a weak context (from 2-categories to pseudo-double categories), one needs to be more careful with the technical details in the proof, 
but the analogous result is directly proved. Observe that the Definition 2.4 in \cite{FGK} of the double category $\Mnd(\Dd)$ of monads in a double category 
$\Dd$ is such that if $\Dd$ is a pseudo-double category, then $\Mnd(\Dd)$ is a pseudo-double category, too.

\begin{prop} \prlabel{equiv monads}
If pseudo-double categories $\Dd$ and $\Ee$ are double equivalent, then their pseudo-double categories of monads $\Mnd(\Dd)$ and $\Mnd(\Ee)$ 
are also double equivalent. 
\end{prop}

\begin{proof}
The proof is straightforward, we only type the diagrams for the relevant structure 2-cells. 
$$
\scalebox{0.86}{
\bfig
\putmorphism(-580,500)(0,-1)[\phantom{Y_2}``=]{450}1l
\putmorphism(1520,500)(0,-1)[\phantom{Y_2}``=]{450}1r
\putmorphism(-550,500)(1,0)[F(X)`\phantom{HF(A)}`=]{550}1a
\putmorphism(-550,50)(1,0)[F(X)`\phantom{HF(A)}`=]{550}1a
\put(-450,260){\fbox{$F^X$}}
\put(1060,260){\fbox{$(F^X)^{-1}$}}

 \putmorphism(980,500)(1,0)[\phantom{F(A)}`F(X) `=]{550}1a
 \putmorphism(1000,50)(1,0)[\phantom{F(A)}`F(X) `=]{550}1a

\putmorphism(-20,500)(1,0)[F(X)`F(X)`F(PP)]{1000}1a
 \putmorphism(-20,50)(1,0)[F(X)`F(X)`F(P)]{1000}1a

\putmorphism(-80,500)(0,-1)[\phantom{Y_2}``]{450}1r
\putmorphism(-100,500)(0,-1)[\phantom{Y_2}``F(u)]{450}0r

\putmorphism(1020,500)(0,-1)[\phantom{Y_2}``]{450}1l
\putmorphism(1050,500)(0,-1)[\phantom{Y_2}``F(u)]{450}0l

\put(300,310){\fbox{$F(\mu_P)$}}
\put(300,710){\fbox{$F_{P,P}^{-1}$}}

\putmorphism(-150,900)(1,0)[F(X)`F(X) `F(P)]{600}1a
 \putmorphism(450,900)(1,0)[\phantom{F(X)}` F(X) `F(P)]{600}1a

\putmorphism(-80,950)(0,-1)[\phantom{Y_2}``=]{450}1l
\putmorphism(1020,950)(0,-1)[\phantom{Y_3}``=]{450}1r

\efig}
\quad
\scalebox{0.86}{
\bfig
\putmorphism(450,900)(1,0)[`F(X) `=]{460}1a
\putmorphism(380,920)(0,-1)[F(X)` `=]{420}1l
\putmorphism(920,920)(0,-1)[\phantom{Y_2}``=]{420}1r
\put(550,700){\fbox{$F_X^{-1}$} }

 \putmorphism(-150,500)(1,0)[F(X)`F(X) `=]{500}1a
 \putmorphism(450,500)(1,0)[` `F(id_X)]{380}1a
\putmorphism(-180,500)(0,-1)[\phantom{Y_2}` `=]{450}1l
\put(-120,240){\fbox{$F^X$}}
\putmorphism(-150,50)(1,0)[F(X)`F(X) `=]{500}1a
\putmorphism(360,500)(0,-1)[\phantom{Y_2}` `F(u)]{450}1l
\put(380,250){\fbox{$F(\eta_P)$} } 
\putmorphism(940,500)(0,-1)[F(X)`F(X)`F(u)]{450}1l
\putmorphism(470,50)(1,0)[``F(P)]{360}1a
 \putmorphism(1050,500)(1,0)[` `=]{280}1a
 \putmorphism(1050,50)(1,0)[` `=]{280}1a
\putmorphism(1440,500)(0,-1)[F(X)`F(X)`=]{450}1r
\put(1020,260){\fbox{$(F^X)^{-1}$}}
\efig}
$$

$$
\scalebox{0.86}{
\bfig
\putmorphism(-580,500)(0,-1)[\phantom{Y_2}``=]{450}1l
\putmorphism(1520,500)(0,-1)[\phantom{Y_2}``=]{450}1r
\putmorphism(-550,500)(1,0)[F(X)`\phantom{HF(A)}`=]{550}1a
\putmorphism(-550,50)(1,0)[F(X)`\phantom{HF(A)}`=]{550}1a
\put(-450,260){\fbox{$F^X$}}
\put(1060,260){\fbox{$(F^Y)^{-1}$}}

 \putmorphism(980,500)(1,0)[\phantom{F(A)}`F(Y) `=]{550}1a
 \putmorphism(1000,50)(1,0)[\phantom{F(A)}`F(Y) `=]{550}1a

\putmorphism(-20,500)(1,0)[F(X)`F(Y)`F(QH)]{1000}1a
 \putmorphism(-20,50)(1,0)[F(X)`F(Y)`F(HP)]{1000}1a

\putmorphism(-80,500)(0,-1)[\phantom{Y_2}``]{450}1r
\putmorphism(-100,500)(0,-1)[\phantom{Y_2}``F(id)]{450}0r

\putmorphism(1020,500)(0,-1)[\phantom{Y_2}``]{450}1l
\putmorphism(1050,500)(0,-1)[\phantom{Y_2}``F(id)]{450}0l

\put(300,310){\fbox{$F(\Phi)$}}
\put(300,710){\fbox{$F_{P,P}^{-1}$}}

\putmorphism(-150,900)(1,0)[F(X)`F(Y) `F(H)]{600}1a
 \putmorphism(450,900)(1,0)[\phantom{F(X)}` F(Y) `F(Q)]{600}1a

\putmorphism(-80,950)(0,-1)[\phantom{Y_2}``=]{450}1l
\putmorphism(1020,950)(0,-1)[\phantom{Y_3}``=]{450}1r
\putmorphism(-80,50)(0,-1)[\phantom{Y_2}``=]{450}1r

\putmorphism(1020,50)(0,-1)[\phantom{Y_2}``=]{450}1l

\putmorphism(-110,-400)(1,0)[F(X)`F(X) `F(P)]{560}1a
 \putmorphism(450,-400)(1,0)[\phantom{F(X)}` F(Y) `F(H)]{560}1a

\put(300,-170){\fbox{$F_{P,P}$}}
\efig}
\quad
\scalebox{0.86}{
\bfig
 \putmorphism(380,500)(1,0)[F(X)` `F(P)]{400}1a
\putmorphism(350,500)(0,-1)[\phantom{Y_2}` `F(u)]{450}1l
\put(490,270){\fbox{$F(\crta u)$} } 
\putmorphism(890,500)(0,-1)[F(X)`F(X)`F(u)]{450}1r
\putmorphism(380,50)(1,0)[F(X)``F(P)]{400}1a
\efig}
$$
\end{proof}

From the above said we obtain:

\begin{cor} \colabel{deduce int-en}
In the conditions of \prref{double equiv} the categories $\V\x\Cat$ of $\V$-enriched categories and $\Cat_d(\V)$ of discretely internal categories in $\V$ 
are equivalent. 
\end{cor}

\bigskip

\section{$(1\times 2)$-category of spans in a tricategory} \selabel{3-spans}

The term {\em $(1\times 2)$-category} is due to \cite{Shul} (see the top of page 2). It is a 1-category internal to a tricategory. 

In this Section we define a structure that would be a tricategorical version of the pseudo-double category 
$\Span_d(\V)$ from \ssref{double} for a category $\V$. We will build a 
category internal in a certain 1-strict tricategory $T$, which to shorten we will call a 
$(1\times 2)$-category $\S$. Then we will also have the horizontal tricategory $\HH(\S)$. We will do this   
gradually, by introducing first two bicategories $C_0$ (the bicategory of objects) and $C_1$ (the bicategory of morphisms). 
Thus, the 0-cells of the tricategory $T$ will be bicategories, then 1-cells will be some kind of 2-dimensional functors between them, 
and we have that $T$ is 1-strict. The precise notion of a category internal in a 1-strict tricategory we introduced in \cite{Fem}.

\subsection{Bicategories $C_0$ and $C_1$ in $\S$}  \sslabel{3-spans}

Assume that $V$ is a 1-strict tricategory with a terminal object 1, small 3-coproducts and 3-pullbacks. 
Let $C_0=\Cat_2$, the 2-category of {\em small} categories. 
Before we define $C_1$ we observe the following. 

\bigskip

In \rmref{prod-coprod} we saw that morphisms between spans in a 1-category $\V$ could be described in two equivalent ways. 
Already turning to dimension 2 (if $\V$ were a 2-category) requires  
the involvement of two 2-cells, in the case of \equref{span square}, respectively of one 2-cell, in the case of the squares 
\equref{one-sided span} and \equref{coproduct square}. 
These two approaches yield two ways of defining 
1-cells in the bicategory of spans in $V$.  
However, the former is more suitable for defining their composition. At the end of the Subsection 
we will show that the two ways of defining these 1-cells are equivalent.

\bigskip

Let $\C$ and $\D$ be small categories (0-cells in $C_0$). 
The 0-cells of $C_1$ are {\em spans in $V$}, given by 1-cells $\u\C \bullet 1 \longleftarrow A \longrightarrow \u\D \bullet 1$ in $V$. Here 
$\u\C:=\Ob\C$ denotes the set of objects of $\C$ and $\u\C \bullet 1$ the $\u\C$-fold copower of 1. 

Given two spans $\u\C \bullet 1 \stackrel{a_1}{\longleftarrow} A \stackrel{a_2}{\longrightarrow} \u\D \bullet 1$ and 
$\u\E \bullet 1 \stackrel{b_1}{\longleftarrow} B \stackrel{b_2}{\longrightarrow} \u\HH \bullet 1$ and 
two functors $F:\C\to\E$ and $G:\D\to\HH$ (1-cells in $C_0$) a 1-cell in $C_1$ is given by a 1-cell $f:A\to B$ and two equivalence 2-cells 
$\alpha:(F\bullet 1)a_1\Rightarrow b_1f$ and $\beta:(G\bullet 1)a_2\rightarrow b_2f$ in $V$: 
\begin{equation}\eqlabel{1-cells in C_1}
\end{equation}
$$
\scalebox{0.86}{\bfig
\putmorphism(-150,250)(1,0)[\u\C \bullet 1`A`a_1]{600}{-1}a
 \putmorphism(400,250)(1,0)[\phantom{\u\C \bullet 1}`\u\D \bullet 1 `a_2]{650}1a

 \putmorphism(-150,-200)(1,0)[\u\E\bullet 1` B` b_1]{600}{-1}a
 \putmorphism(400,-200)(1,0)[\phantom{\u\C \bullet 1}`\u\HH\bullet 1. ` b_2]{650}1a

\putmorphism(-180,250)(0,-1)[\phantom{Y_2}``F\bullet 1]{450}1l
\putmorphism(450,250)(0,-1)[\phantom{Y_2}``]{450}1r
\putmorphism(300,250)(0,-1)[\phantom{Y_2}``f]{450}0r
\putmorphism(1020,250)(0,-1)[\phantom{Y_2}``G\bullet 1]{450}1r
\put(100,20){\fbox{$\alpha$}}
\put(650,20){\fbox{$\beta$}}
\efig}
$$
Given another such 1-cell in $C_1$ with the same 0-cells: 
$$
\scalebox{0.86}{\bfig
\putmorphism(-150,250)(1,0)[\u\C \bullet 1`A`a_1]{600}{-1}a
 \putmorphism(400,250)(1,0)[\phantom{\u\C \bullet 1}`\u\D\bullet 1 `a_2]{650}1a

 \putmorphism(-150,-200)(1,0)[\u\E\bullet 1` B` b_1]{600}{-1}a
 \putmorphism(400,-200)(1,0)[\phantom{\u\C \bullet 1}`\u\HH\bullet 1 ` b_2]{650}1a

\putmorphism(-180,250)(0,-1)[\phantom{Y_2}``F'\bullet 1]{450}1l
\putmorphism(450,250)(0,-1)[\phantom{Y_2}``]{450}1r
\putmorphism(300,250)(0,-1)[\phantom{Y_2}``g]{450}0r
\putmorphism(1020,250)(0,-1)[\phantom{Y_2}``G'\bullet 1]{450}1r
\put(100,20){\fbox{$\gamma$}}
\put(650,20){\fbox{$\delta$}}
\efig}
$$
and natural transformations $\lambda:F\Rightarrow F':\C\to \E$ and $\rho: G\Rightarrow G':\D\to \HH$ (2-cells in $C_0$), 
a 2-cell in $C_1$ between them is given by a 2-cell $\xi: f\Rightarrow g$ and two 3-cells 
\begin{equation} \eqlabel{3-cells}
\Sigma: \frac{\alpha}{[\xi\vert\Id_{b_1}]} \Rrightarrow \frac{[\Id_{a_1}\vert\lambda\bullet 1]}{\gamma} 
\qquad\text{and}\qquad
\Omega: \frac{\beta}{[\xi\vert\Id_{b_2}]} \Rrightarrow \frac{[\Id_{a_1}\vert\rho\bullet 1]}{\delta} 
\end{equation}
in $V$, which are to be considered in the transversal direction, perpendicular to the parallel planes of the 1-cells $(\alpha,f,\beta)$ and 
$(\delta, g,\gamma)$. 
The 2-cells in $C_1$ we think as transversal prisms whose bases are 1-cells in $C_1$. 

Composition of 1-cells in $C_1$, which contains the horizontal composition of 2-cells in $\HH(\S)$, the underlying horizontal 
tricategory of $\S$, 
is defined by the 3-pullback. Given 1-cells $(\alpha,f,\beta)$ and $(\gamma, g,\delta)$ as below:
\begin{equation} \eqlabel{comps of 2-spans}
\end{equation}
$$
\scalebox{0.86}{\bfig
\putmorphism(180,140)(2,1)[``a_1]{50}{-1}l
 \putmorphism(400,300)(2,-1)[\phantom{\u\C \bullet 1}`\u\D \bullet 1 `a_2]{650}1r

 \putmorphism(180,-300)(2,1)[` ` b_1]{50}{-1}r
 \putmorphism(400,-100)(2,-1)[\phantom{\u\C \bullet 1}`\u\HH\bullet 1 ` ]{650}0a
 \putmorphism(400,-140)(2,-1)[\phantom{\u\C \bullet 1}` ` b_2]{630}1l

\putmorphism(-80,0)(0,-1)[\u\C \bullet 1`\u\E\bullet 1`F\bullet 1]{370}1l
\putmorphism(450,250)(0,-1)[A`B`f]{450}1r
\putmorphism(1020,0)(0,-1)[\phantom{Y_2}``G\bullet 1]{450}1r
\put(100,-90){\fbox{$\alpha$}}
\put(650,-90){\fbox{$\beta$}}
\putmorphism(1280,120)(2,1)[``]{50}{-1}l
\putmorphism(1280,140)(2,1)[``a'_1]{50}{0}l
 \putmorphism(1500,280)(2,-1)[\phantom{\u\C \bullet 1}`\u\J \bullet 1 `a'_2]{650}1r

 \putmorphism(1280,-300)(2,1)[` ` ]{50}{-1}r
 \putmorphism(1380,-250)(1,0)[` ` b'_1]{50}0r
 \putmorphism(1500,-100)(2,-1)[\phantom{\u\C \bullet 1}`\u\K\bullet 1 ` ]{650}0a
 \putmorphism(1500,-140)(2,-1)[\phantom{\u\C \bullet 1}` ` b'_2]{630}1l

\putmorphism(1550,250)(0,-1)[A'`B'`g]{450}1r
\putmorphism(2140,-20)(0,-1)[\phantom{Y_2}``H\bullet 1]{420}1r
\put(1300,-80){\fbox{$\gamma$}}
\put(1750,-100){\fbox{$\delta$}}
\efig}
$$
Taking for $\sigma$ in \deref{3-pb} the composite equivalence 2-cell on the right below, there are: a 1-cell $h$, 
equivalence 2-cells $\zeta_1, \zeta_2$ 
and an isomorphism 3-cell $\Sigma$ in $V$: 
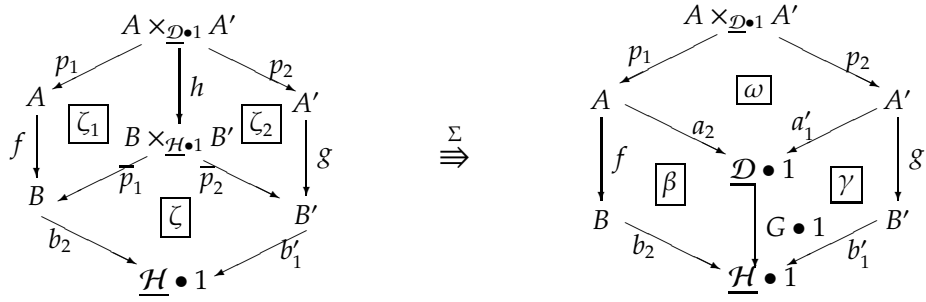
$$
\scalebox{0.86}{\bfig
\putmorphism(120,360)(2,1)[``p_1]{50}{-1}l
 \putmorphism(420,510)(2,-1)[\phantom{\u\C \bullet 1}` `p_2]{600}1r
 \putmorphism(140,-60)(2,1)[` ` \crta p_1]{50}{-1}r
 \putmorphism(400,90)(2,-1)[\phantom{\u\C \bullet 1}` ` \crta p_2]{560}1l
\putmorphism(-80,260)(0,-1)[A`B`f]{370}1l
\putmorphism(450,510)(0,-1)[A\times_{\u\D\bullet 1}A'` B\times_{\u\HH\bullet 1}B' ` h]{420}1r
\putmorphism(920,240)(0,-1)[A'`B'`g]{420}1r
\put(40,130){\fbox{$\zeta_1$}}
\put(680,130){\fbox{$\zeta_2$}}
 \putmorphism(-200,-120)(2,-1)[\phantom{\u\C \bullet 1}`\u\HH\bullet 1 ` b_2]{630}1l
 \putmorphism(720,-340)(2,1)[` ` ]{60}{-1}r
 \putmorphism(840,-250)(1,0)[` ` b'_1]{50}0r
\put(380,-220){\fbox{$\zeta$}}
\efig}
\quad\qquad\stackrel{\Sigma}{\Rrightarrow}
\scalebox{0.86}{\bfig
 \putmorphism(970,550)(2,-1)[A\times_{\u\D\bullet 1}A'` `p_2]{650}1r

\putmorphism(660,380)(2,1)[``]{80}{-1}l
\putmorphism(680,400)(2,1)[``p_1]{50}{0}l
\put(950,260){\fbox{$\omega$}}

 \putmorphism(400,300)(2,-1)[\phantom{\u\C \bullet 1}`\u\D \bullet 1 `a_2]{650}1r

 \putmorphism(400,-100)(2,-1)[\phantom{\u\C \bullet 1}`\u\HH\bullet 1 ` ]{650}0a
 \putmorphism(400,-140)(2,-1)[\phantom{\u\C \bullet 1}` ` b_2]{630}1l

\putmorphism(450,250)(0,-1)[A`B`f]{450}1r
\putmorphism(1020,0)(0,-1)[\phantom{Y_2}``G\bullet 1]{450}1r
\put(650,-90){\fbox{$\beta$}}
\putmorphism(1280,120)(2,1)[``]{50}{-1}l
\putmorphism(1280,140)(2,1)[``a'_1]{50}{0}l

 \putmorphism(1280,-300)(2,1)[` ` ]{50}{-1}r
 \putmorphism(1380,-250)(1,0)[` ` b'_1]{50}0r

\putmorphism(1550,250)(0,-1)[A'`B'`g]{450}1r
\put(1300,-80){\fbox{$\gamma$}}
\efig}
$$
For the desired 1-cell $A\times_{\D\bullet 1} A'\to B\times_{\HH\bullet 1}B'$ in $V$ we take this $h$, 
and for the desired pair of equivalence 2-cells $(\alpha',\beta')$ in $V$ we set the horizontal juxtapositions: $\alpha'=(\alpha\vert\zeta_1)$ and 
$\beta'=(\zeta_2\vert\delta)$ of 2-cells in $V$. 
Properly speaking, $\alpha'=\frac{[\zeta_1\vert\Id_{b_1}]}{[\Id_{p_1}\vert\alpha]}$ and $\beta'=\frac{[\zeta_2\vert\Id_{b_2}]}{[\Id_{p_2}\vert\delta]}$, 
where $p_1,p_2$ are the projections of the 3-pullback $A\times_{\D\bullet 1}A'$. 

Since the composition of spans of the underlying 1-category of $V$ (and consequently of 1-cells in $\HH(\S)$) is not strictly associative, 
the same holds for the 2-cells of $\HH(\S)$ and 1-cells of $C_1$. This is why $C_1$ is a bicategory, and not a 2-category. 

Vertical composition of horizontal 2-cells in $\S$ (and 2-cells in $\HH(\S)$) - as in \equref{1-cells in C_1} and below it - is given in the obvious way: 
$(\frac{[\alpha\vert\Id_{F'\bullet 1}]}{[\Id_f\vert\gamma]}, gf, 
\frac{[\beta\vert\Id_{G'\bullet 1}]}{[\Id_f\vert\delta]})$. It is not strictly associative: interchange law of $V$ must be used, as well as 
the following isomorphisms between 2-cells in $V$: $\Id_{X\ot Y}\iso[\Id_X\vert\Id_Y]$ and the one for the associativity of the horizontal composition in $V$. 
These isomorphisms 
can be expressed in terms of 3-cells in $\HH(\S)$. 

The interchange law in $C_1$ is expressed in terms of 3-cells of $V$ and it holds strictly, as one hoped. 

The vertical composition of 2-cells in $C_1$ (and of 3-cells in $\HH(\S)$ and horizontal 3-cells in $\S$) 
is given by obvious vertical juxtaposition of prisms. In the case of their horizontal 
composition 
we proceed as follows. Suppose we are given two composable horizontal 1-cells in $\S$ as in \equref{comps of 2-spans} 
(horizontally composable 2-cells in $\HH(\S)$) and another such a pair 
with the same 1-cells in $V$, so that only the vertically denoted 1-cells and 2-cells in $V$ are different: $F'\bullet 1, \alpha', f', \beta', 
G'\bullet 1, \gamma', g', \delta', H'\bullet 1$. Suppose that we are given two horizontally composable 2-cells in $C_1$, each of which we present by a pair of prisms. 
Concretely, for simplicity reasons, the pair of 3-cells as in \equref{3-cells} we will identify with prisms which we will write in the present case as: 
$\A:\alpha\Rrightarrow\alpha', \B:\beta\Rrightarrow\beta', \Gamma:\gamma\Rrightarrow\gamma'$ and $\Delta:\delta\Rrightarrow\delta'$. We think $\A,\B,\Gamma$ 
and $\Delta$ as transversal prisms going from the base squares $\alpha, \beta, \gamma, \delta$ in the back towards the base squares $\alpha', \beta', \gamma', 
\delta',$ in the front. On the vertically transversal planes of these four prisms are the following 2-cells in $V$: 
$\lambda\bullet 1:F\bullet 1\Rightarrow F'\bullet 1, 
\xi:f\Rightarrow f', \rho:G\bullet 1\Rightarrow G'\bullet 1, \xi':g\Rightarrow g', \rho':H\bullet 1\Rightarrow H'\bullet 1$, with the obvious meanings. 

We know that the horizontal composition of the bases (back) 2-cells $(\alpha, f, \beta)$ and $(\gamma, g, \delta)$ is given by 
$((\alpha\vert\zeta_1), h, (\zeta_2\vert\delta))$ and we have the isomorphism 3-cell $\Sigma$ (as we explained below \equref{comps of 2-spans}). 
Analogously, at the front we have $((\alpha'\vert\zeta'_1), h', (\zeta'_2\vert\delta'))$ and we have an analogous isomorphism 3-cell $\Sigma'$:  
$$
\scalebox{0.86}{\bfig
\putmorphism(120,360)(2,1)[``p_1]{50}{-1}l
 \putmorphism(420,510)(2,-1)[\phantom{\u\C \bullet 1}` `p_2]{600}1r

 \putmorphism(140,-60)(2,1)[` ` \crta p_1]{50}{-1}r
 \putmorphism(400,90)(2,-1)[\phantom{\u\C \bullet 1}` ` \crta p_2]{560}1l

\putmorphism(-80,260)(0,-1)[A`B`f']{370}1l
\putmorphism(450,510)(0,-1)[A\times_{\u\D\bullet 1}A'` B\times_{\u\HH\bullet 1}B' ` h']{420}1r
\putmorphism(940,240)(0,-1)[A'`B'`g']{420}1r
\put(40,130){\fbox{$\zeta'_1$}}
\put(680,130){\fbox{$\zeta'_2$}}
 \putmorphism(-200,-120)(2,-1)[\phantom{\u\C \bullet 1}`\u\HH\bullet 1 ` b_2]{630}1l
 \putmorphism(720,-340)(2,1)[` ` ]{60}{-1}r
 \putmorphism(840,-250)(1,0)[` ` b'_1]{50}0r
\put(380,-220){\fbox{$\zeta$}}
\efig}
\quad\qquad\stackrel{\Sigma'}{\Rrightarrow}
\scalebox{0.86}{\bfig
 \putmorphism(970,550)(2,-1)[A\times_{\u\D\bullet 1}A'` `p_2]{650}1r

\putmorphism(660,380)(2,1)[``]{80}{-1}l
\putmorphism(680,400)(2,1)[``p_1]{50}{0}l
\put(950,260){\fbox{$\omega$}}

 \putmorphism(400,300)(2,-1)[\phantom{\u\C \bullet 1}`\u\D \bullet 1 `a_2]{650}1r

 \putmorphism(400,-100)(2,-1)[\phantom{\u\C \bullet 1}`\u\HH\bullet 1 ` ]{650}0a
 \putmorphism(400,-140)(2,-1)[\phantom{\u\C \bullet 1}` ` b_2]{630}1l

\putmorphism(450,250)(0,-1)[A`B`f']{450}1r
\putmorphism(1020,0)(0,-1)[\phantom{Y_2}``G'\bullet 1]{450}1r
\put(650,-90){\fbox{$\beta'$}}
\putmorphism(1280,120)(2,1)[``]{50}{-1}l
\putmorphism(1280,140)(2,1)[``a'_1]{50}{0}l

 \putmorphism(1280,-300)(2,1)[` ` ]{50}{-1}r
 \putmorphism(1380,-250)(1,0)[` ` b'_1]{50}0r

\putmorphism(1550,250)(0,-1)[A'`B'.`g']{450}1r
\put(1300,-80){\fbox{$\gamma'$}}
\efig}
$$
Observe the following (transversal) composition of 3-cells: 
$$
\scalebox{0.7}{
\bfig
\putmorphism(-550,260)(0,-1)[\u\C\bullet 1`\u\E\bullet 1`F\bullet 1]{370}1l
\putmorphism(-460,250)(1,0)[``a_1]{340}{-1}a
\putmorphism(-460,-110)(1,0)[``b_1]{340}{-1}a
\putmorphism(120,360)(2,1)[``p_1]{50}{-1}l
 \putmorphism(420,510)(2,-1)[\phantom{\u\C \bullet 1}` `p_2]{560}1r
 \putmorphism(140,-60)(2,1)[` ` \crta p_1]{50}{-1}r
 \putmorphism(400,90)(2,-1)[\phantom{\u\C \bullet 1}` ` \crta p_2]{560}1l
\putmorphism(-80,260)(0,-1)[A`B`f]{370}1l
\putmorphism(450,510)(0,-1)[A\times_{\u\D\bullet 1}A'` B\times_{\u\HH\bullet 1}B' ` h]{420}1r
\putmorphism(900,240)(0,-1)[A'`B'`g]{420}1r
\put(40,130){\fbox{$\zeta_1$}}
\put(680,130){\fbox{$\zeta_2$}}
 \putmorphism(-200,-120)(2,-1)[\phantom{\u\C \bullet 1}`\u\HH\bullet 1 ` b_2]{630}1l
 \putmorphism(720,-340)(2,1)[` ` ]{60}{-1}r
 \putmorphism(840,-250)(1,0)[` ` b'_1]{50}0r
\put(380,-220){\fbox{$\zeta$}}
\putmorphism(1310,240)(0,-1)[\phantom{Y_2}``H\bullet 1]{420}1r
\putmorphism(940,230)(1,0)[`\u\J\bullet 1`a'_2]{400}1a
\putmorphism(940,-180)(1,0)[`\u\K\bullet 1`b'_2]{400}1a
\put(-400,80){\fbox{$\alpha$}}
\put(1050,60){\fbox{$\delta$}}
\efig}
\stackrel{(\Id_\alpha\vert\Sigma\vert\Id_\delta)}{\Rrightarrow}
\scalebox{0.7}{
\bfig
 \putmorphism(-550,260)(0,-1)[\u\C\bullet 1`\u\E\bullet 1`F\bullet 1]{370}1l
\putmorphism(-460,250)(1,0)[``a_1]{340}{-1}a
\putmorphism(-460,-110)(1,0)[``b_1]{340}{-1}a
\putmorphism(450,550)(2,-1)[A\times_{\u\D\bullet 1}A'` `p_2]{650}1r
 \putmorphism(140,380)(2,1)[``p_1]{50}{-1}l 
\put(430,260){\fbox{$\omega$}}
 \putmorphism(-120,300)(2,-1)[\phantom{\u\C \bullet 1}`\u\D \bullet 1 `a_2]{650}1r
 \putmorphism(-120,-100)(2,-1)[\phantom{\u\C \bullet 1}`\u\HH\bullet 1 ` ]{650}0a
 \putmorphism(-120,-140)(2,-1)[\phantom{\u\C \bullet 1}` ` b_2]{630}1l
\putmorphism(-70,250)(0,-1)[A`B`f]{450}1r 
\putmorphism(500,0)(0,-1)[\phantom{Y_2}``G\bullet 1]{450}1r
\put(130,-90){\fbox{$\beta$}}
\putmorphism(760,120)(2,1)[``]{50}{-1}l
\putmorphism(760,140)(2,1)[``a'_1]{50}{0}l
 \putmorphism(760,-300)(2,1)[` ` ]{50}{-1}r
 \putmorphism(860,-250)(1,0)[` ` b'_1]{50}0r
\putmorphism(1030,250)(0,-1)[A'`B'`g]{450}1l
\put(760,-80){\fbox{$\gamma$}}
\putmorphism(1420,240)(0,-1)[\phantom{Y_2}``H\bullet 1]{420}1r
\putmorphism(1060,250)(1,0)[`\u\J\bullet 1`a'_2]{400}1a
\putmorphism(1060,-190)(1,0)[`\u\K\bullet 1`b'_2]{400}1a
\put(-400,80){\fbox{$\alpha$}}
\put(1150,60){\fbox{$\delta$}}
\efig}
$$

$\hspace{9cm} \Downarrow\displaystyle{(\A\vert\frac{\Id_\omega}{(\B\vert\Gamma)}}\vert\Delta) $

$$\hspace{-0,2cm}
\scalebox{0.7}{
\bfig
\putmorphism(-550,260)(0,-1)[\u\C\bullet 1`\u\E\bullet 1`F'\bullet 1]{370}1l
\putmorphism(-460,250)(1,0)[``a_1]{340}{-1}a
\putmorphism(-460,-110)(1,0)[``b_1]{340}{-1}a
\putmorphism(120,360)(2,1)[``p_1]{50}{-1}l
 \putmorphism(420,510)(2,-1)[\phantom{\u\C \bullet 1}` `p_2]{600}1r
 \putmorphism(140,-60)(2,1)[` ` \crta p_1]{50}{-1}r
 \putmorphism(400,90)(2,-1)[\phantom{\u\C \bullet 1}` ` \crta p_2]{560}1l
\putmorphism(-80,260)(0,-1)[A`B`f']{370}1l
\putmorphism(450,510)(0,-1)[A\times_{\u\D\bullet 1}A'` B\times_{\u\HH\bullet 1}B' ` h']{420}1r
\putmorphism(940,240)(0,-1)[A'`B'`g']{420}1r
\put(40,130){\fbox{$\zeta'_1$}}
\put(680,130){\fbox{$\zeta'_2$}}
 \putmorphism(-200,-120)(2,-1)[\phantom{\u\C \bullet 1}`\u\HH\bullet 1 ` b_2]{630}1l
 \putmorphism(720,-340)(2,1)[` ` ]{60}{-1}r
 \putmorphism(840,-250)(1,0)[` ` b'_1]{50}0r
\put(380,-220){\fbox{$\zeta$}}
\putmorphism(1310,240)(0,-1)[\phantom{Y_2}``H'\bullet 1]{420}1r
\putmorphism(940,230)(1,0)[`\u\J\bullet 1`a'_2]{400}1a
\putmorphism(940,-180)(1,0)[`\u\K\bullet 1`b'_2]{400}1a
\put(-400,80){\fbox{$\alpha'$}}
\put(1070,60){\fbox{$\delta'$}}
\efig}
\stackrel{(\Id_{\alpha'}\vert\crta{\Sigma'}\vert\Id_{\delta'})}{\Lleftarrow}
\scalebox{0.7}{
\bfig
 \putmorphism(-550,260)(0,-1)[\u\C\bullet 1`\u\E\bullet 1`F'\bullet 1]{370}1l
\putmorphism(-460,250)(1,0)[``a_1]{340}{-1}a
\putmorphism(-460,-110)(1,0)[``b_1]{340}{-1}a
\putmorphism(450,550)(2,-1)[A\times_{\u\D\bullet 1}A'` `p_2]{650}1r
 \putmorphism(140,380)(2,1)[``p_1]{50}{-1}l 

\put(430,260){\fbox{$\omega$}}
 \putmorphism(-120,300)(2,-1)[\phantom{\u\C \bullet 1}`\u\D \bullet 1 `a_2]{650}1r
 \putmorphism(-120,-100)(2,-1)[\phantom{\u\C \bullet 1}`\u\HH\bullet 1 ` ]{650}0a
 \putmorphism(-120,-140)(2,-1)[\phantom{\u\C \bullet 1}` ` b_2]{630}1l
\putmorphism(-70,250)(0,-1)[A`B`f']{450}1r 
\putmorphism(500,0)(0,-1)[\phantom{Y_2}``G'\bullet 1]{450}1r
\put(130,-90){\fbox{$\beta'$}}
\putmorphism(760,120)(2,1)[``]{50}{-1}l
\putmorphism(760,140)(2,1)[``a'_1]{50}{0}l
 \putmorphism(760,-300)(2,1)[` ` ]{50}{-1}r
 \putmorphism(860,-250)(1,0)[` ` b'_1]{50}0r
\putmorphism(1030,250)(0,-1)[A'`B'`g']{450}1l
\put(700,-80){\fbox{$\gamma'$}}
\putmorphism(1420,240)(0,-1)[\phantom{Y_2}``H'\bullet 1]{420}1r
\putmorphism(1060,250)(1,0)[`\u\J\bullet 1`a'_2]{400}1a
\putmorphism(1060,-190)(1,0)[`\u\K\bullet 1`b'_2]{400}1a
\put(-400,80){\fbox{$\alpha'$}}
\put(1170,60){\fbox{$\delta'$}}
\efig}
$$
Here $\crta{\Sigma'}$ denotes the inverse of $\Sigma'$. 
Compose the domain and codomain 2-cells above vertically with $(\Id_{b_1}\vert\crta\zeta\vert\Id_{b'_2})$ from below, 
where $\crta\zeta$ is a quasi-inverse of $\zeta$, and compose the above composition 3-cell with the according 3-cell induced by the 
identity on $\crta\zeta$. From the result one obtains  
the wanted 3-cell $((\alpha\vert\zeta_1), h, (\zeta_2\vert\delta)) \Rrightarrow ((\alpha'\vert\zeta'_1), h', (\zeta'_2\vert\delta'))$. 

The transversal composition of 3-cells $(\Sigma, \xi, \Omega)$ and $(\Sigma', \xi', \Omega')$ in $\HH(\S)$ (as in \equref{3-cells}) 
is given by the transversal composition of the 3-cell components and the vertical composition of the 2-cells: $\frac{\xi}{\xi'}$. 

\medskip

This finishes the construction of 0-cells $C_0$ and $C_1$ of $T$, the first step to define 
a $(1\times 2)$-category $\S$ of spans in $V$. We finish this Subsection with the promised result. 

\medskip

\begin{prop} \prlabel{2-cell-prop}
Provided the existence of 3-products, 
a 1-cell \equref{1-cells in C_1} in the bicategory $C_1$ of spans in $V$ can equivalently be described by an equivalence 2-cell 
$\gamma$ in $V$: 
$$\scalebox{0.86}{\bfig
\putmorphism(0,400)(1,0)[\phantom{A} `\phantom{(\u{\C}\times \u{\D})\bullet 1}`\langle a_1,a_2\rangle]{900}1a
\putmorphism(0,400)(0,-1)[A`B`f]{400}1l
\putmorphism(950,400)(0,-1)[\u{\C}\bullet 1\times \u{\D}\bullet 1`\u{\E}\bullet 1\times \u{\HH}\bullet 1.`]{400}0r
\putmorphism(770,400)(0,-1)[``F\bullet 1\times G\bullet 1]{400}1r
\putmorphism(0,0)(1,0)[\phantom{A}`\phantom{V}`\langle b_1,b_2\rangle]{700}1a
\put(320,240){\fbox{$\gamma$}}
\efig }
$$
A 2-cell \equref{3-cells} in $C_1$ can equivalently be described by a 3-cell \\
$\Gamma: \displaystyle{\frac{[\Id\vert\gamma]}{[\xi\vert\Id_{\langle b_1,b_2\rangle}]}} \Rrightarrow 
\displaystyle{\frac{[\Id_{\langle a_1,a_2\rangle}\vert\lambda\bullet 1\times\rho\bullet 1]}{[\gamma'\vert\Id]}}$. 
\end{prop}

\begin{proof}
This is \leref{lema} c) and d). The converse holds by \rmref{converse for 3-pr}. 
The claim for 2-cells in $C_1$ follows analogously by \prref{cutting 3-cells}. 
\end{proof}

\subsection{The 1-cells $u,s,t,c$ in $T$}

For the sake of saving space we just state that $c:C_1\times_{C_0}C_1\to C_1$ and $u:C_0\to C_1$ are pseudofunctors and the source and target are strict 2-functors. 
On 0-cells $c:C_1\times_{C_0}C_1\to C_1$ is given by the 3-pullback and on 1- and 2-cells we defined it in the previous Section. 
The pseudofunctor $u:C_0\to C_1$ we define using initial object and cells (see \ssref{terminal}). 
We leave the construction of pseudonatural transformations $a^*, l^*$ and $r^*$ and modifications $\pi^*, \mu^*, \lambda^*, \rho^*$ and 
$\epsilon^*$ from \cite[Definition 6.2]{Fem} to the reader. Accordingly, we take for $T$, the tricategory from the beginning of the Section, to be 
the tricategory $\Bicat_3$ of bicategories, pseudofunctors, pseudonatural transformations and modification.

\section{$(1\times 2)$-category of matrices in a tricategory} \selabel{3-matrices}

In this Section we are going to construct another category internal in the tricategory $T=\Bicat_3$, 
which is a tricategorification of the pseudo-double category $\V\x\Mat$ of matrices from \ssref{double}. We will denote this 
$(1\times 2)$-category by $\M$. 

Suppose that $V$ is a 1-strict tricategory with a terminal object 1, 3-products  
and small 3-coproducts. Let $D_0=\Cat_2$ be the 2-category of {\em small} categories, as in the case of spans in $V$. Before 
defining the bicategory $D_1$, we first set up some notational conventions. 

For a small category $\C$ and $A\in\C$, 1-cells $1\stackrel{\sigma_A}{\to}\u{\C}\bullet 1$ will stand for coprojections. 
Given two small categories $\C$ and $\D$, we will write for short $\crta\C\times\crta\D:=\u\C\bullet 1\times\u\D\bullet 1$. 
A 1-cell to the 3-product $\crta\C\times\crta\D$ induced by coprojections $\sigma_A$ and $\sigma_B$ with $B\in\D$ we will denote by 
$1\stackrel{\langle A, B\rangle}{\longrightarrow}\crta\C\times\crta\D$. Accordingly, for small categories 
$\C,\D,\E, \HH$ and functors $F:\C\to\E$ and $G:\D\to\HH$ we will denote $\crta F\times\crta G:=F\bullet 1\times G\bullet 1$, 
and $\crta\lambda\times\crta\rho:=\lambda\bullet 1\times \rho\bullet 1$ for natural transformations $\lambda$ and $\rho$.

\subsection{Defining the bicategory $D_1$} \sslabel{3-matrices}

We now define the bicategory $D_1$ of matrices in $V$ as follows. 

\medskip

Given small categories $\C$ and $\D$, a {\em matrix in $V$} is a 1-cell 
$\amalg_{\substack{A\in\u\C \\ B\in\u\D}}M(A,B) \stackrel{m}{\to} \crta\C\times\crta\D$ in $V$, 
where the domain 0-cell we think as a matrix of objects $M(A,B)$ in $V$ indexed by the objects of the categories $\C$ and $\D$. 
The 1-cell $m$ is the unique 1-cell to the terminal object followed by 
$1\stackrel{\langle A, B\rangle}{\longrightarrow}\crta\C\times\crta\D$. 
Matrices in $V$ are 0-cells of $D_1$. 

Given two matrices $\amalg_{\substack{A\in\u\C \\ B\in\u\D}}M(A,B) \stackrel{m}{\to} \crta\C\times\crta\D$ and 
$\amalg_{\substack{A'\in\u\E \\ B'\in\u\HH}} N(A',B') \stackrel{n}{\to} \crta\E\times\crta\HH$, for small categories 
$\C,\D,\E, \HH$, and given two functors $F:\C\to\E$ and $G:\D\to\HH$ (1-cells in $D_0$) a {\em morphism of matrices} and a 1-cell of $D_1$ 
is a square: 
$$\scalebox{0.86}{\bfig
\putmorphism(0,400)(1,0)[\amalg_{\substack{A\in\u\C \\ B\in\u\D}}M(A,B) `\phantom{\crta\C\times\crta\D}`m]{950}1a
\putmorphism(-30,440)(0,-1)[\phantom{\amalg_{\substack{A\in\u\C \\ B\in\u\D}}M(A,B)}`
   \phantom{\amalg_{\substack{A'\in\u\E \\ B'\in\u\HH}} N(A',B')} ` f]{400}1r
\putmorphism(950,400)(0,-1)[\crta\C\times\crta\D`\crta\E\times\crta\HH`]{400}0r
\putmorphism(860,400)(0,-1)[``\crta F\times\crta G]{400}1r
\putmorphism(20,0)(1,0)[\amalg_{\substack{A'\in\u\E \\ B'\in\u\HH}} N(A',B')`\phantom{\crta\C\times\crta\D}`n]{920}1a
\put(440,200){\fbox{$\nu$}}
\efig }
$$
which consists of a family of pairs of equivalence 2-cells in $V$: 
\begin{equation} \eqlabel{nu}
\scalebox{0.86}{\bfig
\putmorphism(0,200)(1,0)[\phantom{\big(M(A,B)\big)_{\substack{A\in\u\C \\ B\in\u\D}}} `\phantom{\crta\E\times\crta\HH}`
   {{[A,B]}}]{1100}1a
\putmorphism(130,200)(0,-1)[M(A,B)`N(A',B')`   f_{A,B}]{400}1l
\putmorphism(1100,200)(0,-1)[\crta\C\times\crta\D`\crta\E\times\crta\HH`]{400}0r
\putmorphism(1020,200)(0,-1)[``\crta F\times\crta G]{400}1r
\putmorphism(30,-200)(1,0)[\phantom{\big(M(A,B)\big)_{\stackrel{A\in\u\C}{B\in\D}}}`\phantom{\crta\E\times\crta\HH}`
    {{[A',B']}}]{1060}1a
\put(440,40){\fbox{$\nu_{A,B}$}}
\efig }
\qquad
\scalebox{0.86}{\bfig
 \putmorphism(0,200)(1,0)[1`\crta\C\times\crta\D ` \langle A,B\rangle]{750}1a
 \putmorphism(-60,210)(1,-1)[\phantom{A}` `\langle A',B'\rangle]{520}1l
 \putmorphism(50,230)(1,-1)[\phantom{A}`\crta\E\times\crta\HH `]{550}0l
\putmorphism(650,200)(0,-1)[``\crta F\times\crta G]{500}1r
\put(320,-20){\fbox{$\chi^{F,G}_{A,B}$}}
\efig}
\end{equation}
where $\chi^{F,G}_{A,B}$ is induced, according to \leref{lema} d), by equivalence 2-cells 
\begin{equation} \eqlabel{one-sided chi's}
\scalebox{0.86}{\bfig
 \putmorphism(0,200)(1,0)[1`\u{\C}\bullet 1 ` \sigma_A]{550}1a
 \putmorphism(-20,210)(1,-1)[\phantom{A}` `\sigma_{A'}]{520}1l
 \putmorphism(0,230)(1,-1)[\phantom{A}`\u{\E}\bullet 1 `]{550}0l
\putmorphism(550,200)(0,-1)[``F\bullet 1]{500}1r
\put(320,-20){\fbox{$\chi^F_A$}}
\efig}
\qquad
\scalebox{0.86}{\bfig
 \putmorphism(0,200)(1,0)[1` \u{\D}\bullet 1 ` \sigma_B]{550}1a
 \putmorphism(-20,210)(1,-1)[\phantom{A}` `\sigma_{B'}]{520}1l
 \putmorphism(0,230)(1,-1)[\phantom{A}`\u{\HH}\bullet 1, `]{550}0l
\putmorphism(550,200)(0,-1)[``G\bullet 1]{500}1r
\put(320,-20){\fbox{$\chi^G_B$}}
\efig}
\end{equation}
indexed by pairs $(A,B)\in\u\C\times\u\D$, whenever there exist isomorphisms $F(A)\iso A'$ and $G(B)\iso B'$. 
(These isomorphisms condition also the existence of the 1-cells $f_{A,B}$.) 
In the case that $F$ and $G$ are identities we consider $A=A', B=B'$. 
The 1-cells $[A,B]$ are the unique morphism to 1 followed by $1\stackrel{\langle A, B\rangle}{\longrightarrow}\crta\C\times\crta\D$.

If $\chi^{F^{-1},G^{-1}}_{A,B}$ is an equivalence 2-cell corresponding to 
$\crta{F^{-1}}\times\crta{G^{-1}}$, then for quasi-inverses of $\chi^{F,G}_{A,B}$ one has: $(\chi^{F,G}_{A,B})^{-1}\iso
[\chi^{F^{-1},G^{-1}}_{A',B'}\vert\Id_{\crta F\times\crta G}]$.

Finally, given another morphism of matrices among the same matrices: 
$$\scalebox{0.86}{\bfig
\putmorphism(0,400)(1,0)[\amalg_{\substack{A\in\u\C \\ B\in\u\D}}M(A,B) `\phantom{\crta\C\times\crta\D}`m]{950}1a
\putmorphism(-30,440)(0,-1)[\phantom{\amalg_{\substack{A\in\u\C \\ B\in\u\D}}M(A,B)}`
   \phantom{\amalg_{\substack{A'\in\u\E \\ B'\in\u\HH}} N(A',B')} ` f']{400}1r
\putmorphism(950,400)(0,-1)[\crta\C\times\crta\D`\crta\E\times\crta\HH`]{400}0r
\putmorphism(860,400)(0,-1)[``\crta{F'}\times\crta{G'}]{400}1r
\putmorphism(20,0)(1,0)[\amalg_{\substack{A'\in\u\E \\ B'\in\u\HH}} N(A',B')`\phantom{\crta\C\times\crta\D}`n]{920}1a
\put(440,200){\fbox{$\nu'$}}
\efig }
$$
and natural transformations $\lambda:F\Rightarrow F':\C\to \E$ and $\rho: G\Rightarrow G':\D\to \HH$ (2-cells in $D_0$), 
a 2-cell in $D_1$ between $\nu$ and $\nu'$ is given by a 2-cell $\xi: f\Rightarrow f'$ and a 
transversal prism whose bases are vertical squares of the two 1-cells, where this prism consists of a family 
of prisms, {\em i.e.} 3-cells in $V$: 
\begin{equation} \eqlabel{3-cell in matrix}
\Sigma_{A,B}: \frac{\nu_{A,B}}{[\xi_{A,B}\vert\Id_{[A',B']}]} \Rrightarrow \frac{[\Id_{[A,B]}\vert\crta\lambda\times\crta\rho]}{\nu'_{A,B}} 
\end{equation}
for every $(A,B)\in\u\C\times\u\D$. 
Vertical composition of 2-cells in $D_1$ is clear: it is induced by vertical concatenation of the corresponding 
2-cells $\nu_{A,B}$.

\subsection{Composition of 1-cells in $D_1$}

The composition of matrices in $V$ is analogous to that of matrices in a 1-category, namely in the bicategory $\V\x Mat$ from \ssref{bicats}. 
Given matrices $\big(M(A,B)\big)_{\substack{A\in\u\C \\ B\in\u\D}}$ and $\big(N(B,C)\big)_{\substack{B\in\u\D \\ C\in\u\J}}$ their composition is 
given by the matrix $\big(\amalg_{B\in\u\D}M(A,B)\times N(B,C)\big)_{\substack{A\in\u\C \\ C\in\u\J}}$ and the corresponding 1-cell 
to $\crta\C\times\crta\J$. This defines the composition of 1h-cells 
in the $(1\times 2)$-category $\M$ of matrices. We now define the composition of 1-cells in the bicategory $D_1$.

\medskip

Given 1-cells $\nu$ and $\nu'$ with their respective families of 2-cells, we consider the following diagrams: 
$$\scalebox{0.86}{\bfig \hspace{-1cm}
\putmorphism(0,200)(1,0)[\phantom{\big(M(A,B)\big)_{\substack{A\in\u\C \\ B\in\u\D}}} `\phantom{\crta\C\times\crta\D}`
   {{[A,B]}}]{940}1a
\putmorphism(910,200)(1,0)[\phantom{\crta\C\times\crta\D}`\u{\C}\bullet 1 ` p_1]{640}1a
\putmorphism(40,660)(4,-1)[``]{1760}1r
\putmorphism(40,700)(4,-1)[``{{[A]}}]{1760}0r
\putmorphism(130,580)(0,-1)[M(A,B)``   =]{400}1l
\putmorphism(130,200)(0,-1)[M(A,B)`M'(A',B')`   f_{A,B}]{400}1l
\putmorphism(130,-200)(0,-1)[`M'(A',B')` =  ]{400}1l
\put(210,380){\fbox{$\tilde\theta_1$}}
\putmorphism(640,-520)(4,1)[``]{620}1r
\putmorphism(340,-740)(4,1)[``{{[A']}}]{1760}0l
\put(210,-420){\fbox{$\tilde\theta'_1$}}
\putmorphism(930,200)(0,-1)[\crta\C\times\crta\D`\crta\E\times\crta\HH`]{400}0r
\putmorphism(1020,200)(0,-1)[``\crta F\times\crta G]{400}1l
\putmorphism(1550,200)(0,-1)[``F\bullet 1]{400}1l
\putmorphism(30,-200)(1,0)[\phantom{\big(M(A,B)\big)_{\stackrel{A\in\u\C}{B\in\D}}}`\phantom{\crta\C\times\crta\D}`
    {{[A',B']}}]{910}1a
\putmorphism(910,-200)(1,0)[\phantom{\crta\C\times\crta\D}`\u{\E}\bullet 1 `p'_1]{640}1a
\put(160,20){\fbox{$\nu_{A,B}$}}
\put(1080,10){\fbox{$\omega_1$}}
\efig} 
\hspace{-2cm}
\scalebox{0.86}{\bfig
\putmorphism(0,200)(1,0)[\phantom{\big(M(A,B)\big)_{\substack{A\in\u\C \\ B\in\u\D}}} `\phantom{\crta\C\times\crta\D}`
   {{[B,C]}}]{940}1a
\putmorphism(910,200)(1,0)[\phantom{\crta\C\times\crta\D}`\u{\J}\bullet 1 ` p_2]{640}1a
\putmorphism(420,400)(4,1)[`` {{[C]}}  ]{660}0l
\putmorphism(420,360)(4,1)[``   ]{660}1l
\putmorphism(130,200)(0,-1)[N(B,C)`N'(B', C')`   g_{B,C}]{400}1l
\putmorphism(-150,-200)(4,-1)[``   ]{1820}1l
\putmorphism(-150,-240)(4,-1)[`` {{[C']}} ]{1820}0l
\put(1160,380){\fbox{$\tilde\theta_2$}}
\put(1160,-400){\fbox{$\tilde\theta'_2$}}
\putmorphism(930,200)(0,-1)[\crta\D\times\crta\J`\crta\HH\times\crta\K`]{400}0r
\putmorphism(1020,200)(0,-1)[``\crta G\times\crta H]{400}1l
\putmorphism(1550,600)(0,-1)[\u{\J}\bullet 1``=]{400}1l
\putmorphism(1550,200)(0,-1)[``H\bullet 1]{400}1l
\putmorphism(1550,-200)(0,-1)[`\u{\K}\bullet 1`=]{400}1l
\putmorphism(10,-200)(1,0)[\phantom{\big(M(A,B)\big)_{\stackrel{A\in\u\C}{B\in\D}}}`\phantom{\crta\C\times\crta\D}`
    {{[B',C']}}]{910}1a
\putmorphism(910,-200)(1,0)[\phantom{\crta\C\times\crta\D}`\u{\K}\bullet 1 `p'_2]{640}1a
\put(160,20){\fbox{$\nu'_{B,C}$}}
\put(1080,10){\fbox{$\omega_2$}}
\efig}
$$
where the 2-cells $\omega_1,\omega_2$ are the ones from \leref{lema} b) and $[A]=\sigma_A !$. 
From \leref{lema} a) we have an equivalence 2-cell $\theta_1:\sigma_A\Rightarrow p_1\sigma_{A,B}$, 
then let $\tilde\theta_1:=[\Id_{!}\vert\theta_1]$. Similarly we define $\tilde\theta'_1, \tilde\theta_2, \tilde\theta'_2$. 
To simplify the notations, 
let us denote the above composite equivalence 2-cells by: $\tilde\nu_{B,C}: (F\bullet 1)[A]\Rightarrow [A']f_{A,B}$ and 
$\tilde\nu'_{A,B}: (H\bullet 1)[C]\Rightarrow [C']g_{B,C}$. They induce an equivalence 2-cell $\tilde\nu_{A,B}\times\tilde\nu'_{B,C}$ 
in the middle of the diagram: 
\begin{equation} \eqlabel{nu-nu}
\scalebox{0.86}{\bfig
\putmorphism(270,330)(2,1)[``]{380}1l
\putmorphism(30,230)(2,1)[``!]{1100}0l
\put(750,410){\fbox{$\gamma$}}
\putmorphism(860,890)(0,-1)[`1 `]{300}0a
\putmorphism(800,650)(2,-1)[``]{850}1r
\putmorphism(560,780)(2,-1)[``\langle A, C\rangle]{1160}0r
\putmorphism(220,200)(1,0)[\phantom{\big(M(A,B)\big)_{\substack{A\in\u\C \\ B\in\u\D}}} `\phantom{\crta\C\times\crta\D}`
   {{[A]\times [C]}}]{1200}1a
\putmorphism(130,200)(0,-1)[M(A,B)\times N(B,C)`M'(A',B')\times N'(B',C')`   f_{A,B}\times g_{B,C}]{400}1l
\putmorphism(1420,200)(0,-1)[\crta\C\times\crta\J`\crta\E\times\crta\K`]{400}0r
\putmorphism(1480,200)(0,-1)[``\crta F\times\crta H]{400}1r
\putmorphism(220,-200)(1,0)[\phantom{M(A,B)\times N(B,C)}`\phantom{\crta\C\times\crta\D}`
    {{[A']\times [C']}}]{1200}1b
\put(550,0){\fbox{$\tilde\nu_{A,B}\times\tilde\nu'_{B,C}$}}
\putmorphism(860,-570)(0,-1)[1` `]{300}0a
\putmorphism(1050,-490)(2,1)[``]{280}1r
\putmorphism(170,-950)(2,1)[``\langle A', C'\rangle]{1760}0r
%
\putmorphism(10,-190)(2,-1)[``]{890}1l
\putmorphism(-530,30)(2,-1)[``!]{2050}0l
\put(760,-460){\fbox{$\gamma'$}}
\efig }
\end{equation}
It is easily seen how the 2-cells $\gamma$ and $\gamma'$ are induced. 

\medskip

We next observe the following diagram: 
\begin{equation} \eqlabel{total 2-cell}
\scalebox{0.86}{\bfig
\putmorphism(-830,930)(4,1)[``]{500}1l
\putmorphism(-830,960)(4,1)[``!]{1100}0l
\put(-50,940){\fbox{$\theta^B$}}
\putmorphism(-1500,800)(1,0)[M(A,B)\times N(B,C)`\phantom{\amalg_{B\in\u\D}M(A,B)\times N(B,C)}`\iota^{B}]{1400}1a
\putmorphism(-100,800)(1,0)[\phantom{\amalg_{B\in\u\D}M(A,B)\times N(B,C)}``{{[A,C]}}]{1100}1a
\putmorphism(-1500,800)(0,-1)[`` f_{A,B}\times g_{B,C}]{400}1r
\putmorphism(20,1150)(0,-1)[1` `]{300}0a
\putmorphism(-170,1200)(4,-1)[``]{1580}1r
\putmorphism(-60,1230)(4,-1)[``\langle A, C\rangle]{1160}0r
\putmorphism(-100,800)(0,-1)[\amalg_{B\in\u\D}M(A,B)\times N(B,C)`\amalg_{B'\in\u\HH}M'(A',B')\times N'(B',C')`  ]{400}0l
\putmorphism(30,800)(0,-1)[``  h_{A,C}]{400}1l
\putmorphism(1240,800)(0,-1)[\crta\C\times\crta\J` \crta\E\times\crta\K`\crta F\times\crta H]{400}1l
\putmorphism(-1500,400)(1,0)[M'(A',B')\times N'(B',C')`\phantom{\amalg_{B'\in\u\HH}M'(A',B')\times N'(B',C')}`\iota^{B'}]{1400}1a
\putmorphism(20,380)(0,-1)[`1 `]{300}0a
\putmorphism(410,130)(4,1)[``]{430}1r
\putmorphism(-370,-120)(4,1)[``\langle A', C'\rangle]{1760}0r
\put(-800,580){\fbox{$\zeta^B$}}
\putmorphism(-40,400)(1,0)[\phantom{\amalg_{B\in\u\D}M(A,B)\times N(B,C)}``{{[A',C']}}]{1040}1a
\putmorphism(-1410,390)(4,-1)[``]{1620}1l
\putmorphism(-1280,320)(4,-1)[``!]{2050}0l
\put(-80,200){\fbox{$\theta^{B'}$}}
\efig}
\end{equation}
Taking into account every $B$ for which there exist the 1-cells $f_{A,B}$ and $g_{B,C}$, the composite 1-cells 
$\iota^{B'}(f_{A,B}\times g_{B,C})$ on the left induce a 1-cell $h_{A,C}$ and an equivalence 2-cell $\zeta^B$. 
We set for the total 2-cell $(\crta F\times\crta H)\langle A, C\rangle!\Rightarrow \langle A', C'\rangle!(f_{A,B}\times g_{B,C})$ 
to be the above composite equivalence 2-cell \equref{nu-nu}. 
We now apply the dual of \leref{lema} d) by making the following dual correspondence:
$$p_2\mapsto \iota^{B'}, p'_2\mapsto \iota^B,  f\mapsto \crta F\times\crta H, g\times h\mapsto h_{A,C}, $$
$$s_2\mapsto \langle A, C\rangle!, q_2\mapsto \langle A', C'\rangle!, $$
$$s\mapsto [A,C], t\mapsto [A',C'], $$
$$\omega_2\mapsto\zeta^B, \zeta_2\mapsto\theta^{B'}, \theta_2\mapsto\theta^B$$
$$\beta\mapsto \text{2-cell} \hspace{0,3cm} \equref{nu-nu}.$$
Then there exists an equivalence 2-cell $(\crta F\times\crta H)[A, C]\Rightarrow [A', C']h_{A,C}$ 
corresponding to the right rectangular in the diagram above (the dual of $\gamma$ from \leref{lema}). 
We take this 2-cell for the desired equivalence 2-cell $\nu_{A,C}$ as in \equref{nu} for the composition of matrices 
$\big(M(A,B)\big)_{\substack{A\in\u\C \\ B\in\u\D}}$ and $\big(N(B,C)\big)_{\substack{B\in\u\D \\ C\in\u\J}}$.

\medskip

To get an equivalence 2-cell $\chi^{F,H}_{A,C}$ from \equref{nu} is easy: it is induced by the given equivalence 2-cells $\chi^F_A$ and $\chi^H_C$.

\subsection{Composition of 2-cells in $D_1$}

Let two 2-cells in $D_1$ 
$$\Sigma_{A,B}^1: \frac{\nu_{A,B}^1}{[\xi\vert\Id_{[A',B']}]} \Rrightarrow \frac{[\Id_{[A,B]}\vert\crta\lambda\times\crta\rho]}{\nu^{'1}_{A,B}} 
\quad\text{and}\quad
\Sigma_{B,C}^2: \frac{\nu_{B,C}^2}{[\xi\vert\Id_{[A',B']}]} \Rrightarrow \frac{[\Id_{[A,B]}\vert\crta\lambda\times\crta\rho]}{\nu^{'2}_{B,C}} $$
be given. Recall how the 2-cell $\tilde\nu_{A,B}$ in $V$ above \equref{nu-nu} is induced by $\nu_{A,B}$. Let us denote the 2-cells induced analogously by 
$\nu_{A,B}^1, \nu^{'1}_{A,B}, \nu_{B,C}^2, \nu^{'2}_{B,C}$ as follows: $\tilde\nu_{A,B}^1, \tilde\nu^{'1}_{A,B}, \tilde\nu_{B,C}^2, \tilde\nu^{'2}_{B,C}$, 
respectively. Consider the prism $P_{\tilde\nu_{A,B}^1}$ with basis $\tilde\nu_{A,B}^1$ (and analogously the prism $P_{\tilde\nu_{B,C}^2}$ with basis 
$\tilde\nu_{B,C}^2$) obtained by concatenation of the following prisms: 
$\Sigma_{A,B}^1, C^1_{F\times G}$ from \coref{alfa-x-beta}, identity 3-cells on $\tilde\theta_1$ and $\tilde\theta_1^i$. Analogously, 
the prism $P_{\tilde\nu_{B,C}^2}$ is obtained by concatenation of the prisms 
$\Sigma_{B,C}^2, C^2_{F\times G}$, identity 3-cells on $\tilde\theta_2$ and $\tilde\theta_2^i$. Seeing $P_{\tilde\nu_{A,B}^1}$ and $P_{\tilde\nu_{B,C}^2}$ 
as $\tilde P^1$ and $\tilde P^2$ in \coref{middle prism}, we obtain a unique 3-cell 
$$\Gamma': \displaystyle{\frac{[\Id\vert\tilde\nu_{A,B}^1\times\tilde\nu_{B,C}^2]}{[\xi\times\xi'\vert\Id_{[A', B']\times [B',C']}]}} \Rrightarrow 
\displaystyle{\frac{[\Id_{[A, B]\times [B,C]}\vert\crta\lambda\times\crta\sigma]}{[\tilde\nu^{'1}_{A,B}\times\tilde\nu^{'2}_{B,C}\vert\Id]}}.$$
This is a prism with basis $\tilde\nu_{A,B}^1\times\tilde\nu_{B,C}^2$, as in the middle of \equref{nu-nu}. Concatenate to it the identity 3-cells on 
$\gamma$ and $\gamma'$ from \equref{nu-nu} and consider the obtained prism $P_{total}$ as the prism whose basis is the total 2-cell in \equref{total 2-cell}. 
This total 2-cell we treat as the 2-cell $\alpha_i$ in \prref{dual cutting} (in the reversed direction) and it induces the equivalence 2-cell $\gamma: 
(\crta F\times\crta H)[A, C]\Rightarrow [A', C']h_{A,C}$, which is taken for $\nu_{A,C}$. 
Observe that the 1-cell $h_{A,C}$ can be written as $\amalg_{B\in\u\D}f_{A,B}\times g_{B,C}$. 
This 2-cell $\gamma=\nu_{A,C}$ corresponds to ``$\gamma$ with reversed order'' in \prref{dual cutting}. 
The prism $P_{total}$ corresponds to the 3-cell $P_{\alpha_i}$ 
in there (in the corresponding mapping direction), and we finally obtain a unique 3-cell $\Gamma$ with basis $\gamma$, that is, 
a 3-cell $\Gamma: \nu_{A,C}\Rrightarrow \nu_{A,C}'$. This $\Gamma$ is the horizontal composition 2-cell of the 2-cells $\Sigma_{A,B}^1$ and 
$\Sigma_{B,C}^2$ in $D_1$.

\section{Relating matrices and spans in a tricategory} 

In this Section we are going to construct functors between the $(1\times 2)$-categories of matrices $\M$ and spans $\S$ in a 1-strict tricategory $V$ with a terminal object, small 3-coproducts, 3-products and 3-pullbacks. Such a functor is internal in $\Bicat_3$, so for that purpose we will define pseudofunctors between the bicategories $C_1$ of spans from \seref{3-spans} and $D_1$ of matrices from \seref{3-matrices}, and additionally check their compatibility with the (horizontal)
composition on the 3-pullback. Recall that $C_0=D_0$ is the 2-category of small categories. 
We will obtain a lax internal functor $\S\to\M$ and a colax internal functor $\M\to\S$. 
For completeness and the sake of the next Section we introduce two formal definitions.

\subsection{Internal and enriched functors in 1-strict tricategories} \sslabel{int-enrich-fun}

In this Subsection we only give the two definitions. Referring to the notion and notation from \cite[Definition 6.2]{Fem} we introduce:

\begin{defn}
Let $\C,\D$ be categories internal in a 1-strict tricategory $V$. We say that $F:\C\to\D$ is a (pseudo-/lax/colax) functor internal in $V$ 
if it consists of: 
\begin{enumerate}
\item pseudofunctors $F_0: C_0\to D_0$ and $F_1: C_1\to D_1$ such that $s\comp F_1=F_0\comp s, t\comp F_1=F_0\comp t$;
\item pseudonatural transformations \vspace{-0,3cm}
$$F_\times(f,g): F_1(g)\times_{D_0}F_1(f)\Rightarrow F_1(g\times_{C_0}f)\quad\text{ and}\quad F_u(A): u_{F_0(A)}\Rightarrow F_1(u_A) 
\vspace{-0,3cm}$$ 
in the lax functor case ($F_\times(f,g): F_1(g\times_{C_0}f) \Rightarrow F_1(g)\times_{D_0}F_1(f)$ and  
$F_u(A): F_1(u_A)\Rightarrow u_{F_0(A)}$ in the colax functor case, and for a pseudofunctor require $F_\times(f,g)$ and $F_u(A)$ to be 
equivalence 2-cells) 
for objects $A\in C_0$ and 1h-cells $f,g\in C_1$, whose components are globular equivalences, and 
\item modifications $\Omega_{a^*}, \Omega_{l^*}, \Omega_{r^*}$:
\[\adjustbox{scale=0.88,center}{
\begin{tikzcd}[row sep = huge, column sep = huge]
(F_1(R) \hspace{-0,06cm}\times\hspace{-0,03cm} F_1(S)) \hspace{-0,06cm}\times\hspace{-0,03cm} F_1(T)
	\ar[r, Rightarrow, "{a^*_{F_1(R), F_1(S), F_1(T)}}"]
	\ar[d, Rightarrow, "{F_{\times} \times \Id_{F_1(T)}}"']
&F_1(R) \hspace{-0,06cm}\times\hspace{-0,03cm} (F_1(S) \hspace{-0,06cm}\times\hspace{-0,03cm} F_1(T))
	\ar[r, Rightarrow, "{\Id_{F_1(R)} \times F_{\times}}", "{}"{name=S, below}]
&F_1(R) \hspace{-0,06cm}\times\hspace{-0,03cm} F_1(S \hspace{-0,06cm}\times\hspace{-0,03cm} T)
	\ar[d, Rightarrow, "{F_{\times}}"] \\
F_1(R \hspace{-0,06cm}\times\hspace{-0,03cm} S) \hspace{-0,06cm}\times\hspace{-0,03cm} F_1(T)
	\ar[r, Rightarrow, "{F_{\times}}"{name=T, above}]
&F_1((R \hspace{-0,06cm}\times\hspace{-0,03cm} S) \hspace{-0,06cm}\times\hspace{-0,03cm} T)
	\ar[r, Rightarrow, "{F_1(\alpha_{R, S, T})}"]
&F_1(R \hspace{-0,06cm}\times\hspace{-0,03cm} (S \hspace{-0,06cm}\times\hspace{-0,03cm} T))
	\ar[from=S, to=T, triple, "{\Omega_{a^*}}"description]
\end{tikzcd}
}\]

\[
\begin{tikzcd}[row sep = huge, column sep = huge]
F(R) \times F(u)
	\ar[d, Rightarrow, "{F_{\times}}"', "{}"{name=SR}]
&F(R)
	\ar[d, equal, "{}"'{name=TR}, "{}"{name=TL}]
	\ar[l, Rightarrow, "{\Id_{F_1(R)} \times F_{u}}"']
	\ar[r, Rightarrow, "{F_{u} \times \Id_{F_1(R)}}"]
&F(u) \times F(R)
	\ar[d, Rightarrow, "{F_{\times}}", "{}"'{name=SL}] \\
F(R \times u)
	\ar[r, Rightarrow, "{F(r^*)}"]
&F(R)
&F(u \times R)
	\ar[l, Rightarrow, "{F(l^*)}"']
	\ar[from=SR, to=TR, triple, "{\Omega_{r^*}}"]
	\ar[from=SL, to=TL, triple, "{\Omega_{l^*}}"']
\end{tikzcd}
\]
which satisfy the two diagrammatic equations (A1) and (A2) in the Appendix. 
\end{enumerate}
\end{defn}

\bigskip

Referring to the notion and notation from \cite[Definition 8.1]{Fem} we introduce:

\begin{defn}
Let $\Tau,\Tau'$ be categories enriched in a 1-strict tricategory $V$. We say that $F:\Tau\to\Tau'$ is a functor enriched in $V$ 
if it consists of: 
\begin{enumerate}
\item an assignment $F_0:\Ob\Tau\to\Ob\Tau'$;
\item a 1-cell $F_1:\Tau(A,B)\to\Tau'(F(A),F(B))$ in $V$ for all $A,B\in\Ob\Tau$; 
\item equivalence 2-cells $F_c: F_1\cdot\circ\Rightarrow\circ'\cdot(F_1\times F_1)$ and $F_{I_A}: F_1\cdot I_A\Rightarrow I'_{F(A)}$, and
\item bijective 3-cells $\Omega_{a^\dagger}, \Omega_{l^\dagger}, \Omega_{r^\dagger}$:
\[
\begin{tikzcd}
{F_1(-\circ(-\circ-))}
	\ar[r, Rightarrow, "{F_1(a^\dagger)}"]
	\ar[d, Rightarrow, "{F_c^L}"']
&{F_1((-\circ-)\circ-)}
	\ar[d, Rightarrow, "{F_c^R}"] \\
{F_1(-)\circ(F_1(-)\circ F_1(-))}
	\ar[r, Rightarrow, "{a^{'\dagger}}"]
	\ar[ru, triple, "{\Omega_{a^\dagger}}" description]
&{(F(-)\circ F(-))\circ F(-)}
\end{tikzcd}
\]
where $F_c^L$ and $F_c^R$ are the obvious 2-cells induced by $F_c$ and $\Id\times F_C$, and by $F_c$ and $F_c\times\Id$, respectively, 
\[
\begin{tikzcd}
{F_1(I_B\circ -)}
	\ar[r, Rightarrow, "{F_1(l^\dagger)}"]
	\ar[d, Rightarrow, "{F_c(\Id)}"']
&{F_1(-)}
	\ar[d, Rightarrow, "{l^{'\dagger}}"] \\
{(F_1\cdot I_B)\circ' F_1(-)}
	\ar[r, Rightarrow, "{F_{I_B}\circ'\Id}"]
	\ar[ru, triple, "{\Omega_{l^\dagger}}" description]
&{I'_{F(B)}\circ' F_1(-)}
\end{tikzcd}
\]

\[
\begin{tikzcd}
{F_1(-\circ I_A)}
	\ar[r, Rightarrow, "{F_1(r^\dagger)}"]
	\ar[d, Rightarrow, "{F_c(\Id)}"']
&{F_1(-)}
	\ar[d, Rightarrow, "{r^{'\dagger}}"] \\
{F_1(-)\circ' (F_1\cdot I_A)}
	\ar[r, Rightarrow, "{\Id\circ' F_{I_A}}"]
	\ar[ru, triple, "{\Omega_{r^\dagger}}" description]
&{F_1(-)\circ' I'_{F(A)} }
\end{tikzcd}
\]
which satisfy axioms analogous to those from the Appendix (substitute the horizontal composition on pullbacks, there denoted by $\times$ for short, with 
$\circ$, the pseudonatural transformations $F_\times$ and $F_u$ by the 2-cells $F_c^{-1}$ and $F_{I}^{-1}$, respectively, and the 
modifications $\Omega_{a^*}, \Omega_{l^*}, \Omega_{r^*}$ by 3-cells $\Omega_{a^\dagger}, \Omega_{l^\dagger}, \Omega_{r^\dagger}$, respectively. 
\end{enumerate}
\end{defn}

\subsection{From spans to matrices} \sslabel{span->matrix}

We start by defining a pseudo functor $En: C_1\to D_1$. It maps a span $\u\C\bullet 1 \stackrel{r_1}{\longleftarrow} R \stackrel{r_2}{\longrightarrow} \u\D\bullet 1$ 
to the matrix $En(R)$ whose $(A,B)$-component is given by the 3-pullback 
\begin{equation} \eqlabel{matrix out of span}
\end{equation}
$$\scalebox{0.86}{\bfig
\putmorphism(30,800)(1,0)[\phantom{En(V)(i,j)}`\phantom{B}`!]{830}1a
\putmorphism(840,800)(0,-1)[1`(\u\C \bullet 1) \times (\u\D \bullet 1)`\langle A,B\rangle]{500}1r
\putmorphism(0,800)(0,-1)[En(R)(A,B)`R`\iota_{A,B}^R]{500}1l
\putmorphism(0,300)(1,0)[\phantom{V}`\phantom{(I \bullet 1) \times (J \bullet 1)}`\langle r_1, r_2\rangle]{800}1a
\put(300,550){\fbox{$\omega^R_{A,B}$}}
\efig} \vspace{-1,4cm}
$$
for each $A\in\C, B\in\D$. 

We next compose the above equivalence 2-cell with the equivalence 2-cells inducing morphisms into 3-products (so far they are well-known 
and we do not label them) and set for the composite equivalence 2-cells $\omega^R_A$ and $\omega^R_B$:
$$\omega^R_A:=
\scalebox{0.8}{
\bfig
\putmorphism(30,550)(1,0)[\phantom{En(V)(i,j)}`\phantom{B}`!]{1040}1a
\putmorphism(1040,550)(0,-1)[1``\langle A,B\rangle]{500}1l
\putmorphism(0,550)(0,-1)[En(R)(A,B)`R`\iota_{A,B}^R]{500}1l
\putmorphism(0,50)(1,0)[\phantom{V}`(\u\C \bullet 1) \times (\u\D \bullet 1)`\langle r_1, r_2\rangle]{800}1a
\put(300,300){\fbox{$\omega^R_{A,B}$}}
\putmorphism(1050,580)(1,-2)[``\sigma_A]{430}1r
\putmorphism(750,50)(2,-1)[`\u\C\bullet 1`p_1]{700}1r
\putmorphism(-210,60)(4,-1)[``]{1780}1l
\putmorphism(-210,20)(4,-1)[``r_1]{1780}0l
\efig}
\quad
\omega^R_B:=
\scalebox{0.8}{
\bfig
\putmorphism(30,550)(1,0)[\phantom{En(V)(i,j)}`\phantom{B}`!]{1040}1a
\putmorphism(1040,550)(0,-1)[1``\langle A,B\rangle]{500}1l
\putmorphism(0,550)(0,-1)[En(R)(A,B)`R`\iota_{A,B}^R]{500}1l
\putmorphism(0,50)(1,0)[\phantom{V}`(\u\C \bullet 1) \times (\u\D \bullet 1)`\langle r_1, r_2\rangle]{800}1a
\put(300,300){\fbox{$\omega^R_{A,B}$}}
\putmorphism(1050,580)(1,-2)[``\sigma_B]{430}1r
\putmorphism(750,50)(2,-1)[`\u\D\bullet 1`p_2]{700}1r
\putmorphism(-210,60)(4,-1)[``]{1780}1l
\putmorphism(-210,20)(4,-1)[``r_2]{1780}0l
\efig}
$$

\begin{lma} \lelabel{span->chi}
Given a morphism of spans \equref{1-cells in C_1} (we write $R$ and $S$ instead of $A$ and $B$ there), together with the above equivalence 
2-cells $\omega^R_{A,B}$ and $\omega^S_{A,B}$ it induces equivalence 2-cells $\chi^F_A, \chi^G_B$ as in \equref{one-sided chi's}, which in turn induce 
an equivalence 2-cell $\chi^{F,G}_{A,B}$ as in \equref{nu}. 
\end{lma}

\begin{proof}
Consider the composite equivalence 2-cell
$$\alpha_A:=
\scalebox{0.86}{\bfig
 \putmorphism(-200,500)(1,0)[En(R)(A,B)``\iota^R_{A,B}]{630}1a
 \putmorphism(450,500)(1,0)[R`1 `!]{400}1a
 \putmorphism(780,500)(1,0)[\phantom{F(B)}` \u\C \bullet 1`\sigma_A]{550}1a

\putmorphism(0,520)(0,-1)[\phantom{Y_2}``=]{430}1l
\putmorphism(850,520)(0,-1)[\phantom{Y_2}``=]{430}1l
\putmorphism(1340,520)(0,-1)[\phantom{Y_2}``=]{430}1r
\put(350,300){\fbox{$\kappa$}}

 \putmorphism(0,120)(1,0)[``!]{850}1a
 \putmorphism(770,120)(1,0)[\phantom{F(B)}` `\sigma_A]{570}1a

\putmorphism(0,140)(0,-1)[\phantom{Y_2}``=]{430}1l
\putmorphism(1340,140)(0,-1)[\phantom{Y_2}``=]{430}1r
\put(600,-100){\fbox{$\omega^R_A$}}
 \putmorphism(0,-260)(1,0)[``\iota^R_{A,B}]{470}1a
 \putmorphism(340,-260)(1,0)[\phantom{F(B)}` `r_1]{1000}1a
\efig}
$$
and similarly $\alpha_B$ (changing to $\sigma_B, \omega^R_B$ and $r_2$). By the 3-coproduct property there are equivalence 2-ells 
$\gamma_A:\sigma_A!\Rightarrow r_1$ and $\gamma_B:\sigma_B!\Rightarrow r_2$. Then set 
$$\chi^F_A:=
\scalebox{0.86}{\bfig
\putmorphism(850,670)(-3,-1)[``]{1180}1l
\putmorphism(850,690)(-3,-1)[``\sigma_A]{1180}0l

\putmorphism(-150,250)(1,0)[\u\C \bullet 1`R`r_1]{600}{-1}a
 \putmorphism(540,380)(1,1)[` `!]{80}1r

 \putmorphism(-150,-200)(1,0)[\u\E\bullet 1` S` s_1]{600}{-1}a
 \putmorphism(370,-200)(1,-1)[\phantom{\u\C \bullet 1}`1 ` !]{350}1r

\putmorphism(-180,250)(0,-1)[\phantom{Y_2}``F\bullet 1]{450}1l
\putmorphism(450,250)(0,-1)[\phantom{Y_2}``]{450}1r
\putmorphism(300,250)(0,-1)[\phantom{Y_2}``f]{450}0r
\putmorphism(730,600)(0,-1)[1``=]{1130}1r
 \putmorphism(430,-470)(-3,1)[``\sigma_{A'}]{400}1l
\put(250,370){\fbox{$\gamma_A$}}
\put(100,20){\fbox{$\alpha$}}
\put(550,0){\fbox{$\kappa$}}
\put(230,-330){\fbox{$\gamma_{A'}^{-1}$}}
\efig}
$$ 
and similarly for $\chi^G_B$. Finally apply \leref{lema} d) to get an equivalence 2-cell $\chi^{F,G}_{A,B}$. 
\end{proof}

Given a 1-cell \equref{1-cells in C_1} in $C_1$. 
To define its image by $En$ we consider the following cube: 
\begin{equation} \eqlabel{cube}
\scalebox{0.86}{\bfig
\putmorphism(-650,400)(1,0)[En(R)(A,B)`1`!]{1100}1a
\putmorphism(-650,400)(0,-1)[` `En(f)(A,B)]{800}1l
\putmorphism(-800,400)(3,-1)[`` ]{1210}1l
\putmorphism(-750,350)(3,-1)[`` \iota_{A,B}^R]{1270}0l
\putmorphism(450,400)(0,-1)[` `]{800}1r
\putmorphism(460,330)(0,-1)[` `=]{800}0l
\putmorphism(260,40)(1,0)[R`\crta\C\times\crta\D`\langle r_1,r_2\rangle]{1100}1a
\putmorphism(250,40)(0,-1)[` S`]{800}1l
\putmorphism(270,120)(0,-1)[` `f]{800}0l
\putmorphism(1440,40)(0,-1)[` `\crta F\times\crta G]{800}1r
\putmorphism(260,-750)(1,0)[`\crta\E\times\crta\HH`\langle s_1,s_2\rangle]{1100}1a
\put(150,220){\fbox{$\omega^R_{A,B}$}}
\putmorphism(280,450)(3,-1)[`` ]{1190}1r
\putmorphism(360,460)(3,-1)[`` \langle A,B\rangle]{1330}0r
\putmorphism(-650,-400)(1,0)[En(S)(A',B')`1`!]{1100}1a
\putmorphism(-830,-410)(3,-1)[`` \iota_{A',B'}^S]{1210}1l
\putmorphism(280,-340)(3,-1)[`` ]{1190}1r
\putmorphism(360,-320)(3,-1)[`` \langle A',B'\rangle]{1330}0r
\put(-50,-50){\fbox{$\kappa$}}
\put(280,-560){\fbox{$\omega^S_{A',B'}$}}
\put(800,-260){\fbox{$\chi^{F,G}_{A,B}$}}
\efig}
\end{equation}
The 1-cell $En(f)(A,B)$ and the left 2-cell are to be defined. 
The top and bottom are the equivalence 2-cells \equref{matrix out of span}, the right 2-cell is the equivalence 2-cell from \leref{span->chi}, 
the back is a terminal (equivalence) 2-cell, and the front one is the equivalence 2-cell from \prref{2-cell-prop}, 
let us denote it by $\gamma^f$. By the 3-pullback property of $En(S)(A',B')$ the composite equivalence 2-cell
$$
\scalebox{0.86}{\bfig
 \putmorphism(0,900)(1,0)[``!]{450}1a
 \putmorphism(310,900)(1,0)[\phantom{F(B)}` `\langle A',B'\rangle]{1200}1a

\putmorphism(0,900)(0,-1)[\phantom{Y_2}``=]{400}1r
\putmorphism(440,900)(0,-1)[\phantom{Y_2}``=]{400}1r
\putmorphism(1510,900)(0,-1)[\phantom{Y_2}``=]{400}1r
\put(900,750){\fbox{$\chi_{A,B}^{-1}$}}
\put(350,280){\fbox{$\omega^R_{A,B}$}}

 \putmorphism(0,500)(1,0)[``!]{450}1a
 \putmorphism(320,500)(1,0)[\phantom{F(B)}` `\langle A,B\rangle]{640}1a
 \putmorphism(830,500)(1,0)[\phantom{F(B)}` `\crta F\times\crta G]{700}1a

\putmorphism(0,520)(0,-1)[\phantom{Y_2}``=]{450}1l
\putmorphism(940,520)(0,-1)[\phantom{Y_2}``=]{450}1l
\putmorphism(1510,520)(0,-1)[\phantom{Y_2}``=]{450}1r

 \putmorphism(0,100)(1,0)[``\iota^R]{470}1a 
 \putmorphism(320,100)(1,0)[\phantom{F(B)}` `\langle r_1,r_2\rangle]{640}1a
 \putmorphism(830,100)(1,0)[\phantom{F(B)}` `\crta F\times\crta G]{700}1a

\putmorphism(0,120)(0,-1)[\phantom{Y_2}``=]{430}1l
\putmorphism(440,120)(0,-1)[\phantom{Y_2}``=]{400}1r
\putmorphism(1510,120)(0,-1)[\phantom{Y_2}``=]{430}1r
\put(830,-120){\fbox{$\gamma^f$}}

 \putmorphism(0,-280)(1,0)[``\iota^R]{450}1a
 \putmorphism(320,-280)(1,0)[\phantom{F(B)}` `f]{640}1a
 \putmorphism(830,-280)(1,0)[\phantom{F(B)}` `\langle s_1,s_2\rangle]{700}1a
\efig}
$$
induces a 1-cell $En(f)(A,B):En(R)(A,B)\to En(S)(A,B)$, an equivalence 2-cell $\zeta^f:f\iota^R_{A,B}\Rightarrow \iota_{A',B'}^S En(f)(A,B)$, 
the left face of the cube, and a bijective 3-cell 
$\Sigma: \threefrac{[\kappa\vert\Id_{\langle A',B'\rangle}]}{[\Id_{En(f)}\vert\omega^S_{A',B'}]}{[\zeta^f\vert\Id_{\langle s_1,s_2\rangle}}
\Rrightarrow
\threefrac{[\Id_!\vert\chi^{-1}_{A,B}]}{[\omega^R_{A,B}\vert\Id_{\crta F\times\crta G}]}{[\Id_{\iota^R}\vert\gamma^f]}$.  
Then we set for the equivalence 2-cell \equref{nu} determining a 1-cell in $D_1$ 
the concatenation of the four faces of the cube: 
\begin{equation} \eqlabel{En(f)}
\nu_{A,B}:=\threefrac{[\omega^R_{A,B}\vert\Id_{\crta F\times\crta G}]}{(\zeta^f\vert\gamma^f)}{[\Id_{En(f)(A,B)}\vert(\omega^S_{A',B'})^{-1}]}.
\end{equation} 
(We recall that 
by $(\zeta^f\vert\gamma^f)$ we mean the horizontal concatenation of 2-cells, simplifying the writing of the horizontal composition of 2-cells.)
Thus $En(\alpha\vert f\vert\beta)$ is the family of these $\nu_{A,B}$ and  $\chi^{F,G}_{A,B}$ from \leref{span->chi}, for $(A,B)\in\u\C\times\u\D$.

\bigskip

For the image of 2-cells in $D_1$ we procced as follows. 
Given a 2-cell \equref{3-cells} in $C_1$. 
Observe that for the domain and codomain 1-cells in the corresponding cubes 
\equref{cube} three faces remain the same ($\omega$'s and $\kappa$) and the resting three faces differ. Let 
$\nu_{A,B}=\threefrac{[\omega^R_{A,B}\vert\Id_{\crta F\times\crta G}]}{(\zeta^f\vert\gamma^f)}{[\Id_{En(f)(A,B)}\vert(\omega^S_{A',B'})^{-1}]}$ 
and 
$\nu_{A,B}':=\threefrac{[\omega^R_{A,B}\vert\Id_{\crta F'\times\crta G'}]}{(\zeta^g\vert\gamma^g)}{[\Id_{En(g)(A,B)}\vert(\omega^S_{A',B'})^{-1}]}$ 
denote the images of the 1-cells $(\alpha\vert f\vert\beta)$ and $(\gamma\vert g\vert\delta)$ by $En$. 

The 3-cells $\Sigma$ and $\Omega$ from \equref{3-cells} induce, on one hand, 3-cells (prisms) over $\chi^F_A$ and $\chi^G_B$ in \leref{span->chi} 
and then by \prref{2-cell-prop} a prism between $\chi^{F,G}_{\alpha,\beta}$ and $\chi^{F',G'}_{\gamma,\delta}$, and on the other hand, 
a 3-cell (prism) between $\gamma^f$ and $\gamma^g$ also by \prref{2-cell-prop}. Taking identity 3-cell over $\omega^R$ and concatenating the obtained prisms, 
we obtain a 3-cell 
$\upkappa_0$ as in \leref{pullb-lema}, and thus also a 2-cell $\gamma: En(f)(A,B)\Rightarrow En(g)(A,B)$ and a prism with basis $\zeta^f$ and opposite 
face $\zeta^g$, so that all the 1-cells joining their corresponding vertexes are identities, with $\alpha=\Id_!$ and $\beta=[\Id_{\iota^R}\vert\xi]$. 
Next, apart from the prism over $\gamma^f$, between the resting 2-cells comprising $\nu_{A,B}$ and $\nu_{A,B}'$ take identity 3-cells and their corresponding 
prisms. Then for the obvious concatenation of these prisms we set to be the image by $En$ of the 2-cell \equref{3-cells}.

For the sake of saving space we skip the proof of compatibility of $En$ with the composition of 1- and 2-cells. We only record that 
we constructed a 1-cell $w:\amalg_{B\in\u\D}En(R)(A,B)\times En(S)(B,C)\to En(R\times_{\u\D\bullet 1} S)(A,C)$ for spans 
$\u\C\bullet 1 \stackrel{r_1}{\longleftarrow} R \stackrel{r_2}{\longrightarrow} \u\D\bullet 1$ and 
$\u\D\bullet 1 \stackrel{s_1}{\longleftarrow} S \stackrel{s_2}{\longrightarrow} \u\E\bullet 1$, proving the laxity of an internal functor 
$\S\to\M$.


\subsection{From matrices to spans}

We define here the pseudofunctor $Int: D_1\to C_1$. 
Given a matrix $\big(M(A,B)\big)_{\substack{A\in\u\C \\ B\in\u\D}}$ in $V$ (with the corresponding 1-cell to $\crta\C\times\crta\D$), the pseudofunctor $Int$ maps it into the span 
$\u\C \bullet 1 \stackrel{m_1}{\longleftarrow} \amalg_{\substack{A\in\u\C \\ B\in\u\D}}M(A,B) \stackrel{m_2}{\longrightarrow} \u\D \bullet 1$. 
The 1-cells $m_1,m_2$ are induced by the following 1-cells on $M(A,B)$, for fixed $A\in\u\C, B\in\u\D$: 
the unique morphism to 1 followed by the coprojections to $\u\C \bullet 1$ and $\u\D \bullet 1$, respectively.

Given a 1-cell in $D_1$, with a family of 2-cells $\nu_{A,B}$ as in the most left rectangular diagram below, we are going to define a 1-cell 
$\amalg f_{A,B}$ and 2-cells $\alpha$ and $\beta$ as on the right in: 
\begin{equation} \eqlabel{Int-image of nu}
\scalebox{0.88}{
\bfig
\putmorphism(220,530)(3,1)[``]{300}1l
\putmorphism(-80,460)(3,1)[`\amalg_{\substack{A\in\u\C \\ B\in\u\D}}M(A,B)`\iota^M_{A,B}]{1000}0l
\put(750,540){\fbox{$\zeta^M_1$}}
\putmorphism(730,770)(3,-1)[``]{1170}1r
\putmorphism(640,830)(3,-1)[``m_1]{1160}0r
\putmorphism(0,400)(1,0)[\phantom{M(A,B)} `\phantom{(\u{\C}\bullet 1}`{{[A,B]}}]{900}1a
\putmorphism(0,400)(0,-1)[M(A,B)`N(A',B')`f_{A,B}]{400}1l
\putmorphism(900,400)(0,-1)[\crta\C\times\crta\D`\crta\E\times\crta\HH`]{400}0r
\putmorphism(870,400)(0,-1)[``\crta F\times\crta G]{400}1r 
\putmorphism(870,400)(1,0)[\phantom{(\u{\C}\times \bullet 1}``p_1]{780}1a
\putmorphism(1800,400)(0,-1)[\u{\C}\bullet 1`\u{\E}\bullet 1`]{400}0r
\putmorphism(1820,400)(0,-1)[``F\bullet 1]{400}1r
\putmorphism(910,0)(1,0)[\phantom{(\u{\C}\bullet 1}``p_1']{760}1a
\putmorphism(0,0)(1,0)[\phantom{N(A',B')}`\phantom{V}`{{[A',B']}}]{800}1a
\put(260,240){\fbox{$\nu_{A,B}$}}
\put(1460,200){\fbox{$\omega_1$}}
\putmorphism(1190,-220)(3,1)[``]{260}1r
\putmorphism(330,-520)(3,1)[``n_1]{1760}0r
\putmorphism(-210,-10)(3,-1)[``]{1000}1l
\putmorphism(-150,-70)(3,-1)[`\amalg_{\substack{A'\in\u\E \\ B'\in\u\HH}} N(A',B')`\iota^N_{A',B'}]{1050}0l
\put(720,-230){\fbox{$\zeta^N_1$}}
\efig }
\quad
\scalebox{0.84}{
\bfig
\putmorphism(-250,400)(1,0)[\u\C \bullet 1`\amalg_{\substack{A\in\u\C \\ B\in\u\D}}M(A,B)`m_1]{700}{-1}a
 \putmorphism(700,400)(1,0)[\phantom{\u\C \bullet 1}`\u\D \bullet 1 `m_2]{500}1a

 \putmorphism(-250,-50)(1,0)[\u\E\bullet 1` \amalg_{\substack{A'\in\u\E \\ B'\in\u\HH}} N(A',B')` n_1]{700}{-1}a
 \putmorphism(700,-50)(1,0)[\phantom{\u\C \bullet 1}`\u\HH\bullet 1. ` n_2]{500}1a

\putmorphism(-240,400)(0,-1)[\phantom{Y_2}``F\bullet 1]{400}1l
\putmorphism(420,400)(0,-1)[\phantom{Y_2}``]{400}1r
\putmorphism(410,400)(0,-1)[\phantom{Y_2}``\amalg f_{A,B}]{450}0r
\putmorphism(1200,400)(0,-1)[\phantom{Y_2}``G\bullet 1]{400}1r
\put(0,170){\fbox{$\alpha$}}
\put(750,170){\fbox{$\beta$}}
\efig}
\end{equation}
The 2-cells $\omega_1, \zeta^M_1$ and $\zeta^N_1$ are from (the dual of) \leref{lema}. Observe now the next diagram, where 
the upper 2-cell is the composite 2-cell from the left diagram above:

$$
\bfig
 \putmorphism(-400,900)(1,0)[M(A,B)`\amalg_{\substack{A\in\u\C \\ B\in\u\D}}M(A,B)`\iota^M_{A,B}]{980}1a
 \putmorphism(790,900)(1,0)[\phantom{F(B)}`\u{\C}\bullet 1 `m_1]{640}1a
 \putmorphism(1280,900)(1,0)[\phantom{(\u{\C}\times\u\D)\bullet 1}`\u{\E}\bullet 1 `F\bullet 1]{730}1a
\put(70,700){\fbox{left composite 2-cell from above}}

\putmorphism(-400,900)(0,-1)[\phantom{Y_2}``=]{450}1r
\putmorphism(2020,900)(0,-1)[\phantom{Y_2}``=]{400}1r

 \putmorphism(-400,500)(1,0)[`N(A',B')`f_{A,B}]{700}1a
 \putmorphism(500,500)(1,0)[`\amalg_{\substack{A'\in\u\E \\ B'\in\u\HH}} N(A',B')`\iota^N_{A',B'}]{820}1a
 \putmorphism(1560,500)(1,0)[\phantom{F(B)}` `n_1]{480}1a

\putmorphism(-400,520)(0,-1)[\phantom{Y_2}``=]{450}1l
\putmorphism(1370,520)(0,-1)[\phantom{Y_2}``=]{370}1l
\putmorphism(2020,520)(0,-1)[\phantom{Y_2}``=]{430}1r
\put(250,310){\fbox{$\zeta_{A,B}^{-1}$}}

 \putmorphism(-400,100)(1,0)[`\amalg_{\substack{A\in\u\C \\ B\in\u\D}}M(A,B)`\iota^M_{A,B}]{680}1a
 \putmorphism(600,100)(1,0)[`\amalg_{\substack{A'\in\u\E \\ B'\in\u\HH}} N(A',B')`\amalg f_{A,B}]{740}1a
 \putmorphism(1560,100)(1,0)[\phantom{F(B)}` `n_1]{480}1a
\efig
$$
and $\zeta_{A,B}^{-1}$ is a quasi-inverse of a 2-cell $\zeta_{A,B}$, which together with a 1-cell $\amalg f_{A,B}$, is induced 
by the 1-cell $\iota^N_{A',B'}f_{A,B}$ and the 3-coproduct property. Since this new composite 2-cell is an equivalence 2-cell, 
by the dual of \leref{equiv cells} we obtain an equivalence 2-cell $\alpha$, as desired. Doing the same as above, but projecting to the 
second coordinate, one obtains a desired equivalence 2-cell $\beta$. 

\bigskip


Let a 2-cell in $D_1$ be given, which in turn is given by a family of 3-cells (prisms) 
$$\Sigma_{A,B}: \frac{\nu_{A,B}}{[\xi_{A,B}\vert\Id_{[A',B']}]} \Rrightarrow \frac{[\Id_{[A,B]}\vert\crta\lambda\times\crta\rho]}{\nu'_{A,B}} $$
as in \equref{3-cell in matrix}. Let $\alpha, \beta$ and $\alpha', \beta'$ be the 2-cells from the image by $Int$ of the 1-cells given by 
$\nu_{A,B}$ and $\nu'_{A,B}$, respectively. 
Denote the total 2-cell in the left hand-side of \equref{Int-image of nu} by $\alpha_0$ and the 
analogous opposite 2-cell related to $\nu'_{A,B}$ by $\alpha_0'$. Concatenate the equivalence 2-cell $\zeta^{-1}_{A,B}$ (the one defining 
$\amalg f_{A,B}$) to $\alpha_0$, and analogously $(\zeta'_{A,B})^{-1}$ to $\alpha_0'$. Observe that $\xi_{A,B}: f_{A,B}\Rightarrow f_{A,B}'$  induces an equivalence 2-cell 
$\amalg\xi_{A,B}:\amalg f_{A,B}\Rightarrow\amalg f_{A,B}'$ and that moreover there is an invertible 3-cell 
$P_\zeta: \frac{\zeta^{-1}_{A,B}}{[\Id_{\iota^M}\vert\amalg\xi_{A,B}]} \Rrightarrow \frac{[\xi_{A,B}\vert\Id_{\iota^N}]}{(\zeta'_{A,B})^{-1}}$ 
- the both come from the dual of \coref{alfa-x-beta}. In \coref{cutting 3-cell-coproduct} set for $\xi$
to be $\crta\lambda\ot\Id_{m_1}:\crta F m_1\Rightarrow\crta F' m_1$ and for $\zeta$ in there to be 
$\Id_{n_1}\ot(\amalg\xi_{A,B}): n_1(\amalg f_{A,B})\Rightarrow n_1(\amalg f_{A,B}')$. Apart from the 3-cell $P_\zeta$, 
we also have all the prisms whose bases are the constituting 2-cells in $\alpha_0$ and whose opposite faces make $\alpha_0'$. Now by 
\coref{cutting 3-cell-coproduct} there is a unique 3-cell $\frac{\alpha}{[\amalg\xi_{A,B}\vert\Id_{n_1}]}\Rrightarrow\frac{[\Id_{m_1}\vert \crta\lambda]}{\alpha'}$. Similarly one obtains a unique 3-cell 
$\frac{\beta}{[\amalg\xi_{A,B}\vert\Id_{n_2}]}\Rrightarrow\frac{[\Id_{m_2}\vert \crta\rho]}{\beta'}$.

\subsection{Compatibility of $Int$ with 0-cells}

Before checking the compatibility of $Int$ with the composition of 1- and 2-cells, that we will not type here for the sake of saving space, 
one first needs to prove that there is 
a 1-cell $v: \amalg_{\substack{A\in\u\C \\ C\in\u\E}}\big(\amalg_{B\in\u\D}M(A,B)\times N(B,C)\big)
\to \big(\amalg_{\substack{A\in\u\C \\ B\in\u\D}}M(A,B)\big)\times_{\u\D\bullet 1}\big(\amalg_{\substack{B\in\u\D \\ C\in\u\E}} N(B,C)\big)$ 
in $V$. We show this as we will use it in the last Section. This also leads to a colax internal functor $\M\to\S$. 

To construct $v$ we need to find an equivalence 2-cell 
\begin{equation} \eqlabel{sigma za Int}
\sigma: n_1q_N\Rightarrow m_2q_M
\end{equation} 
(the total 2-cell in \equref{Int colax}), 
then apart from $v$ we will get also equivalence 2-cells $\lambda$ and $\rho$ as in \equref{Int colax} and an invertible 3-cell 
$\Sigma: \threefrac{\Id_{n_1}\ot\rho}{\omega\ot\Id_v}{\Id_{m_2}\ot\lambda} \Rrightarrow\sigma$. 
In the diagram we set for short $(\amalg M) \times_{\crta\D} (\amalg N) =
\big(\amalg_{\substack{A\in\u\C \\ B\in\u\D}}M(A,B)\big)\times_{\u\D\bullet 1}\big(\amalg_{\substack{B\in\u\D \\ C\in\u\E}} N(B,C)\big)$:
\begin{equation} \eqlabel{Int colax}
\scalebox{0.86}{\bfig
\putmorphism(-400,450)(1,0)[\amalg_{\substack{A\in\u\C \\ C\in\u\E}}\big(\amalg_{B\in\u\D}M(A,B)\times N(B,C)\big)``]{1100}0a
\putmorphism(-340,440)(1,-2)[` `q_M]{500}1l

\putmorphism(-180,450)(1,-1)[`(\amalg M) \times_{\crta\D} (\amalg N)` v]{420}1l 
\putmorphism(560,40)(1,0)[` \amalg_{\substack{B\in\u\D \\ C\in\u\E}} N(B,C)`p_2]{600}1a
\putmorphism(250,40)(0,-1)[` \amalg_{\substack{A\in\u\C \\ B\in\u\D}}M(A,B)`]{600}0l
\putmorphism(230,40)(0,-1)[` `]{560}1l
\putmorphism(220,60)(0,-1)[` `p_1]{600}0r
\putmorphism(1140,40)(0,-1)[` `n_1]{600}1r
\putmorphism(560,-550)(1,0)[`\u\D\bullet 1.`m_2]{560}1b
\put(250,220){\fbox{$\rho$}}
\putmorphism(40,460)(3,-1)[`` ]{1210}1r
\putmorphism(40,490)(3,-1)[`` q_N]{1210}0r
\put(0,-150){\fbox{$\lambda$}}
\put(680,-260){\fbox{$\omega$}}
\efig}
\end{equation} 
We first explain how to get $q_M$ and $q_N$. In the next diagram the 1-cell $\iota^M_Bp_1$ induces a 1-cell $h_M$ and an 
equivalence 2-cell $\zeta^M_1$. In turn, $\iota^M_Ah_M$ similarly induces $q_M$ and $\zeta^M_2$. Let us denote the composite 
equivalence 2-cell $(\zeta^M_1\vert\zeta^M_2)$ by $\zeta_{q_M}$. Analogous 1-cell $q_N$ and 
equivalence 2-cell $\zeta_{q_N}$ 
are obtained similarly, by projecting to the second coordinate on the most left. 
$$\scalebox{0.86}{
\bfig

\putmorphism(-1500,800)(1,0)[M(A,B)\times N(B,C)`\phantom{\amalg_{B\in\u\D}M(A,B)\times N(B,C)}`\iota^{B}]{1350}1a
\putmorphism(-150,800)(1,0)[\phantom{\amalg_{B\in\u\D}M(A,B)\times N(B,C)}``\iota^{A,C}]{950}1a
\putmorphism(-1500,800)(0,-1)[`` p_1]{400}1r
\putmorphism(-140,800)(0,-1)[\amalg_{B\in\u\D}M(A,B)\times N(B,C)`\amalg_{B\in\u\D}M(A,B)`  ]{400}0l
\putmorphism(-20,800)(0,-1)[``  h_M]{400}1l
\putmorphism(1560,800)(0,-1)[\amalg_{\substack{A\in\u\C \\ C\in\u\E}}\big(\amalg_{B\in\u\D}M(A,B)\times N(B,C)\big)` `q_M]{390}1l
\putmorphism(1560,760)(0,-1)[`     \amalg_{\substack{A\in\u\C \\ B\in\u\D}}M(A,B).`]{380}0l
\putmorphism(-1500,400)(1,0)[M(A,B)`\phantom{\amalg_{B'\in\u\HH}M'(A',B')\times N'(B',C')}`\iota^M_B]{1670}1a
\put(-990,610){\fbox{$\zeta^M_1$}}
\put(600,610){\fbox{$\zeta^M_2$}}
\putmorphism(-400,400)(1,0)[\phantom{\amalg_{B\in\u\D}M(A,B)\times N(B,C)}``\iota^M_A]{1640}1a
\efig}
$$
We set for short $f=\iota^M_A\iota^M_Bp_1, g=\iota^N_C\iota^N_B p_2$ and $\iota=\iota^{A,C}\iota^B$, then so far we have equivalence 
2-cells $\zeta_{q_M}$ and $\zeta_{q_N}$ in the diagram:
\begin{equation} \eqlabel{3-pullb}
\scalebox{0.86}{\bfig
\putmorphism(-130,450)(1,0)[\amalg_{\substack{A\in\u\C \\ C\in\u\E}}\big(\amalg_{B\in\u\D}M(A,B)\times N(B,C)\big)``]{1100}0a
\putmorphism(-420,440)(1,-2)[` `q_M]{480}1l

\putmorphism(-180,450)(1,-1)[`M(A,B)\times N(B,C)` \iota]{420}{-1}l 
\putmorphism(670,40)(1,0)[` \amalg_{\substack{B\in\u\D \\ C\in\u\E}} N(B,C)`g]{600}1a
\putmorphism(250,40)(0,-1)[` \amalg_{\substack{A\in\u\C \\ B\in\u\D}}M(A,B)`]{600}0l
\putmorphism(230,40)(0,-1)[` `]{560}1l
\putmorphism(220,60)(0,-1)[` `f]{600}0r
\putmorphism(1250,40)(0,-1)[` `n_1]{600}1r
\putmorphism(560,-550)(1,0)[`\u\D\bullet 1`m_2]{700}1b
\put(250,220){\fbox{$\zeta_{q_N}$}}
\putmorphism(160,460)(3,-1)[`` ]{1210}1r
\putmorphism(130,490)(3,-1)[`` q_N]{1210}0r
\put(-60,-150){\fbox{$\zeta_{q_M}$}}
\efig}
\end{equation}
and we define a 2-cell $n_1g\Rightarrow m_2f$ to be the following composite equivalence 2-cell:
$$
\scalebox{0.86}{\bfig
 \putmorphism(-350,900)(1,0)[M(A,B)\times N(B,C)`N(B,C)`p_2]{850}1a
 \putmorphism(550,900)(1,0)[\phantom{F(B)}`\amalg_{\substack{B\in\u\D \\ C\in\u\E}}N(B,C) `\iota^N_C\iota^N_B]{730}1a
 \putmorphism(1600,900)(1,0)[`\u{\D}\bullet 1 `n_1]{450}1a
\put(1240,700){\fbox{$\zeta_{n_1}$}}

\putmorphism(-350,900)(0,-1)[\phantom{Y_2}``=]{400}1r
\putmorphism(520,900)(0,-1)[\phantom{Y_2}``=]{400}1r
\putmorphism(2040,900)(0,-1)[\phantom{Y_2}``=]{400}1r

 \putmorphism(-350,500)(1,0)[``p_2]{860}1a
 \putmorphism(500,500)(1,0)[``!]{730}1a
 \putmorphism(1230,500)(1,0)[` `\sigma_B]{800}1a

\putmorphism(-350,520)(0,-1)[\phantom{Y_2}``=]{450}1l
\putmorphism(1220,520)(0,-1)[\phantom{Y_2}``=]{450}1l
\putmorphism(2040,520)(0,-1)[\phantom{Y_2}``=]{450}1r
\put(450,280){\fbox{$\kappa$}}

 \putmorphism(-350,100)(1,0)[``!]{1560}1a
 \putmorphism(1110,100)(1,0)[\phantom{F(B)}` `\sigma_B]{900}1a

\putmorphism(-350,120)(0,-1)[\phantom{Y_2}``=]{430}1l
\putmorphism(1220,120)(0,-1)[\phantom{Y_2}``=]{430}1r
\putmorphism(2040,120)(0,-1)[\phantom{Y_2}``=]{450}1r
\put(200,-120){\fbox{$\kappa^{-1}$}}

 \putmorphism(-350,-280)(1,0)[`M(A,B)`p_1]{740}1a
 \putmorphism(390,-280)(1,0)[\phantom{M(A,B)}` `!]{800}1a
 \putmorphism(1110,-280)(1,0)[\phantom{F(B)}` `\sigma_B]{900}1a

\putmorphism(-350,-260)(0,-1)[\phantom{Y_2}``=]{430}1l
\putmorphism(380,-260)(0,-1)[\phantom{Y_2}``=]{430}1l
\putmorphism(2040,-260)(0,-1)[\phantom{Y_2}``=]{430}1r
\put(1050,-480){\fbox{$\zeta_{m_2}$}}

 \putmorphism(-350,-670)(1,0)[``p_1]{760}1a
 \putmorphism(260,-670)(1,0)[\phantom{F(B)}`\amalg_{\substack{A\in\u\C \\ B\in\u\D}}M(A,B)`\iota^M_A\iota^M_B]{1000}1a
 \putmorphism(1480,-670)(1,0)[\phantom{F(B)}` `m_2]{570}1a
\efig}
$$
where $\zeta_{n_1}$ and $\zeta_{m_2}$ are the obvious equivalence 2-cells. 
Now we set the composite equivalence 2-cell $n_1q_N\Rightarrow m_2q_M$ in \equref{3-pullb} to be the desired 2-cell $\sigma$ in \equref{sigma za Int}. 
Finally, by the 3-pullback property of  
$\big(\amalg_{\substack{A\in\u\C \\ B\in\u\D}}M(A,B)\big)\times_{\u\D\bullet 1}\big(\amalg_{\substack{B\in\u\D \\ C\in\u\E}} N(B,C)\big)$, 
we get a 1-cell $v$ and equivalence 2-cells $\lambda$ and $\rho$ in \equref{Int colax} and an isomorphism 3-cell $\Sigma$, as announced.

\section{Equivalence of matrices and spans in a tricategory}

In this final Section we examine equivalence conditions for the $(1\times 2)$-categories of matrices $\M$ and spans $\S$ in 
a 1-strict tricategory $V$ 
to be equivalent, and also for the 1-categories of discretely internal and enriched categories in $V$ to be equivalent. Inspired by the ideas that we exposed in 
\seref{Power} we start by introducing monads in $(1\times 2)$-categories and then summarize our findings in the tricategorical setting.

\subsection{Monads in $(1\times 2)$-categories} \sslabel{3-double monads}

In this Subection we are going to introduce monads and vertical monad morphisms in a $(1\times 2)$-category $\Vv$. 
In the analogy with the definition of a monad in a double category \cite[Definition 2.4]{FGK} we introduce:

\begin{defn}
A monad in a $(1\times 2)$-category $\Vv$ is a monad in the horizontal tricategory $\HH(\Vv)$ of $\Vv$ (see \deref{3-monad}). 
\end{defn}

Whereas:

\begin{defn} \delabel{3-monad}
A monad in a tricategory $V$ is given by a 1-endocell $T:\A\to\A$ with two 2-cells $\mu:TT\Rightarrow T$ and $\eta:\Id_\A\Rightarrow T$ 
and 3-cells  $\alpha:\frac{\Id_T\ot\mu}{\mu}\Rrightarrow\frac{\mu\ot\Id_T}{\mu}, 
\lambda:\frac{\Id_T\ot\eta}{\mu}\Rrightarrow\Id_T, \rho:\frac{\Id_T\ot\eta}{\mu}\Rrightarrow\Id_T$ 
which satisfy the usual five axioms that expressed in terms of equations of the transversal compositions of 3-cells have the form:
$$\frac{\alpha\ot\Id}{\Id}\cdot\frac{\Id}{\alpha}\cdot\frac{\Id\ot\alpha}{\Id}=\frac{\Id}{\alpha}\cdot\frac{inter}{\Id}\cdot\frac{\Id}{\alpha}$$
$$\frac{\Id\ot\lambda}{\Id}\cdot\frac{\Id}{\alpha}=\frac{\rho}{\Id}$$
$$\frac{\Id}{\lambda}\cdot inter\cdot\frac{\Id}{\alpha}=\frac{\lambda}{\Id}$$
$$\frac{\Id}{\rho}\cdot inter=\frac{\rho\ot\Id}{\Id}\cdot\frac{\Id}{\alpha}$$
$$\frac{\Id}{\rho}\cdot inter=\frac{\Id}{\lambda}.$$ 
\end{defn}

We are interested in monads in the $(1\times 2)$-categories of matrices $\M$ and spans $\S$ in a 1-strict tricategory $V$. Being monads in their 
respective horizontal tricategories $\HH(\M)$ and $\HH(\S)$, observe that their 3-cells $\alpha, \lambda, \rho$ for associativity and unitality 
are given through both 2-cells and 3-cells in $V$. We explain now how 
the five axioms for $\alpha, \lambda, \rho$ come down to 3-cells analogous to those in \cite[Definition 6.2, 4)]{Fem}. 
We will restrict to a particular kind of monads. Namely, we consider those monads 
in $\HH(\S)$ whose 2-cells $\mu$ are given by two identity 2-cells and identity vertical 1-cells ({\em i.e.} identity functors), 
see \equref{1-cells in C_1}. We will refer to such monads {\em strict monads}. Let a strict monad be given by a cospan 
$\crta\C\stackrel{s}{\longleftarrow}T\stackrel{t}{\longrightarrow}\crta\C$ and a 1-cell $c:T\times_{\u\C}T\to T$ in $V$ determining the 2-cell $\mu$ 
for the monad. Then we have $sc=sp_1$ and $tc=tp_2$ 
(we restrict to strict monads precisely in order to have the latter identities hold strictly, for the purpose of \prref{monad-int}). 
The associativity 3-cell $\alpha$ is then given by (a pair of prisms determined by) a 2-cell $a^*:c(\Id_T\times_{\u\C} c)\Rightarrow c(c\times_{\u\C}\Id_T)$ and 
a pair of 3-cells of the form $\frac{\Id_{(sp_1)p_1=sc(c\times\Id)}}{\Id_s\ot a*}\Rrightarrow \Id_{\Id_{(sp_1)p_1=sc(\Id\times c)}}$ and 
an analogous one for $t$. To study the first of the five axioms above, observe that each of the 3-cells that are being composed in the equation 
comes down to the horizontal composition of $a^*$ and the relating identity 2-cells ($\ot$ becomes $\times_{\u\C}$ and $\frac{\bullet}{\bullet}$ becomes 
horizontal composition of 2-cells in $V$), while the transversal composition of those 3-cells comes down to the vertical composition of 
the obtained 2-cells. Thus the first of the five axioms means that two pairs of 3-cells (prisms) 
$(P^\Lambda_s, P^\Lambda_t)$ and $(P^P_s, P^P_t)$ are equal so that $P^\Lambda_s=P^P_s$ 
yields that two 3-cells of the form $\Omega_\Lambda:\Id_s\ot \Lambda\Rrightarrow \Id$ and $\Omega_P:\Id_s\ot P\Rrightarrow \Id$ are equal, 
where 
$$\Lambda= \threefrac{\Id_c \ot(\Id_{id_T}\times_{\u\C} a^*)}{a^*\ot\Id_{1\times_{\u\C} c\times_{\u\C} 1}}{\Id_c \ot (a^*\times_{\u\C} \Id_{id_T})}\quad
\text{and}\quad
P=\threefrac{a^* \ot\Id_{1\times_{\u\C} 1\times_{\u\C} c}}{\Id_c \ot \Nat}{a^* \ot\Id_{c\times_{\u\C} 1\times_{\u\C} 1}}.$$
This implies that the 2-cells $\Id_s\ot \Lambda$ and $\Id_s\ot P$ are equal, yielding 
a bijective 3-cell $\pi^*: \Lambda\Rrightarrow P$. 
(\prref{2-cell-prop} gives a hint that it is sufficient to consider the equality of components related to $s$.) 
The similar reasoning is applied to the other four axioms, and one finds that they induce bijective 3-cells $\mu^*, \lambda^*, \rho^*, \Epsilon^*$, 
respectively, which satisfy the axioms from \cite[Definition 6.2, 4)]{Fem}. 

We do the same for matrices in $V$: we consider strict monads in matrices in $V$ and come to analogous conclusions. A {\em strict monad} 
in $\HH(\M)$ is a monad in $\HH(\M)$ for which $F$ and $G$ are identities and for every $A,B\in\C$ it is 
$\nu_{A,B}\iso\Id_{\Id_{\C\times\C}}\ot\kappa_{f_{A,B}}$ and $\chi^{F,G}_{A,B}$ are identities, see \equref{nu}.

Analogously to the well-known fact that monads in the bicategories of matrices and spans in a 1-category $\C$ which has pullbacks, products and coproducts 
are categories enriched, respectively internal, in $\C$, we have:

\begin{prop} \prlabel{monad-int}
A strict monad in $\S$, the $(1\times 2)$-category of spans in a 1-strict tricategory $V$, is a category discretely internal in $V$ in the sense of 
\cite[Definition 6.2]{Fem}. 

A strict monad in $\M$, the $(1\times 2)$-category of matrices in a 1-strict tricategory $V$, is a category enriched in $V$ in the sense of 
\cite[Definition 8.1]{Fem}. 
\end{prop}

Discretely internal here means that the object of objects is a 3-coproduct of copies of the terminal object. 

Rather than defining a $(1\times 2)$-category of monads in a $(1\times 2)$-category $\Vv$, in analogy to the double category of monads in a 
double category from \cite{FGK}, for simplicity reasons we restrict ourselves to defining only the {\em vertical morphisms} of monads in $\Vv$.

\begin{defn} 
A vertical morphism between monads $T:\A\to\A$ and $L:\A'\to\A'$ in a $(1\times 2)$-category $\Vv$ is a horizontal 2-cell $\delta$ as below together with 3-cells:
$$\scalebox{0.86}{
\bfig
\putmorphism(0,300)(1,0)[``T]{450}1a
\putmorphism(460,300)(1,0)[``T]{450}1a
\putmorphism(0,320)(0,-1)[\phantom{Y_2}``=]{400}1l
\putmorphism(890,320)(0,-1)[\phantom{Y_2}``=]{400}1r
\put(360,130){\fbox{$\mu_T$}}
\putmorphism(0,-20)(1,0)[``T]{900}1b
\putmorphism(0,-10)(0,-1)[\phantom{Y_2}``u]{400}1l
\putmorphism(890,-10)(0,-1)[\phantom{Y_2}``u]{400}1r
\put(380,-240){\fbox{$\delta$}}
\putmorphism(0,-350)(1,0)[``L]{900}1b
\efig}
\quad\stackrel{m^*}{\Rrightarrow}\quad
\scalebox{0.86}{\bfig
\putmorphism(0,300)(1,0)[``T]{450}1a
\putmorphism(460,300)(1,0)[``T]{450}1a
\putmorphism(0,340)(0,-1)[\phantom{Y_2}``u]{400}1l
\putmorphism(450,340)(0,-1)[\phantom{Y_2}``u]{400}1r
\putmorphism(890,340)(0,-1)[\phantom{Y_2}``u]{400}1r
\put(150,110){\fbox{$\delta$}}
\put(600,110){\fbox{$\delta$}}
\putmorphism(0,-30)(1,0)[``L]{450}1b
\putmorphism(460,-30)(1,0)[``L]{450}1b
\putmorphism(0,-10)(0,-1)[\phantom{Y_2}``=]{400}1l
\putmorphism(890,-10)(0,-1)[\phantom{Y_2}``=]{400}1r
\put(380,-200){\fbox{$\mu_L$}}
\putmorphism(0,-350)(1,0)[``L]{900}1b
\efig}
$$
and
$$
\scalebox{0.86}{\bfig
\putmorphism(0,350)(1,0)[``=]{450}1a
\putmorphism(0,370)(0,-1)[\phantom{Y_2}``=]{400}1l
\putmorphism(450,370)(0,-1)[\phantom{Y_2}``=]{400}1r
\put(150,150){\fbox{$\eta_T$}}
\put(150,-210){\fbox{$\delta$}}
\putmorphism(0,10)(0,-1)[\phantom{Y_2}``u]{400}1l
\putmorphism(450,10)(0,-1)[\phantom{Y_2}``u]{400}1r
\putmorphism(0,0)(1,0)[``T]{450}1b
\putmorphism(0,-360)(1,0)[``L]{450}1b
\efig}
\quad\stackrel{i^*}{\Rrightarrow}\quad
\scalebox{0.86}{\bfig
\putmorphism(0,350)(1,0)[``=]{450}1a
\putmorphism(0,370)(0,-1)[\phantom{Y_2}``u]{400}1l
\putmorphism(450,370)(0,-1)[\phantom{Y_2}``u]{400}1r
\put(150,150){\fbox{$\Id$}}
\put(150,-210){\fbox{$\eta_L$}}
\putmorphism(0,10)(0,-1)[\phantom{Y_2}``=]{400}1l
\putmorphism(450,10)(0,-1)[\phantom{Y_2}``=]{400}1r
\putmorphism(0,0)(1,0)[``=]{450}1b
\putmorphism(0,-360)(1,0)[``L]{450}1b
\efig}
$$
satisfying 
the following axioms which we express in terms of equations of the transversal compositions of 3-cells:
$$\frac{\Id}{\alpha}\cdot\frac{\xi^{-1}(\Id\ot m^*)\xi}{\Id}\cdot\frac{\Id}{m^*}=\frac{\xi^{-1}(m^*\ot\Id)\xi}{\Id}\cdot
\frac{\Id}{m^*}\cdot\frac{\alpha}{\Id}$$
$$l\cdot\frac{\Id}{\lambda}\cdot\frac{\xi^{-1}(\Id\ot i^*)\xi}{\Id}=\frac{\lambda}{\Id}\cdot\frac{\Id}{(m^*)^{-1}}$$
$$r\cdot\frac{\Id}{\rho}\cdot\frac{\xi^{-1}(i^*\ot\Id)\xi}{\Id}=\frac{\rho}{\Id}\cdot\frac{\Id}{(m^*)^{-1}}.$$
Here $\xi$ stands for the interchange 3-cell, and $l,r$ are left and right unity constraints for the horizontal composition of 2-cells in $V$. 
\end{defn}

We denote the category of monads and their vertical morphisms in a $(1\times 2)$-category $\Vv$, with the vertical composition of 
horizontal 2-cells in $\Vv$ (that is, of vertical monad morphisms $\delta$), by $\Mnd(\Vv)$. 

Analogously to \prref{equiv monads} 
the following is straightforward to prove:

\begin{prop} \prlabel{lift to monads}
For two equivalent $(1\times 2)$-categories $\Vv_1$ and $\Vv_2$ their respective categories of monads $\Mnd(\Vv_1)$ and $\Mnd(\Vv_2)$ are equivalent. 
\end{prop}

\subsection{Equivalence of $(1\times 2)$-categories of matrices and spans and of discretely internal and enriched categories in a tricategory}

In analogy to bicategorical biequivalence functors, a trifunctor is a triequivalence if and only if it is pseudo by nature ({\em i.e.} 
it is compatible with the composition of 
1-cells up to an equivalence 2-cell), it is essentially surjective on the class of objects and its every component 
bicategorical functor is an equivalence pseudofunctor (a biequivalence). 

By their construction, the internal categories $\M$ and $\S$ are equivalent in any of the two following cases:\\
$\bullet$ the bicategories $C_1$ and $D_1$ from Sections 5 and 6 
are biequivalent; \\
$\bullet$ the horizontal tricategories $\HH(\M)$ and $\HH(\S)$ are triequivalent. \\
For the second case we may consider the trifunctor $\I:\HH(\M)\to\HH(\S)$ which is identity on 0-cells and on 
hom-bicategories for fixed 0-cells, which are small categories $\C$ and $\D$, 
consider the clear restriction $Int_{\C,\D}: V\x\Mat(\C,\D)\to\Span_d(V)(\C,\D)$ of the pseudofunctor $Int: D_1\to C_1$ to the 
hom-bicategories of $\HH(\M)$ and $\HH(\S)$, {\em i.e.} the obvious sub-bicategories $V\x\Mat(\C,\D)$ of $D_1$ and $\Span_d(V)(\C,\D)$ of 
$C_1$. By the above observation $\I$ is a triequivalence if and only if for all 
small categories \(\D\) the 1-cell $v$ in \equref{Int colax} is a biequivalence 1-cell in $V$ and the pseudofunctors $Int_{\C,\D}$ 
are biequivalences for all small categories $\C$ and $\D$.


We next study when $v$ is a biequivalence 1-cell. 
Let us consider the following sub-tricategories of $\HH(\M)$ and $\HH(\S)$. First consider the sub-tricategory of $\HH(\M)$ in which 
at all cell levels matrices indexed over pairs of small categories $(\C,\D)$ are replaced by matrices indexed over pairs $(*,\D)$, that is, 
lists indexed over small categories $\D$. 
Correspondingly, all higher cells on pairs $(\C,\D)$ are replaced by identity higher cells over categories $\D$. Next, take a 
sub-tricategory of the latter sub-tricategory, where we fix a single 0-cell $\D$. We denote this tricategory by $V^\D$. Importantly, 
observe that 
its only hom-bicategory is the sub-bicategory $V\x\Mat(*,\D)$ of $D_1$. 

Fully analogously to $V^\D$, let $V / (\u\D \bullet 1)$ be the sub-tricategory of $\HH(\S)$ with a fixed 0-cell $\D$ and whose 
only hom-bicategory is the sub-bicategory $\Span_d(V)(* ,\D)$ of $C_1$. 


Let us now consider the trifunctor 
\begin{equation} \eqlabel{perp-3}
\amalg : V^{\u\D} \to V / (\u\D \bullet 1) 
\end{equation}
that is identity on the unique 0-cell $\D$ and on hom-bicategories set 
$$\amalg_{\D,\D}=Int_{*,\D}:V\x\Mat(*,\D)\to\Span_d(V)(*,\D),$$ 
the restriction of the pseudofunctor $Int: D_1\to C_1$.    

\begin{rem}
Saying that the trifunctor \(\coprod : V^{\u\D} \to V / (\u\D \bullet 1)\) has a property P for every small category \(\D\) 
is the same as saying that the trifunctor \(\coprod : V^{\u\C\times\u\D} \to V / ((\u\C\times\u\D) \bullet 1)\) has it 
for all small categories $\C$ and $\D$ (replace $\D$ by $\C\times\D$ in one direction, and $\C$ by the trivial category *, in the other). 

If $V$ is 3-Cartesian closed the trifunctors $X\times-$, $-\times X$ preserve 3-coproducts, for any object $X$ of $V$. 
Then there is a natural biequivalence 1-cell $\phi:\u{\C} \bullet 1 \times \u{\D}\bullet 1 \to (\u{\C}\times \u{\D})\bullet 1$ 
in $V$ with (naturality) equivalence 2-cells $\Phi$ below for all functors $F,G$: 
$$\bfig
\putmorphism(950,400)(0,-1)[\u{\C}\bullet 1\times \u{\D}\bullet 1`\u{\E}\bullet 1\times \u{\HH}\bullet 1`]{400}0r
\putmorphism(770,400)(0,-1)[``F\bullet 1\times G\bullet 1]{400}1l
\putmorphism(960,400)(1,0)[\phantom{(\u{\C}\times \u{\D})\bullet 1}``\phi]{600}1a
\putmorphism(1800,400)(0,-1)[(\u{\C}\times \u{\D})\bullet 1`(\u{\E}\times \u{\HH})\bullet 1.`]{400}0r
\putmorphism(1720,400)(0,-1)[``(F\times G)\bullet 1]{400}1r
\putmorphism(980,0)(1,0)[\phantom{(\u{\C}\times \u{\D})\bullet 1}``\phi']{570}1a
\put(1300,200){\fbox{$\Phi$}}
\efig 
$$
In this case by \prref{2-cell-prop} we obtain that there is a biequivalence of bicategories $\Span_d(V)(\C,\D)\simeq\Span_d(V)(*,\C\times\D)$, 
being the latter the hom-bicategory of the (sub-)tricategory $V/((\u\C\times\u\D) \bullet 1)$ 
(concatenate the above equivalence 2-cells $\Phi$ to the 2-cell $\gamma$ in \prref{2-cell-prop}). On the other hand, it is clear that $V\x\Mat(*,\C\times\D)
\simeq V\x\Mat(\C,\D)$. Then we may 
observe that the trifunctor \(\coprod : V^{\u\C\times\u\D} \to V / ((\u\C\times\u\D) \bullet 1)\) on hom-bicategories is given indeed by 
$Int_{1,\C\times\D}=Int_{\C,\D}$. 
\end{rem}

Due to the Remark we may state: 

\begin{prop}
If $V$ is 3-Cartesian closed 
the trifunctor \(\coprod : V^{\u\D} \to V / (\u\D \bullet 1)\) is a triequivalence for all \(\D\) if and only if it is ``pseudo'' and 
the pseudofunctors $Int_{\C,\D}$ are biequivalences of bicategories for all small categories $\C$ and $\D$. 
\end{prop}

The following result is a tricategorification of \cite[Proposition 3.1]{CFP}. 

\begin{prop} \prlabel{perp preserves}
Let \(V\) be 3-Cartesian closed and assume that for every small category \(\D\), the trifunctor \(\coprod : V^{\u\D} \to V / (\u\D \bullet 1)\) preserves 
binary 3-products. Then the 1-cell $v$ in \equref{Int colax} is a biequivalence in $V$ (and consequently, the trifunctor 
$\I:\HH(\M)\to\HH(\S)$ is ``pseudo''). 
\end{prop}

\begin{proof}
Binary 3-products in \(V^{\u\D}\) are given pointwise,
while binary 3-products in \(V / (\u\D \bullet 1)\) are given by 3-pullback.
By the Cartesian closedness of $V$ the 3-product trifunctors $X\times-$ and $-\times X$ commute with 3-coproducts. Thus the 3-coproduct 
$\amalg_{\substack{A\in\u\C \\ C\in\u\E}}\big(M(A,B)\times N(B,C)\big)$ 
can be seen as a 3-product 
$\big(\amalg_{\substack{A\in\u\C}}M(A,B)\big)\times\big(\amalg_{C\in\u\E} N(B,C)\big)$, i.e. a 3-pullback over 1.  
When the trifunctor $\amalg_{B\in\u\D}$ acts on it, by assumption it sands it to a 3-product, 
which in $V / (\u\D \bullet 1)$ is a 3-pullback over $\D$. This means that the outer arrows in \equref{Int colax} 
denote a 3-pullback, while obviously the inner square diagram in there denotes a 3-pullback for the same cospan in $V$. 
%
Since 3-limits are unique up to biequivalence 1-cells (recall \rmref{uniqueness}),
the comparison 1-cell \(v\) is a biequivalence in \(V\), making \(\I\) into a pseudo trifunctor. 
\end{proof}

For the converse we may assume less:

\begin{prop} \prlabel{perp preserves - converse}
Let \(V\) be 3-Cartesian closed and assume that for every small category \(\D\) the 1-cell $v: \amalg_{B\in\u\D}M(*,B)\times N(B,*)
\to \big(\amalg_{\substack{B\in\u\D}}M(*,B)\big)\times_{\u\D\bullet 1}\big(\amalg_{\substack{B\in\u\D}} N(B,*)\big)$ 
(a special case of \equref{Int colax}) is a biequivalence in $V$ (and consequently, the trifunctor 
\(\coprod : V^{\u\D} \to V / (\u\D \bullet 1)\)  is ``pseudo''). 
Then the trifunctor \(\coprod : V^{\u\D} \to V / (\u\D \bullet 1)\) preserves binary 3-products. 
\end{prop}

\begin{proof}
The proof is fully analogous to the direct direction of \cite[Proposition 3.1]{CFP}. 
\end{proof}

\begin{cor}
In a 3-Cartesian closed 1-strict tricategory $V$ with terminal object, 3-(co)products and 3-pullbacks, the following are equivalent:
\begin{enumerate}
\item the trifunctor \(\coprod : V^{\u\D} \to V / (\u\D \bullet 1)\) is ``psuedo'';
\item the trifunctor $\I:\HH(\M)\to\HH(\S)$ is ``psuedo'';
\item the trifunctor \(\coprod : V^{\u\D} \to V / (\u\D \bullet 1)\) preserves binary 3-products.  
\end{enumerate} 
\end{cor}

As a byproduct of the above Corollary one obtains that  
the 1-cell \vspace{-0,3cm}
$$v: \amalg_{\substack{A\in\u\C \\ C\in\u\E}}\big(\amalg_{B\in\u\D}M(A,B)\times N(B,C)\big)
\to \big(\amalg_{\substack{A\in\u\C \\ B\in\u\D}}M(A,B)\big)\times_{\u\D\bullet 1}\big(\amalg_{\substack{B\in\u\D \\ C\in\u\E}} N(B,C)\big)
\vspace{-0,3cm}$$ 
from \equref{Int colax} is a biequivalence for all $\C,\D,\E$ if and only if so is its special case 
1-cell $v: \amalg_{B\in\u\D}M(*,B)\times N(B,*)
\to \big(\amalg_{\substack{B\in\u\D}}M(*,B)\big)\times_{\u\D\bullet 1}\big(\amalg_{\substack{B\in\u\D}} N(B,*)\big)$ for every $\D$. 

From all the above said we obtain:

\begin{cor} \colabel{mainThm}
In a 3-Cartesian closed 1-strict tricategory \(V\) with terminal object, 3-(co)products and 3-pullbacks. The following are equivalent:
 \begin{enumerate}
\item the trifunctor \(\coprod : V^{\u\D} \to V / (\u\D \bullet 1)\) is a triequivalence for all \(\D\); 
\item the trifunctor $\I:\HH(\M)\to\HH(\S)$ is a triequivalence; 
\item the colax internal functor $\M\to\S$ (constructed in Subsections 7.3 and 7.4) is an equivalence of $(1\times 2)$-categories. 
 \end{enumerate}
\end{cor}

\begin{rem}
Analogously to $\I:\HH(\M)\to\HH(\S)$ we may consider the trifunctor $\E:\HH(\S)\to\HH(\M)$ which is identity on 0-cells and on 
hom-bicategories the clear restriction $En_{\C,\D}: \Span_d(V)(\C,\D)\to V\x\Mat(\C,\D)$ of the pseudofunctor $En: C_1\to D_1$. 
By analogy to \cite[Proposition 2.2]{CFP} one has that the trifunctors $\I:\HH(\M)\to\HH(\S)$ and $\E:\HH(\S)\to\HH(\M)$ 
are 3-adjoint. Then the following two equivalent statements can be added as two additional equivalent conditions in the above Corollary: 
$\E:\HH(\S)\to\HH(\M)$ is a triequivalence if and only if the lax internal functor $\S\to\M$ (constructed in Subsection 7.2) 
is an equivalence of $(1\times 2)$-categories. 
\end{rem}


In any of the equivalent conditions of the above Corollary, 
by \prref{lift to monads} the categories of monads of $\M$ and $\S$ are equivalent. As a matter of fact, analogously as 
lax functors of monoidal categories preserve monoids, any lax trifunctor preserves monads in tricategories. Hence $\E$ preserves monads, 
and $\I$ does so if it is ``pseudo'' (for example in the conditions of \prref{perp preserves}). Consequently, in exactly the same conditions 
the internal functors $\M\to\S$ and $\S\to\M$ preserve monads. 
It is easily and directly proved that $\E$ and $\I$ (the latter being ``pseudo'') preserve strict monads. Thus we have:

\begin{prop} \prlabel{I preserves monads}
Under conditions of \prref{perp preserves} the trifunctor $\I$ preserves strict monads, {\em i.e.} due to \prref{monad-int} 
it maps categories enriched in $V$ into categories discretely internal in $V$. 
\end{prop}

Moreover, let us restrict to {\em strict vertical morphisms of monads in $\S$} - those for which the equivalence 2-cells $\alpha$ 
and $\beta$ in  \equref{1-cells in C_1} are identities, and to {\em strict vertical morphisms of monads in $\M$} - those for which 
$\nu_{A,B}\iso\Id_{\Id_{\C\times\C}}\ot\kappa_{f_{A,B}}$ and $\chi^{F,G}_{A,B}$ are identities for all $A,B\in\u\C$ in \equref{nu}. 
Then it is not difficult to see that strict vertical morphisms of monads in $\S$ yield internal functors in $V$  
and strict vertical morphisms of monads in $\M$ yield enriched functors in $V$, recall \ssref{int-enrich-fun}. 
 
Now by \prref{monad-int} and the latter we finally obtain:

\begin{cor} \colabel{equiv-monads}
Let \(V\) be 3-Cartesian closed with terminal object, 3-(co)products and 3-pullbacks and assume that for every small category \(\D\), 
the trifunctor \(\coprod : V^{\u\D} \to V / (\u\D \bullet 1)\) is a triequivalence. Then the 
subcategories of strict monads and strict vertical morhisms of monads $\Mnd_s(\M)$ and $\Mnd_s(\S)$ of $\Mnd(\M)$ and $\Mnd(\S)$, 
respectively, are equivalent. Equivalently, the categories of discretely internal and enriched categories in $V$  
are equivalent.
\end{cor}

\medskip

We finish the paper with a comparison of our results with those from \cite[Section 8]{Fem}. 
We assume Cartesian closedness in \prref{perp preserves} in order to have that the  trifunctors $X\times-$ and $-\times X$ commute with 3-coproducts, 
which is the second preservation assumption in \cite[Proposition 8.4]{Fem}. That the trifunctor \(\coprod : V^{\u\D} \to V / (\u\D \bullet 1)\) preserves 
binary 3-products can be restated so that the 3-coproduct functor maps binary 3-products, which are special kind of 3-pullbacks, into 3-pullbacks in $V$. 
Thus requiring that the 3-coproduct functor preserves 3-pullbacks, which is the first preservation assumption in \cite[Proposition 8.4]{Fem}, is a little bit 
stronger than the latter. The result of \cite[Proposition 8.4]{Fem} concerns an enriched category $\Tau$ in $V$ and $n$-ary 3-products of its endo-hom 0-cells 
in $V$, so that for $n=2$ it claims the existence of biequivalence 1-cells 
\begin{multline*}
\amalg_{B\in Ob\Tau}\big((\amalg_{A\in Ob\Tau}\Tau(A,B))\times(\amalg_{C\in Ob\Tau}\Tau(B,C))\big)\to \\
 (\amalg_{A,B\in Ob\Tau}\Tau(A,B))\times_{\substack{(\amalg_{B\in Ob\Tau}1_B)}}(\amalg_{B,C\in Ob\Tau}\Tau(B,C)).
\end{multline*}
Thus it is a special case of our \prref{perp preserves}. 
Finally, under the above assumptions in \cite[Proposition 8.5]{Fem} it is proved that an enriched category in $V$ is an internal category in $V$, which is 
our \prref{I preserves monads}. It is clear that our present constructions give a much broader picture than the approach that we employed in 
\cite[Section 8]{Fem}, for which we followed \cite{Ehr}.

\bigskip

\bigskip

{\bf Acknowledgmentes.} 
The first author was supported by the Science Fund of the Republic of Serbia, Grant No. 7749891, Graphical Languages - GWORDS. 

\bigskip

\section{Appendix - Axioms for 3-cells of an internal functor}

In these diagrams the symbol $\times$ will stand for short for the pullback $-\times_{C_0}-$ or $-\times_{D_0}-$. 

(A1) 
\[
\adjustbox{scale=0.7,center}{%
\begin{tikzcd}[row sep = huge, column sep = huge]
{}
&{F(R) \hspace{-0,06cm}\times\hspace{-0,03cm} (F(S) \hspace{-0,06cm}\times\hspace{-0,03cm} F(T \hspace{-0,06cm}\times\hspace{-0,03cm} U))}
	\ar[r, Rightarrow, "{\Id_{F(R)} \hspace{-0,06cm}\times\hspace{-0,03cm} F_{\times}^{S, T \times U}}"]
	\ar[ddr, triple, "{\Id_{F(R)} \hspace{-0,06cm}\times\hspace{-0,03cm} \Omega_{\alpha}^{S, T, U}}" description]
&{F(R) \hspace{-0,06cm}\times\hspace{-0,03cm} F(S \hspace{-0,06cm}\times\hspace{-0,03cm} (T \hspace{-0,06cm}\times\hspace{-0,03cm} U))}
	\ar[r, Rightarrow, "{F_{\times}^{R, S\times (T\times U)}}"]
	\ar[dr, phantom, "{\text{nat.}}" description]
&{F(R \hspace{-0,06cm}\times\hspace{-0,03cm} (S \hspace{-0,06cm}\times\hspace{-0,03cm} (T \hspace{-0,06cm}\times\hspace{-0,03cm} U)))} \\
{(F(R) \hspace{-0,06cm}\times\hspace{-0,03cm} F(S)) \hspace{-0,06cm}\times\hspace{-0,03cm} (F(T) \hspace{-0,06cm}\times\hspace{-0,03cm} F(U))}
	\ar[r, Rightarrow, "{\alpha^*_{F(R), F(S), F(T \times U)}}"]
	\ar[ddr, triple, "{\pi^*}" description]
&{F(R) \hspace{-0,06cm}\times\hspace{-0,03cm} (F(S) \hspace{-0,06cm}\times\hspace{-0,03cm} (F(T) \hspace{-0,06cm}\times\hspace{-0,03cm} F(U)))}
	\ar[u, Rightarrow, "{\Id_{F(R)} \hspace{-0,06cm}\times\hspace{-0,03cm} (\Id_{F(S)} \hspace{-0,06cm}\times\hspace{-0,03cm} F_{\times}^{T, U})}"]
&{F(R) \hspace{-0,06cm}\times\hspace{-0,03cm} F((S \hspace{-0,06cm}\times\hspace{-0,03cm} T) \hspace{-0,06cm}\times\hspace{-0,03cm} U)}
	\ar[u, Rightarrow, "{\Id_{F(R)} \hspace{-0,06cm}\times\hspace{-0,03cm} F(\alpha^*_{S, T, U})}" description]
	\ar[r, Rightarrow, "{F_{\times}^{R, (S \times T) \hspace{-0,06cm}\times\hspace{-0,03cm} U}}"]
&{F(R \hspace{-0,06cm}\times\hspace{-0,03cm} ((S \hspace{-0,06cm}\times\hspace{-0,03cm} T) \hspace{-0,06cm}\times\hspace{-0,03cm} U))}
	\ar[u, Rightarrow, "{F(\Id_R \hspace{-0,06cm}\times\hspace{-0,03cm} \alpha^*_{S, T, U})}"'] \\
{}
&{F(R) \hspace{-0,06cm}\times\hspace{-0,03cm} ((F(S) \hspace{-0,06cm}\times\hspace{-0,03cm} F(T)) \hspace{-0,06cm}\times\hspace{-0,03cm} F(U))}
	\ar[u, Rightarrow, "{\Id_{F(R)} \hspace{-0,06cm}\times\hspace{-0,03cm} \alpha^*_{F(S), F(T), F(U)}}" description]
	\ar[r, Rightarrow, "{(\Id_{F(R)} \hspace{-0,06cm}\times\hspace{-0,03cm} F_{\times}^{S, T}) \hspace{-0,06cm}\times\hspace{-0,03cm} \Id_{F(U)}}"]
	\ar[dr, phantom, "{\text{nat.}}" description]
&{F(R) \hspace{-0,06cm}\times\hspace{-0,03cm} (F(S \hspace{-0,06cm}\times\hspace{-0,03cm} T) \hspace{-0,06cm}\times\hspace{-0,03cm} F(U))}
	\ar[u, Rightarrow, "{\Id_{F(R)} \hspace{-0,06cm}\times\hspace{-0,03cm} F_{\times}^{S \times T, U}}" description]
	\ar[ddr, triple, "{\Omega_{\alpha}^{R, S \times T, U}}" description]
&{} \\
{}
&{(F(R) \hspace{-0,06cm}\times\hspace{-0,03cm} (F(S) \hspace{-0,06cm}\times\hspace{-0,03cm} F(T))) \hspace{-0,06cm}\times\hspace{-0,03cm} F(U)}
	\ar[u, Rightarrow, "{\alpha^*_{F(R), F(S) \times F(T), F(U)}}" description]
	\ar[r, Rightarrow, "{\Id_{F(R)} \hspace{-0,06cm}\times\hspace{-0,03cm} (F_{\times}^{S, T} \hspace{-0,06cm}\times\hspace{-0,03cm} \Id_{F(U)})}"]
	\ar[ddr, triple, "{\Omega_{\alpha}^{R, S, T} \hspace{-0,06cm}\times\hspace{-0,03cm} \Id_{F(U)}}" description]
&{(F(R) \hspace{-0,06cm}\times\hspace{-0,03cm} F(S \hspace{-0,06cm}\times\hspace{-0,03cm} T)) \hspace{-0,06cm}\times\hspace{-0,03cm} F(U)}
	\ar[u, Rightarrow, "{\alpha^*_{F(R), F(S \times T), F(U)}}" description]
	\ar[d, Rightarrow, "{F_{\times}^{R, S \times T} \hspace{-0,06cm}\times\hspace{-0,03cm} \Id_{F(U)}}" description]
&{} \\
{((F(R) \hspace{-0,06cm}\times\hspace{-0,03cm} F(S)) \hspace{-0,06cm}\times\hspace{-0,03cm} F(T)) \hspace{-0,06cm}\times\hspace{-0,03cm} F(U)}
	\ar[uuu, Rightarrow, "{\alpha^*_{F(R) \times F(S), F(T), F(U)}}"]
	\ar[ru, Rightarrow, "{\alpha^*_{F(R), F(S), F(T)} \hspace{-0,06cm}\times\hspace{-0,03cm} \Id_{F(U)}}" description]
	\ar[dr, Rightarrow, "{(F_{\times}^{R, S} \hspace{-0,06cm}\times\hspace{-0,03cm} \Id_{F(T)}) \hspace{-0,06cm}\times\hspace{-0,03cm} \Id_{F(U)}}"']
&{}
&{F(R \hspace{-0,06cm}\times\hspace{-0,03cm} (S \hspace{-0,06cm}\times\hspace{-0,03cm} T)) \hspace{-0,06cm}\times\hspace{-0,03cm} F(U)}
	\ar[r, Rightarrow, "{F_{\times}^{R \times (S \times T), U}}"']
&{F((R \hspace{-0,06cm}\times\hspace{-0,03cm} (S \hspace{-0,06cm}\times\hspace{-0,03cm} T)) \hspace{-0,06cm}\times\hspace{-0,03cm} U)}
	\ar[uuu, Rightarrow, "{F(\alpha^*_{R, S \times T, U})}"'] \\
{}
&{(F(R \hspace{-0,06cm}\times\hspace{-0,03cm} S) \hspace{-0,06cm}\times\hspace{-0,03cm} F(T)) \hspace{-0,06cm}\times\hspace{-0,03cm} F(U)}
	\ar[r, Rightarrow, "{F_{\times}^{R \times S, T} \hspace{-0,06cm}\times\hspace{-0,03cm} \Id_{F(U)}}"']
&{F((R \hspace{-0,06cm}\times\hspace{-0,03cm} S) \hspace{-0,06cm}\times\hspace{-0,03cm} T) \hspace{-0,06cm}\times\hspace{-0,03cm} F(U)}
	\ar[u, Rightarrow, "{F(\alpha^*_{R, S, T}) \hspace{-0,06cm}\times\hspace{-0,03cm} \Id_{F(U)}}"']
&{}
\end{tikzcd}
}
\]
\[
\adjustbox{scale=0.7,center}{%
= \begin{tikzcd}[row sep = huge, column sep = huge]
{F(R) \hspace{-0,06cm}\times\hspace{-0,03cm} (F(S) \hspace{-0,06cm}\times\hspace{-0,03cm} (F(T) \hspace{-0,06cm}\times\hspace{-0,03cm} F(U)))}
	\ar[rr, Rightarrow, "{\Id_{F(R)} \hspace{-0,06cm}\times\hspace{-0,03cm} (\Id_{F(S)} \hspace{-0,06cm}\times\hspace{-0,03cm} F_{\times}^{T, U})}"]
&{}
&{F(R) \hspace{-0,06cm}\times\hspace{-0,03cm} (F(S) \hspace{-0,06cm}\times\hspace{-0,03cm} F(T \hspace{-0,06cm}\times\hspace{-0,03cm} U))}
	\ar[r, Rightarrow, "{\Id_{F(R)} \hspace{-0,06cm}\times\hspace{-0,03cm} F_{\times}^{S, T \times U}}"]
	\ar[d, Rightarrow, "{\alpha^{-1}_{F(R), F(S), F(T \times U)}}" description]
	\ar[ddr, triple, "{\Omega_{\alpha^*}^{R, S, T \times U}}" description]
&{F(R) \hspace{-0,06cm}\times\hspace{-0,03cm} F(S \hspace{-0,06cm}\times\hspace{-0,03cm} (T \hspace{-0,06cm}\times\hspace{-0,03cm} U))}
	\ar[dd, Rightarrow, "{F_{\times}^{R, S\times (T\times U)}}"'] \\
{(F(R) \hspace{-0,06cm}\times\hspace{-0,03cm} F(S)) \hspace{-0,06cm}\times\hspace{-0,03cm} (F(T) \hspace{-0,06cm}\times\hspace{-0,03cm} F(U))}
	\ar[urr, phantom, "{\text{nat.}}"]
	\ar[r, equal]
	\ar[u, Rightarrow, "{\alpha^*_{F(R), F(S), F(T \times U)}}"]
	\ar[dr, Rightarrow, "{F_{\times}^{R, S} \hspace{-0,06cm}\times\hspace{-0,03cm} \Id_{F(T) \hspace{-0,06cm}\times\hspace{-0,03cm} F(U)}}"description]
&{(F(R) \hspace{-0,06cm}\times\hspace{-0,03cm} F(S)) \hspace{-0,06cm}\times\hspace{-0,03cm} (F(T) \hspace{-0,06cm}\times\hspace{-0,03cm} F(U))}
	\ar[d, phantom, "{\text{nat.}}"]
	\ar[r, Rightarrow, "{\Id_{F(R) \hspace{-0,06cm}\times\hspace{-0,03cm} F(S)} \hspace{-0,06cm}\times\hspace{-0,03cm} F_{\times}^{T, U}}"]
	\ar[dr, Rightarrow, "{F_{\times}^{R, S} \hspace{-0,06cm}\times\hspace{-0,03cm} F_{\times}^{T, U}}"{name=N, description}]
&{(F(R) \hspace{-0,06cm}\times\hspace{-0,03cm} F(S)) \hspace{-0,06cm}\times\hspace{-0,03cm} F(T \hspace{-0,06cm}\times\hspace{-0,03cm} U)}
	\ar[d, Rightarrow, "{F_{\times}^{R, S} \hspace{-0,06cm}\times\hspace{-0,03cm} \Id_{F(T \hspace{-0,06cm}\times\hspace{-0,03cm} U)}}" description]
	\ar[to=N, phantom, "{\text{nat.}}"]
&{} \\
{}
&{F(R \hspace{-0,06cm}\times\hspace{-0,03cm} S) \hspace{-0,06cm}\times\hspace{-0,03cm} (F(T) \hspace{-0,06cm}\times\hspace{-0,03cm} F(U))}
	\ar[r, Rightarrow, "{\Id_{F(R \hspace{-0,06cm}\times\hspace{-0,03cm} S)} \hspace{-0,06cm}\times\hspace{-0,03cm} F_{\times}^{T, U}}"]
	\ar[dddr, triple, "{\Omega_{\alpha^*}^{R \times S, T, U}}"description]
&{F(R \hspace{-0,06cm}\times\hspace{-0,03cm} S) \hspace{-0,06cm}\times\hspace{-0,03cm} F(T \hspace{-0,06cm}\times\hspace{-0,03cm} U)}
	\ar[d, Rightarrow, "{F_{\times}^{R \times S, T \times U}}" description]
&{F(R \hspace{-0,06cm}\times\hspace{-0,03cm} (S \hspace{-0,06cm}\times\hspace{-0,03cm} (T \hspace{-0,06cm}\times\hspace{-0,03cm} U)))} \\
{}
&{}
&{F((R \hspace{-0,06cm}\times\hspace{-0,03cm} S) \hspace{-0,06cm}\times\hspace{-0,03cm} (T \hspace{-0,06cm}\times\hspace{-0,03cm} U))}
	\ar[ru, Rightarrow, "{F(\alpha^*_{R, S, T \times U})}" description]
	\ar[dr, triple, "{F(\Pi^*)}" description]
&{F(R \hspace{-0,06cm}\times\hspace{-0,03cm} ((S \hspace{-0,06cm}\times\hspace{-0,03cm} T) \hspace{-0,06cm}\times\hspace{-0,03cm} U))}
	\ar[u, Rightarrow, "{F(\Id_R \hspace{-0,06cm}\times\hspace{-0,03cm} \alpha^*_{S, T, U})}"'] \\
{((F(R) \hspace{-0,06cm}\times\hspace{-0,03cm} F(S)) \hspace{-0,06cm}\times\hspace{-0,03cm} F(T)) \hspace{-0,06cm}\times\hspace{-0,03cm} F(U)}
	\ar[uur, phantom, "{\text{nat.}}"]
	\ar[uuu, Rightarrow, "{\alpha^*_{F(R) \times F(S), F(T), F(U)}}"]
	\ar[dr, Rightarrow, "{(F_{\times}^{R, S} \hspace{-0,06cm}\times\hspace{-0,03cm} \Id_{F(T)}) \hspace{-0,06cm}\times\hspace{-0,03cm} \Id_{F(U)}}"']
&{}
&{F(((R \hspace{-0,06cm}\times\hspace{-0,03cm} S) \hspace{-0,06cm}\times\hspace{-0,03cm} T) \hspace{-0,06cm}\times\hspace{-0,03cm} U)}
	\ar[r, Rightarrow, "{F(\alpha^*_{R, S, T} \hspace{-0,06cm}\times\hspace{-0,03cm} \Id_{U})}"']
	\ar[u, Rightarrow, "{F(\alpha^*_{R \times S, T, U})}" description]
	\ar[dr, phantom, "{\text{nat.}}"]
&{F((R \hspace{-0,06cm}\times\hspace{-0,03cm} (S \hspace{-0,06cm}\times\hspace{-0,03cm} T)) \hspace{-0,06cm}\times\hspace{-0,03cm} U)}
	\ar[u, Rightarrow, "{F(\alpha^*_{R, S \times T, U})}"'] \\
{}
&{(F(R \hspace{-0,06cm}\times\hspace{-0,03cm} S) \hspace{-0,06cm}\times\hspace{-0,03cm} F(T)) \hspace{-0,06cm}\times\hspace{-0,03cm} F(U)}
	\ar[r, Rightarrow, "{F_{\times}^{R \times S, T} \hspace{-0,06cm}\times\hspace{-0,03cm} \Id_{F(U)}}"']
	\ar[uuu, Rightarrow, "{\alpha^*_{F(R \times S), F(T), F(U)}}" description]
&{F((R \hspace{-0,06cm}\times\hspace{-0,03cm} S) \hspace{-0,06cm}\times\hspace{-0,03cm} T) \hspace{-0,06cm}\times\hspace{-0,03cm} F(U)}
	\ar[r, Rightarrow, "{F(\alpha^*_{R, S, T}) \hspace{-0,06cm}\times\hspace{-0,03cm} \Id_{F(U)}}"']
	\ar[u, Rightarrow, "{F_{\times}^{(R \times S) \hspace{-0,06cm}\times\hspace{-0,03cm} T, U}}" description]
&{F(R \hspace{-0,06cm}\times\hspace{-0,03cm} (S \hspace{-0,06cm}\times\hspace{-0,03cm} T)) \hspace{-0,06cm}\times\hspace{-0,03cm} F(U)}
	\ar[u, Rightarrow, "{F_{\times}^{R \times (S \hspace{-0,06cm}\times\hspace{-0,03cm} T), U}}"']
\end{tikzcd}
}
\]

(A2)
\[
\begin{tikzcd}[row sep = huge, column sep = huge]
{(F(R) \times F(u)) \times F(S)}
	\ar[d, Rightarrow, "{F_{\times} \times \Id_{F(S)}}"']
&{F(R) \times F(S)}
	\ar[l, Rightarrow, "{(\Id_{F(R)} \times F_{u}) \times \Id_{F(S)}}"']
	\ar[d, equal]
&{} \\
{F(R \times u) \times F(S)}
	\ar[r, Rightarrow, "{F(r^*) \times \Id_{F(S)}}"]
	\ar[d, Rightarrow, "{F_{\times}}"']
	\ar[ur, triple, "{\Omega_{r^*} \times \Id_{F(S)}}" description]
&{F(R) \times F(S)}
	\ar[r, Rightarrow, "{F_{\times}}"]
&{F(R \times S)} \\
{F((R \times u) \times S)}
	\ar[rr, Rightarrow, "{F(\alpha^*)}"']
	\ar[urr, Rightarrow, "{F(r^* \times \Id_{S})}"{description, name=T}]
	\ar[from=2-1, to=T, phantom, "{\text{nat.}}"]
&{}
&{F(R \times (u \times S))}
	\ar[u, Rightarrow, "{F(\Id_{R} \times l^*)}"']
	\ar[from=T, triple, "{F(\mu^*)}" description]
\end{tikzcd}
\]
\[
= \begin{tikzcd}[row sep = huge, column sep = huge]
{}
&{F(R) \times F(S)}
	\ar[dl, Rightarrow, "{(\Id_{F(R)} \times F_{u}) \times \Id_{F(S)}}"']
	\ar[dr, Rightarrow, "{\Id_{F(R)} \times (F_{u} \times \Id_{F(S))}}"]
	\ar[d, phantom, "{\text{nat.}}" description]
&{} \\
{(F(R) \times F(u)) \times F(S)}
	\ar[r, Rightarrow, "{\alpha^*_{F(R), F(u), F(S)}}"]
	\ar[d, Rightarrow, "{F_{\times} \times \Id_{F(S)}}"']
&{F(R) \times (F(u) \times F(S))}
	\ar[d, Rightarrow, "{\Id_{F(S)} \times F_{\times}}" description]
	\ar[ddl, triple, "{\Omega_{\alpha^*}}" description]
&{F(R) \times F(S)}
	\ar[l, Rightarrow, "{(\Id_{F(R)} \times F_{u}) \times \Id_{F(S)}}"']
	\ar[d, equal] \\
{F(R \times u) \times F(S)}
	\ar[d, Rightarrow, "{F_{\times}}"']
&{F(R) \times F(u \times S)}
	\ar[d, Rightarrow, "{F_{\times}}" description]
	\ar[r, Rightarrow, "{\Id_{F(R)} \times F(l^*)}"']
	\ar[dr, phantom, "{\text{nat.}}" description]
	\ar[ur, triple, "{\Id_{F(R)} \times \Omega_{l^*}}" description]
&{F(R) \times F(S)}
	\ar[d, Rightarrow, "{F_{\times}}"] \\
{F((R \times u) \times S)}
	\ar[r, Rightarrow, "{F(\alpha^*)}"']
&{F(R \times (u \times S))}
	\ar[r, Rightarrow, "{F(\Id_{R} \times l^*)}"']
&{F(R \times S)}
\end{tikzcd}
\]

\end{document}